%Paper: funct-an/9211010
%From: lschweit@sol.UVic.CA (Larry Schweitzer)
%Date: Fri, 27 Nov 92 15:04:34 PST

%%%%%%%%%%%%%%%%%%%%%%%%%%%%%%%%%%%%%%%%%%%%%%%%%%%%%
%
% This requires AMS Tex version 2.1, 92 pages.
%
% \hsize set to 6.7truein, but the index has different
% margins.
%
%%%%%%%%%%%%%%%%%%%%%%%%%%%%%%%%%%%%%%%%%%%%%%%%%%%%%
%
%
\input amstex
\documentstyle {amsppt}
\tolerance=10000
\magnification=1200
%
%
% Margins:
%
\hsize=6.7truein %widened (but index differs)
%
% Double spacing:
%
\addto\tenpoint{\normalbaselineskip=18pt\normalbaselines}
\addto\eightpoint{\normalbaselineskip=15pt\normalbaselines}
\def\R{\Bbb R}
\def\Z{\Bbb Z}
\def\N{\Bbb N}
\def\T{\Bbb T}
\def\C{\Bbb C}
\def\Q{\Bbb Q}
\def\s{\sigma}
\def\w{\omega}
\def\ga{\gamma}
\def\be{\beta}

\def\F{\Cal F}
\def\normm{ \pa \quad \pa_{m} }
\def\normdp{\pa f \pa_{d}}
\def\normi{\pa \quad \pa_{\infty}}
\def\normo{\pa \quad \pa_0}
\def\basisq{X_1, ... X_q}
\def\basisp{X_1, ... X_p}
\def\pa{\parallel}
\def\tense{\widehat \otimes_{\epsilon}}
\def\tensp{\widehat \otimes_{\pi}}
\def\vq{\vec q}
\def\vr{\vec r}
\def\schwr{\Cal S (\R)}
\def\schwcm{{C}^{\s} (M)}
\def\schwmh{\Cal S_{H}^{\s}(M)}
\def\schwmns{{\Cal S} (M)}

\def\schwm{{\Cal S}^{\s} (M)}
\def\schwsig{{\Cal S^{\s}_{1}}(G)}
\def\schwsi{{\Cal S^{\s}}(G)}
\def\schw{{\Cal S}(G)}
\def\schwL{\Cal S^{\w}(G)}
\def\so{\Cal S^{\s}_{1}(G)}
\def\soinf{\Cal S^{\s}_{\infty}(G)}
\def\solto{\Cal S^{\s}_{2}(G)}
\def\solp{\Cal S^{\s}_{r}(G)}
\def\cc{C_{c}^{\infty}(G)}
\def\soinfh{\Cal S^{\s}_{\infty}(H)}
\def\soh{\Cal S^{\s}_{1}(H)}

\def\GC{G\rtimes C_{0}(M)}
\def\GSno{G \rtimes \schwmns}
\def\schwgtm{\Cal S (G\times M)}

\def\sonec{\Cal S^{\s}(G, E)}
\def\soneco{\Cal S_{1}^{\s}(G, E)}
\def\sonas{G\rtimes^{\s} A}
\def\sona{G\rtimes A}
\def\sonac{\Cal S^{\s}(G, A)}
\def\sonaco{\Cal S^{\s}_{1}(G, A)}
\def\Einf{E^{\infty}}

\def\EF{E\hat \otimes F}
\def\EinF{\Einf \hat \otimes F}

\def\Cai{C^{\infty}_{\alpha}}

\leftheadtext{Larry B. Schweitzer}
\rightheadtext{Dense $m$-convex Fr\'echet Subalgebras}
\refstyle{A}
\document
\centerline{\bf Dense $m$-convex Fr\'echet Subalgebras of Operator}
\centerline{\bf Algebra Crossed Products by Lie Groups}
\vskip\baselineskip
\centerline{Larry B. Schweitzer}
\vskip\baselineskip
\heading  Abstract  \endheading
\par
Let $A$ be a dense Fr\'echet *-subalgebra of a C*-algebra
$B$. (We do not require Fr\'echet algebras to be $m$-convex.)
Let $G$ be a Lie group, not necessarily connected, which
acts on both $A$ and $B$ by *-automorphisms, and let $\s$ be a
submultiplicative
function from $G$ to the nonnegative real numbers.
If $\s$ and the action of $G$ on $A$
satisfy certain simple properties,  we define a
dense Fr\'echet *-subalgebra $G\rtimes^{\s} A$ of the crossed
product $L^{1}(G, B)$.  Our algebra consists of differentiable
$A$-valued functions on $G$, rapidly vanishing in $\s$.
\par
We give conditions on the action of $G$ on $A$ which imply
the $m$-convexity of the dense subalgebra $G\rtimes^{\s}A$.
A locally convex algebra is said to be $m$-convex if there is
a family of submultiplicative seminorms for the topology
of the algebra.  The property of $m$-convexity is  important  for
a Fr\'echet algebra, and is useful
in modern operator theory.
\par
If $G$ acts as a transformation group on a manifold $M$,
we develop a class of dense subalgebras for the crossed
product
$L^{1}(G, C_{0}(M))$,
where $C_{0}(M)$ denotes the continuous functions on
$M$ vanishing at infinity with the sup norm topology.
We define Schwartz functions $\Cal S(M)$
on $M$, which are differentiable with respect to some group action on
$M$, and are rapidly vanishing with respect to some scale on $M$.
We then form a dense $m$-convex Fr\'echet *-subalgebra
$G\rtimes^{\s}\Cal S (M)$ of rapidly vanishing, $G$-differentiable
functions from $G$ to $\Cal S(M)$.
\par
If the reciprocal of $\s$
is in $L^{p}(G)$ for some $p$, we prove that our group
algebras $\Cal S^{\s}(G)$
are nuclear Fr\'echet spaces, and that $G\rtimes^{\s}A$
is the projective completion $\Cal S^{\s}(G) {\widehat \otimes} A$.
\vskip\baselineskip
\vskip\baselineskip
\vskip\baselineskip
\head Contents \endhead
$$\aligned
\S 0.&\text{   Introduction} \\
\S 1.&\text{   Weights, Gauges and Group Algebras}\\
&\S 1.1.\text{   Scales, Weights, Gauges, Orderings and
Equivalence Relations}\\
&\S 1.2.\text{   Lie Groups and Schwartz Functions}\\
&\S 1.3.\text{   Group Schwartz Algebras}\\
&\S 1.4.\text{   Gauges that Bound $Ad$}\\
&\S 1.5.\text{   Polynomial Growth Groups}\\
&\S 1.6.\text{   Examples of Weights that Bound $Ad$}\\
\S 2.&\text{   Smooth Crossed Products}\\
&\S 2.1.\text{   Fr\'echet Space Valued Schwartz Functions}\\
&\S 2.2.\text{   Conditions for an Algebra}\\
\S 3.&\text{   $m$-convexity}\\
&\S 3.1.\text{   Conditions for $m$-convexity of the Crossed Product}\\
&\S 3.2.\text{   Non $m$-convex Group Algebras}\\
\S 4.&\text{   Conditions for a *-algebra}\\
\S 5.&\text{   Examples and Scaled $G$-spaces}\\
\S 6.&\text{   Sup Norm and $L^{r}$ Norm in Place of $L^{1}$ Norm,
and Nuclearity}\\
\text{Appendix.}&\text{   Sets of $C^{\infty}$-vectors}\\
&\text{   Index.}
\endaligned
$$
\vskip\baselineskip
\heading \S 0 Introduction  \endheading
Let  $G$ be a Lie group, not necessarily connected,
 and let $B$ be a C*-algebra
with strongly continuous action of $G$.
We let $A$ be a dense Fr\'echet *-subalgebra of $B$ (with
Fr\'echet topology possibly stronger than the topology on $B$),
such that $G$ acts on $A$ and the action
of $G$ on $A$ is infinitely
differentiable.
(Fr\'echet algebras are assumed locally convex, but not $m$-convex.)
We construct and study a class of
dense Fr\'echet *-subalgebras
$G\rtimes A$ of
the convolution algebra $L^{1}(G, B)$.
Many such (topological) subalgebras have been
defined in the literature for
various sorts of groups and algebras.
The dense subalgebras in \cite{ENN}, \cite{Jo},
\cite{Ji}, \cite{Vi \S 7}, the Schwartz algebras
in \cite{Bo, Thm 2.3.3}, and the Schwartz algebra of a connected
simply connected nilpotent Lie group \cite{Ho}, \cite{Lu 1},
\cite{Lu 2} are all special cases of the algebras we define here.
Our algebras consist of differentiable rapidly vanishing
functions from $G$ to $A$.
They have some advantages over algebras of compactly supported
differentiable functions.  For example,
the algebras of compactly supported functions
are rarely spectral invariant in the C*-crossed product $G\rtimes B$,
unless $G$ is compact or $G$ acts properly on
a manifold (see \cite{BC, Appendix 1}).
However, the Schwartz algebra of a connected simply connected nilpotent
Lie group is always spectral invariant in the group C*-algebra
$C^{*}(G)$ \cite{Lu 1, Prop 2.2} \cite{Sc 3, Cor 7.17}.
Spectral invariance is an important
property since it implies that the $K$-theories $K_{*}(G\rtimes A)$
and $K_{*}(G \rtimes B)$ are the same \cite{Co, VI.3}\cite{Sc 2, Lemma 1.2,
Cor 2.3}.
In the papers \cite{Sc 1}\cite{Sc 2}\cite{Sc 3}\cite{Sc 4} we
study the representation theory and spectral invariance of
the dense subalgebras we define here.
\par
We consider the following  situation.
We call a
 positive Borel measurable
function $\s\geq 1$ on $G$ a {\it weight which bounds $Ad$} if it satisfies
 the four properties
$$ \aligned  \s(gh) \leq &\s(g) \s(h) \\
\s(g^{-1}) & = \s(g) \\
 \s(e)& = 1\\
\pa Ad_{g}\pa &\leq D \s^{p}(g)\endaligned $$
for some constant $D$ and some positive integer $p$,
where $g, h\in G$,
and $\pa \quad \pa $  denotes
some norm
on the
space of linear operators on the Lie algebra of $G$.
The group algebra $\schwsi$ of $\s$-rapidly vanishing
differentiable functions on $G$
is then well formed, and is a dense Fr\'echet *-subalgebra
of $L^{1}(G)$. (This corresponds to the case $A=B= \C$.)
In \S 2, we define an appropriate notion of a tempered
action of $G$ on the dense subalgebra $A$,
much as in \cite{ENN} for the case $G=\R$,
and \cite{DuC 1, \S 4}\cite{DuC 2} for the case of a simply connected,
connected nilpotent Lie group.
If the action of $G$ on $A$ is tempered,
we define a dense Fr\'echet *-subalgebra
$\sonas$
 of
$L^{1}(G, B)$,
which consists of $\s$-rapidly vanishing $A$-valued differentiable
functions on $G$.
\par
As we noted, these dense subalgebras $\sonas$
 generalize those in
\cite{ENN},  and \cite{Bo, Thm 2.3.3}.
(The smooth crossed products by $\R$ in \cite{ENN, \S 7} and by
$\Z$ in \cite{Ns} may not in general be of the form $\sonas$, but
they are always a projective limit $\varprojlim
G\rtimes^{\s_{n}} A$, where $\s_{n}$ is an increasing sequence
of weights on $G$.)
In the purely group algebra case, our $L^{2}$ version of Schwartz
functions $\solto$ generalizes the spaces of rapidly
vanish functions on a discrete group with a
length function studied in \cite{Ji} \cite{Jo} \cite{Vi, \S 7}
(If $G$ is rapidly decaying \cite{Jo}, these $L^{2}$ Schwartz
functions become algebras.)
See  Remark 6.18 below for a comparison of our algebras with the
group algebra of Rader defined for a reductive Lie group \cite{War},
and  with the zero Schwartz space for $SL_{2}(\R)$
studied in \cite{Bar, \S 19}.
\par
We investigate whether our dense subalgebras are $m$-convex,
that is, whether their topology is given by a family of
submultiplicative seminorms.
The property of
$m$-convexity for Fr\'echet
(or more general locally convex) algebras
is important for several
reasons.  Arens was the first to notice that
$m$-convexity implies the existence and
continuity of entire functions on the algebra \cite{Ar}.
Although a unital $m$-convex algebra may not have an
open group of invertible elements, inversion is always
continuous on the set of invertible elements \cite{Mi, Prop 2.8}.
Every $m$-convex algebra is a projective limit of
Banach algebras \cite{Mi, Thm 5.1}.  Many other
important properties of $m$-convexity are studied in
\cite{Mi} and \cite {Ze}.
The property of $m$-convexity  often arises
in modern operator algebra theory  \cite{Ph},
\cite{Br}, \cite{Da, p. 136}.
For example, in  my thesis \cite{Sc 1} I study
a method for
showing that dense subalgebras are spectral
invariant in their  C*-algebras, which requires the
dense subalgebras to be $m$-convex.
\par
We define what it means for an action of
$G$ on $A$ to be {\it $m$-tempered}
(see \S 3), and show that $\sonas$ is $m$-convex if
the action is $m$-tempered.  As far as I know, no results on the
$m$-convexity of
a Fr\'echet crossed product like $\sonas$ have  appeared
in the literature.  For the case of the standard
Schwartz algebra of a simply connected nilpotent
Lie group, the $m$-convexity is noticed in \cite{Lu 2}.
The $m$-convexity of the group convolution algebra
of $C^{\infty}$-functions with compact support and
inductive limit topology is proved in \cite{Ma, Chap VIII, \S 11}.
In \S 3.2, we also give several interesting examples of {\it non}
$m$-convex group Schwartz algebras.
\par
In \S 4, we show that $\sonas$ has a well defined and
continuous $*$-operation.
In \S 5, we develop the example of a transformation
group, where $B=C_{0}(M)$, the commutative C*-algebra of continuous functions
on a locally compact $G$-space
 $M$ vanishing at infinity, with pointwise multiplication.
We let $\ga$
be a  scale (or nonnegative real-valued function)
 on $M$, and define $C^{\ga}(M)$ to be the Fr\'echet
algebra
of continuous functions on $M$ which vanish $\ga$-rapidly.
Then we take $A=\Cal S^{\ga}_{H}(M)$,
the set of $C^{\infty}$-vectors for the action of a Lie group $H$
on $C^{\ga}(M)$.
If certain simple conditions for the action of $G$ on $M$
are satisfied (see \S 5), we show that
the convolution algebra $G \rtimes^{\s}
\Cal S^{\ga}_{H}(M) $ is a well formed
$m$-convex Fr\'echet algebra, which is dense in $L^{1}(G, C_{0}(M))$.
\par
In \S 6,
we show that the  integrability condition
$$(\exists p\in \N)
\qquad \int_{G}{1\over{\s^{p}(G)}}dg <\infty \tag 0.1$$
implies that our group algebra $\schwsi$
is a nuclear Fr\'echet space, and that $\sonas$ is then isomorphic
to the projective completion $\schwsi {\widehat \otimes} A$.
  We show that condition (0.1) is
satisfied in several important cases,
and compare our algebras with those
in  \cite{Jo}.
\par
We end this introduction by giving a brief overview of some of the dense
subalgebras of crossed products we
are familiar with from the literature.
All of the dense subalgebras we speak of will
be complete locally convex algebras in some topology.
\par
First consider the group algebra case, when $A$ and $B$ are
both the complex numbers.
For a general Lie group, the $C^{\infty}$-functions with
compact support form a dense subalgebra (see \cite{Ma, Chap. VIII \S 11}).
If $G$ is compact, this subalgebra of course reduces to
the convolution algebra of $C^{\infty}$ functions on $G$.
 When $G$ is a simply connected
nilpotent
Lie group, there is the well known algebra of Schwartz functions
on $G$ obtained via the exponential map (see for example
\cite{Lu 1}, \cite{Ho, p. 346})
This generalizes to an algebra of Schwartz functions
for any closed subgroup of a nilpotent Lie group.
If $G$ is a reductive Lie group, Harish-Chandra developed
a space of rapidly vanishing differentiable functions
which is an
algebra (see Rader's theorem \cite{War 1, Prop. 8.3.7.14}) - see
also \cite{Bar}.
For rapidly vanishing functions on a discrete group,
see the definitions in \cite{Ji}, \cite{Jo} and \cite{Vi, \S 7}.
In the appendix of \cite{Jo}, a space of rapidly vanishing
functions for a general locally compact group with a
(subadditive) length function is also studied.  This algebra does not take
into account any $C^{\infty}$-structure of the group.
\par
In the special case of a transformation group, namely when the
algebra $B$  is continuous functions vanishing at infinity
on a manifold $M$, many definitions of a dense
subalgebra have also been made.  We could take the
convolution algebra of compactly supported continuous
functions on $G\times M$ (see \cite{EH, \S 3}).  This
ignores any $C^{\infty}$ structure on $G$ or $M$,
so we could instead take the $C^{\infty}$ functions
with compact support on $G\times M$.
If $G$ is discrete, such an algebra is introduced
in the appendix of \cite{BC}.  Perhaps the most recent
contribution in the case of a transformation group
is the work of F. du Cloux \cite{DuC 4} for
the case where $G$ is an algebraic Lie group and $M$ an
algebraic variety.
\par
Next we consider the general case of dense subalgebras for a
crossed product, where $B$ can be any C*-algebra
and $A$ is any dense Fr\'echet subalgebra of $B$.
The $C^{\infty}$-functions
with compact support from $G$ to $A$ can be used as
a dense subalgebra.
In \cite{ENN} and \cite{Ns} an appropriate subalgebra
of $A$-valued Schwartz functions is defined for crossed products
of $\R$ and $\Z$ respectively with an arbitrary dense Fr\'echet
subalgebra $A$ on which $G$ acts smoothly in an appropriate
sense.
If $G$ is elementary Abelian, a
convolution algebra of Schwartz functions is studied in
\cite{Bo, see \S 2.1.4}, in the case that $A$ is the set
of $C^{\infty}$ vectors for the action of $G$ on $B$.
Also in this work of Bost, a dense subalgebra
of exponentially vanishing (but not differentiable) functions is studied
for
an arbitrary locally compact group
$G$ possessing a subadditive length function -
see \cite{Bo, \S 2.3}.
\par
Throughout this paper, the notations $\N$, $\N^{+}$, $\Z$, $\Q$, $\R$, $\T$
shall be used for the natural numbers with zero, natural numbers
without zero, integers, rationals, reals,
and the circle group respectively.  All of our algebras
will be over $\C$.  The term norm may be used
interchangeably with the term seminorm.
If the positive definiteness
of a norm is  important, we shall state it explicitly.  We shall
use the notation $C_{0}(M)$ ($C_{0}^{\infty}(M)$) to mean
continuous (differentiable) functions
vanishing (along with all their derivatives)
at $\infty$ on a locally compact space (Lie group)
$M$ with sup norm (and sup norm of derivatives).
Compactly supported
continuous (differentiable) functions on a locally compact
space (differentiable manifold) $M$ shall be denote by $C_{c}(M)$
($C_{c}^{\infty}(M)$).
The term differentiable will always mean infinitely
differentiable.
We shall often refer to the standard set of
Schwartz functions $\Cal S(\R)$
on $\R$.
This means the set of differentiable
functions on $\R$ that vanish
at infinity (along with their derivatives)
faster than the reciprocal of any polynomial.
\par
All groups
will be assumed Hausdorff, and,  unless otherwise stated,
locally compact.
We work with Fr\'echet spaces.  However,
some variant of what is said ought to work if the dense subalgebra
$A$ were a general complete locally convex algebra.
\par
This paper is a portion of my Ph.D. thesis, written at the University
of California at Berkeley, under the supervision of Marc A. Rieffel.
I would like to thank Bill Arveson, Joeseph Wolf,
and especially Marc Rieffel for many helpful
comments and suggestions.
\heading \S 1 Weights, Gauges and Group Algebras \endheading
\par
In this section, we define what it means for a positive real valued
function $\s$ on a Lie group $G$ to be  a weight or
a gauge
 and state an appropriate
additional condition
on $\s$, namely that $\s$ bound $Ad$,  so that
the Fr\'echet space $\schwsi$ of differentiable
rapidly vanishing functions on $G$ will be  a *-subalgebra
of the convolution
algebra $L^{1}(G)$.
We give general conditions on $G$ which imply the
existence of weights or gauges that bound $Ad$.
\heading \S 1.1 Scales, Weights, Gauges, Orderings and Equivalence Relations
\endheading
\par
We define scales, weights, and gauges, and introduce orderings and
equivalence relations on them.  We introduce the notion of a compactly
generated group, and show that such groups are characterized by
having a largest gauge (which we call the word gauge)
and a largest weight defined on them.  The word gauge will be the analog of
the length function on a discrete finitely generated group.
\subhead Definition 1.1.1  \endsubhead
Let $M$ be any topological space.
Let $\s$
be a real valued Borel measurable
function which
maps $M$ to the half open interval $[0, \infty ) $.
We call $\s$ a {\it scale} on the space $M$.
We define some special types of scales.
Let $G$ be any (Hausdorff) topological group, and let $e$ denote the
identity element of $G$.
We say that a scale $\tau$ on $G$ is a {\it gauge}  if
$$\tau(gh) \leq \tau(g)  + \tau(h) \tag 1.1.2 $$
$$\tau(g^{-1}) =\tau(g) \tag 1.1.3  $$
$$\tau(e) = 0 \tag 1.1.4 $$
for $g, h \in G$.
An example of a gauge on the real numbers is given by
$\tau(r) = |r|$.  The trivial scale $\tau\equiv 0$ defines
a gauge on any group.
In \cite{Jo} and \cite{Ji}, gauges are also called
{\it length functions}, but this differs from the common
usage of the term length function in the works
\cite{Pr} \cite{Ch} \cite{Ha} \cite{Ly}.
In \cite{Pr}, what we call a gauge is called a
(real valued) semigauge.
\par
We say that a scale $\w$ is a weight \cite{Py} \cite{Dz 1} on $G$
if $\w \geq 1$ and
$$\w(gh) \leq \w(g)  \w(h) \tag 1.1.5 $$
$$\w(g^{-1}) =\w(g) \tag 1.1.6  $$
$$\w(e) = 1 \tag 1.1.7  $$
for $g, h \in G$.
By  \cite{Dz 1, Prop 2.1} or Theorem 1.2.11 below, if $G$ is locally compact,
then weights are bounded on compact sets.
Since any gauge $\tau$ defines a weight on $G$ by $\w(g)
=e^{\tau (g)}$, gauges are also bounded on compact sets.
Two examples of weights on $\R$ are given by
$\w(r) = e^{|r|}$ and $\w(r) =1 +|r|$. The trivial scale
$\w\equiv 1$ defines
a weight on any group.  The correspondences
$\tau \mapsto e^{\tau}$ and
$\w \mapsto \log(\w)$ give a bijection
between gauges and weights on a group $G$.
We call this bijection the {\it exponential bijection } between gauges
and weights.
\subheading{Definition 1.1.8}
We say that a scale  $\s_{2}$ {\it dominates }
another scale $\s_{1}$ ($\s_{2}
\succcurlyeq
\s_{1}$)  if
there are constants $C, D\geq 0$ and $m\in \N$ such that
$$ \s_{1}(g) \leq C\s_{2}^{m}(g) + D, \qquad g \in G.
\tag 1.1.9 $$
We say that
$\s_{1}$ and $\s_{2}$ are {\it equivalent} ($\s_{1} \thicksim \s_{2}$) if
$\s_{1} \preccurlyeq \s_{2}$ and $\s_{2} \preccurlyeq \s_{1}$.
Since $\preccurlyeq$ is transitive,
we have a partial order on equivalence classes
of scales.
 We say that
$\s_{2}$ {\it strictly dominates} $\s_{1}$ if $\s_{1}\preccurlyeq \s_{2}$
but $\s_{1}$ and $\s_{2}$ are not equivalent.
\par
On a compact group, since weights and gauges are bounded,
every weight and every gauge is equivalent to the constant scale $\s \equiv
1$.
\par
The weight $\w(r) = 1+|r|$ is equivalent to the gauge $\tau(r)
= |r|$ on $\R$.  In fact, every gauge $\tau$ on an
arbitrary topological group
is equivalent to the (subadditive)
weight $1 + \tau$. Conversely, if $\w$ is a
subadditive weight, we can define
an equivalent gauge $\tau$ by
$$\tau(g) = \cases 0 & g=e \\ \w(g) & g\not= e.\endcases
$$
\subheading{Example 1.1.10}
In general, if $\w$ is a weight which is not subadditive,
there is no way of defining a gauge which
is equivalent to $\w$.  For example, let $G= \Z$ and
define $\w(n) = e^{|n|}$.  Then
if $\tau$ is a gauge on $\Z$, the inequality (1.1.9)
cannot be satisfied for $\s_{1}= \w$ and $\s_{2}=\tau$,
since for positive $n$, $\w(n) = e^{|n|}$ and
$\tau(n) \leq n\tau(1) $.   So $\w$ cannot be equivalent to
$\tau$.
We thus  have the proper inclusion
$$ \{\text{gauges}\} = \{\text{subadditive weights}\}
\subsetneq \{ \text{weights} \} $$
of equivalence classes. In fact, since $\w$ dominates
any linear function on $\Z$, it strictly dominates every gauge on
$\Z$.
\subheading{Example 1.1.11}
On groups with more than two generators, it is easy to come up with
weights which are not equivalent to any gauge, but also do not dominate
every gauge.  For example, let $G= \Z^{2}$, and define $\w(n_{1},
n_{2}) = e^{|n_{1}|}(1+\log(1+|n_{2}|))$.  Then by Example 1.1.10,
$\w$ restricted to the first copy of $\Z$ strictly dominates
every gauge on $\Z$, so $\w$ (as a weight on $G$) cannot be equivalent to
any gauge.  However, the weight $1+\log(1 + |n_{2}|)$
on the second copy of $\Z$
is strictly dominated
by the gauge $|n_{2}|$, so $\w$ cannot dominate every gauge on $G$.
This example shows that the partial ordering
on equivalence classes of weights may not
be total, since we may have $\w_{1} \npreceq \w_{2}$
and $\w_{2} \npreceq \w_{1}$.
\subheading{Definition 1.1.12}
We say that a gauge $\tau_{2}$ {\it strongly dominates}
a gauge $\tau_{1}$
($\tau_{2}
\succcurlyeq_{s} \tau_{1}$) if
there are constants $C, D\geq 0$ such that
$$ \tau_{1}(g) \leq C\tau_{2}(g) + D, \qquad g \in G.
\tag 1.1.13 $$
We call $\succcurlyeq$ the {\it usual ordering} on gauges,
in order to distinguish it from the strong ordering $\succcurlyeq_{s}$.
We say that $\tau_{1}$ and $\tau_{2}$ are {\it strongly equivalent}
($\tau_{1} \thicksim_{s} \tau_{2}$) if
$\tau_{1} \preccurlyeq_{s} \tau_{2}$
and $\tau_{2} \preccurlyeq_{s} \tau_{1}$.
\par
The natural exponential bijection
between  gauges and
weights given by
$\tau \mapsto e^{\tau}$ and $\w \mapsto \log\w$ is order preserving,
when the strong ordering is placed on
gauges and the usual ordering on weights.
On any compact group, every gauge is bounded by a constant and hence is
strongly equivalent to the trivial
gauge $\tau \equiv 1$.
\subheading{Example 1.1.14}
We give an example of two gauges which are equivalent but not strongly
equivalent.
Let $G= \Z$ and let $\tau_{1}(n) = |n|$, $\tau_{2}(n) = |n|^{1/2}$.
We have
$\tau_{2} \leq \tau_{1}$ and  $\tau_{1} \leq \tau_{2}^{2}$,  so these
two gauges are equivalent.  However $\tau_{2}$ cannot strongly dominate
its square $\tau_{1}$, so they are not strongly equivalent.
%
\subheading{Example 1.1.15} We give an
example showing that for some groups $G$ there is no largest weight,
and no largest gauge (in either ordering).
Let $G=F_{\infty}$ be the free group on
countably infinitely many generators.
Because the exponential bijection between gauges and weights
is order preserving, it suffices to show
that
for any gauge $\tau$ there is another gauge which strongly
dominates $\tau$, and which is not
equivalent to $\tau$ (in the usual ordering).
Let $c_{n}$ be the value
of $\tau$ on the $n$th generator $g_{n}$, where $n\in \N^{+}$.
Define another gauge
$\ga$ by letting
$\ga(g_{n_{1}}^{\alpha_{1}}\dots g_{n_{k}}^{\alpha_{k}})$ be
the sum $\sum_{i=1}^{k} |\alpha_{i}|e^{c_{n_{i}}+i}$.
Since $\tau(g_{n_{1}}^{\alpha_{1}}\dots g_{n_{k}}^{\alpha_{k}})\leq
\sum_{i=1}^{k} { |\alpha_{i}|c_{n_{i}}}$, we have $\tau \leq \ga$ and
so $\tau \preccurlyeq_{s} \ga$. On the other hand, we cannot have $\ga
\preccurlyeq_{s} \tau$ or even
$\ga \preccurlyeq \tau$, since then we would have
$d \in \N$ and $C, D>0$ such that
$$ \ga(g_{n}) = e^{c_{n} +n} \leq Cc_{n}^{d} +D, \qquad\,\, n \in \N^{+}.
\tag 1.1.16  $$
But $e^{c_{n}}$ by itself could only be bounded by
$Cc_{n}^{d} +D$ if $c_{n}$ is a bounded sequence.  And if $c_{n}$
is bounded, then the right hand side of (1.1.16) is bounded, so the
inequality cannot hold.  So $\ga$ stricly dominates $\tau$.  Hence there
is no  largest gauge on $G$, in either of the orderings
on gauges, and consequently no largest weight on $G$.
\subheading{Example 1.1.17. Word Gauge}
Let $G$ be a topological
group.  Let $U$ be any subset
containing the identity of $G$ (not necessarily
relatively compact), which satisfies $U=U^{-1}$
(we can always replace $U$
with $U\cup U^{-1}$ to achieve this) and such that
$$ \cup_{n=0}^{\infty} U^{n} = G. \tag 1.1.18$$
We call
such a set $U$ a {\it generating set} for $G$.
We define a gauge $\tau_{U}$ on $G$ by
$$ \tau_{U}(g) = \min \{ n \, | \, g \in U^{n} \}, \tag 1.1.19 $$
where $U^{0}=\{e\}$ (see \cite{Py}).
If the generating set $U$ is understood, we simply write $\tau$.
Note that the set $U$ could be the entire group, in which
case $\tau\equiv 0$.
\par
If $G$ is a discrete group and $U$ is a set of generators together
with their inverses, then $\tau$ is simply the word length. (This
follows straight from the definition (1.1.19).) So in general, $\tau$ can
be regarded as a kind of generalized word length.
We call
$\tau$ the {\it  word gauge with respect to the
generating set $U$}.  We may simply refer to $\tau$ as the
word gauge, if $U$ is understood.
\par
We call the weight $\w=e^{\tau}$ the
{\it exponentiated word weight} with respect to $U$.
\subheading{Definition 1.1.20. Compactly Generated Groups}
We say that a topological group
$G$ is {\it compactly generated } if there is some
relatively compact (rel. comp.)
generating set $U$ for $G$\cite{HR, Def 5.12}.
If $G$ is locally compact and compactly generated, then $G$ has
an open rel. comp. generating set.
\par
Let $\tau$ be the word gauge with respect to an open rel. comp.
generating set $U$.  The set of $g$ with $\tau (g) \leq n$ is the
open set $U^{n}$, so the
function $\tau $ is measurable (in fact upper-semicontinuous) on $G$.
In general, $\tau$ is not continuous, even if $U$ is open.
For example take $G=\R$,
$U= (-1, 1)$, or let $G$ be any connected non-compact group.
\par
Examples of compactly generated groups are given by discrete finitely
generated groups and connected locally compact groups.  The free
group on infinitely many generators is not compactly generated in the
discrete topology. (Note that
therefore the class of compactly generated groups
is not closed under subgroups.)
In Example 1.1.15,  we showed
that there exists no largest equivalence class
of weight or gauge on the free
(non-compactly generated) group $F_{\infty}$.
More generally, we have the following theorem.
\proclaim{Theorem 1.1.21} Let $G$ be a locally compact group.
If $G$ is  compactly generated, the strong
equivalence class of the
word gauge
does not depend on the choice of rel. comp. generating
set $U$.
The group $G$ is compactly generated
iff $G$ is a countable union of compact sets
and there exists a gauge (namely the word gauge)
on $G$ which strongly dominates every other
gauge on $G$.   A compactly generated
group $G$ is compact iff every gauge on $G$ is bounded.
Similar statements are true for weights.
\par
Also,  a quotient of a compactly generated group by a
closed normal subgroup is compactly generated.
If the connected component $G_{0}$ of the identity of $G$ is open,
then $G$ is compactly generated if and only if the discrete group $G/G_{0}$
is finitely generated.
\endproclaim
\par
We remark that if $G_{0}$ is not open in the last paragraph of the theorem,
the
forward implication is false.  For example, let $G$ be the product
of countably many copies of the compact group $\Z/2\Z$.  Then
$G$ is compact by Tychonoff's theorem, $G_{0}=(0,0,\dots)$ is not
open, and $G=G/G_{0}$ is not finitely generated.
\demo{Proof of Theorem 1.1.21}
Assume that $G$ is compactly generated.
Given any open rel. comp. generating sets
$U, V$ (such sets exist since $G$ is locally compact),
there exist positive integers
$k$ and $l$ such that $U \subseteq V^{k}$ and
$V\subseteq U^{l}$.
If $g\in U^{n}$, then $g\in V^{kn}$.  It follows that $\tau_{V}(g)
\leq k\tau_{U}(g)$.
Thus the equivalence class
of the gauge does not depend on the choice of
generating set.
\par
We show that $\tau$ strongly dominates every gauge on $G$.
Let $\ga$ be another gauge.  Let $C= \sup_{g \in U} \ga(g)$.
(Here we have
used $G$ locally compact and Theorem 1.2.11 below or \cite{Dz 1, Prop 2.1}
to see that $C<\infty$.)
Then if $\tau_{U}(g) = n$, let $g= g_{1} \dots g_{n}$ with
$g_{i} \in U$.  We have
$$\ga(g) \leq \ga(g_{1}) + \dots \ga(g_{n} ) \leq n C = C \tau(g).$$
Hence $\tau $ strongly dominates $\ga$.
\par
It follows immediately that the exponentiated word weight $\w$ dominates
every weight on $G$.  Also, since $G = \cup_{k=0}^{\infty} U^{k}$, $G$ is
a countable union of compact subsets.
\par
If $G$ is compact, clearly every gauge
and every weight on $G$ is bounded, since
the word gauge is bounded.
Conversely, assume that $G$ is compactly generated and that
either every gauge on $G$ is bounded, or that
every weight on $G$ is bounded.
Then the word gauge $\tau$ is a bounded
function on $G$, so $G= U^{N}$
for some $N\in \N$.  Hence $G$ is compact, since
$U$ is relatively compact. (The converse did not use $G$ locally compact.)
\par
We show that if $G$ is not compactly generated and is
a countable union of compact sets, then $G$ has no
gauge that dominates every other gauge.  Let $\tau$ be a gauge on $G$.
To facilitate our proof, we use the following lemma.
\proclaim{Lemma 1.1.22} Let $\tau$ be a gauge on a topological group
$G$. Then $\tau$ is strongly equivalent to a gauge $\tau^{\prime}$ on $G$
such that $\tau^{\prime}$ is integer valued and $\tau^{\prime}(g)=0$
implies $g = e$. \endproclaim
\demo{proof} Define $\tau_{int}(g) $ to be the least integer greater than or
equal to $\tau(g)$.  Then $\tau \leq \tau_{int}$ and $\tau_{int} \leq \tau
+1$ so $\tau \thicksim_{s} \tau_{int}$.  Define
$$\tau^{\prime}(g) = \cases 1 & \tau(g) = 0, \, g \not= e \\
\tau_{int}(g) & \text{otherwise} \endcases
$$
Then $\tau^{\prime}$ is a
gauge since $\tau^{\prime}(gh) \leq \tau^{\prime}(g)
+\tau^{\prime}(h)$ still holds if $\tau(gh)=0$ and $gh\not= e$ (since then
$\tau^{\prime}(g)$ or $\tau^{\prime}(h)$ is greater than or equal to one).
Clearly $\tau \thicksim_{s} \tau^{\prime}$.  \qed \enddemo
So replace $\tau $ with the gauge $\tau^{\prime}$ in the lemma.
We construct a gauge $\ga$ on $G$ such that $\tau \leq \ga$ and
$\ga \npreceq \tau $.
Define
$$U_{m} = \{ \, g \, | \quad \tau(g) \leq m \quad \}, \tag 1.1.23$$
and let $V_{m}$ be increasing compact sets containing $e$ such that
$\cup_{m=0}^{\infty} V_{m} = G$.
Define increasing rel. comp. sets $W_{m}= U_{m}\cap V_{m}$.
Then $\cup_{m=0}^{\infty} W_{m} = G$,
and  since $G$ is not compactly generated, we have $\cup_{k=0}^{\infty}
W^{k}_{m} \not= G$ for each $m$.
Inductively choose sequences $m_{p}\in \N^{+}$
and $g_{p}\in G$ such that $m_{0}= 1$, $g_{0}=e$ and
$$  g_{p} \in W_{m_{p}} -
\biggl( \cup_{k=0}^{\infty} W_{m_{p-1}}^{k} \biggr)
\tag 1.1.24 $$
for $p\geq 1$.  Note that $1=m_{0}<m_{1}<\dots $.
Define $\ga$ by
$$  \ga(g) = \inf_{h_{1}\dots h_{k}=g,\, h_{i}\in W_{m_{p_{i}}}-
W_{m_{p_{i}-1}}}   \biggl( \sum_{i=1}^{k} e^{m_{p_{i}}} \tau(h_{i})\biggr),
\tag 1.1.25$$
where the $p_{i}$'s are the least integers such that
$h_{i} \in W_{m_{p_{i}}}$.
Since the $m_{p}$'s tend to infinity as $p$ goes to infinity, every $g\in G$
lies in some $W_{m_{p}}$.  Hence $\ga(g) \leq e^{m_{p}}\tau(g)\leq
e^{m_{p}} m_{p}$, so $\ga $ is a well defined function from $G$ to
$[0, \infty)$.
\par
We show that $\ga$ is gauge.   Let $g, h \in G$.  Let $h_{1} \dots h_{k}=g$
and ${\tilde h_{1}} \dots {\tilde h_{\tilde k}} = h$ be such that
$$   \sum_{i=1}^{k} e^{m_{p_{i}}} \tau(h_{i}) \leq \ga(g) +\epsilon
\tag 1.1.26 $$
and
$$   \sum_{i=1}^{\tilde k} e^{m_{{\tilde p_{i}}}} \tau({\tilde h_{i}})
\leq \ga(h) +\epsilon. \tag 1.1.27 $$
Then
$$ \ga(gh) \leq \ga(g) +\ga(h) +2\epsilon \tag 1.1.28 $$
by (1.1.26),(1.1.27) and the definition of $\ga(gh)$ (1.1.25).
Letting $\epsilon$ tend to zero, we see that $\ga$ is  subadditive.
The remaining properties of a gauge are
not too difficult to show.
\par
We show that $\tau \leq \ga$.    Let $g \in G$.  Let $h_{1} \dots h_{k}=g$
 be such that
(1.1.26) is satisfied.  Then
$$\tau(g) \leq \sum_{i=1}^{k} \tau(h_{i}) \leq \sum_{i=1}^{k}
e^{m_{p_{i}}}\tau(h_{i}) \leq \ga(g) + \epsilon. $$
We let $\epsilon$ tend to zero to obtain $\tau \leq \ga$.
\par
Finally, we show that $\tau$ does not dominate $\ga$.  Recall $\tau(g_{p})
\leq m_{p}$ since $g_{p} \in W_{m_{p}}$.  It suffices to show that
$\ga(g_{p}) \geq e^{m_{p}}$.  But if we write $g_{p}= h_{1}\dots h_{k}$,
then by (1.1.24) at least one $h_{i}$ is {\it not} in $W_{m_{p-1}}$.
Then $h_{i} \not= e$ so $\tau(h_{i})\geq 1$.
Thus
$$\ga(g_{p}) \geq e^{m_{p}}1.$$
This proves  the first paragraph of Theorem 1.21.
\par
Let $G$ be compactly generated with generating
set $U$, and let $H$ be any closed normal subgroup.
Let $\pi \colon G \rightarrow G/H$ be the canonical map.
Then $\pi(U)$ is a relatively compact neighborhood of the identity
in $G/H$ and $\pi(U^{n})= \pi(U)^{n}$, so $\pi(U)$ generates $G/H$.
Note that this implies that $G/G_{0}$ is a compactly generated
group, which has the discrete topology if $G_{0}$ is open,
and so is finitely generated in that case.
Conversely, assume that $G/G_{0}$ is finitely
generated.  Let $g_{1}=e,g_{2}, \dots g_{k} \in G$ be such that
$\pi(g_{i})$ generate $G/G_{0}$.  Let $V$ be any relatively
compact neighborhood of the identity of $G_{0}$. (Here we use
the local compactness of $G$.)  Let $U = \cup_{i=1}^{k} g_{i} V$.
Since $G_{0}$ is connected, $\cup_{n=0}^{\infty} U^{n} \supseteq
\cup_{n=0}^{\infty} V^{n}
= G_{0}$\cite{HR, Thm 5.7}.  Also, $\cup_{n=0}^{\infty} U^{n}$
is a subgroup of $G$ which contains an
element of every coset of $G/G_{0}$.
Hence $\cup_{n=0}^{\infty} U^{n} = G$ and
$G$ is compactly generated.
\qed
\enddemo
If $G$ is compactly generated, we will say that the {\it  word gauge}
on $G$ is the  word gauge defined above with a relatively
compact generating set.
\vskip\baselineskip
\heading \S 1.2 Lie Groups and Schwartz functions \endheading
\par
If $G$ is a Lie group and $\s$ is a scale on $G$, we define
the $\s$-rapidly vanishing Schwartz functions $\schwsi$ on $G$.  We find
an appropriate condition (namely if $\s$ is equivalent to
every one of its left translates) on $\s$ so that $\schwsi$ is a Fr\'echet
space on which $G$ acts differentiably by left translation.
\par
Assume that $G$ is a real Lie group, not necessarily connected.
That is, assume that the connected component $G_{0}$
of the identity of $G$ is open and is a real connected Lie group.
Our Lie group $G$ is locally compact, so left Haar measure on $G$ exists.
Let $\frak G$ be the Lie algebra of $G$, and let $q$ be
the dimension of $\frak G$.
Let $\basisq$ be a basis for $\frak G$, and
let $G$ act strongly continuously on a Fr\'echet space $E$, with action
denoted by $\be$.
 If
$\ga \in \N^{q}$, we define the differential operator $X^{\ga}$
by
$$ X^{\ga}e = \biggl( {d\over {dt_{1}}}\biggr)^{\ga_1}
\dots \biggl( {d\over {dt_{q}}}\biggr)^{\ga_q}
\be_{exp(t_{1}X_{1}) \dots exp(t_{q}X_{q})}(e)
\restriction_{t_{1}=\dots t_{q}=0},
\tag 1.2.1
$$
where $e$ is a $C^{\infty}$-vector for the action of $G$ on $E$.
If $\pa \quad \pa_{m}$ is a family of seminorms topologizing $E$,
then we topologize the set $E^{\infty}$ of $C^{\infty}$-vectors
for the action of $G$ on $E$ by $\pa e \pa_{d, \ga}= \pa X^{\ga}e\pa_{d}$.
In this topology, $\Einf$ is a Fr\'echet space, and is
a dense $G$-invariant
subspace of $E$ (see the appendix, Theorem A.2).
\subheading{Definition 1.2.2}  Let $\s$ be a scale on $G$.
Let $\schwsig$ be the set of differentiable functions $\varphi \colon
G \rightarrow \C$ satisfying
$$ \pa \varphi \pa_{m,\ga} =
\pa \s^{m}X^{\ga}\varphi \pa_{1} =
\int_{G} \s^{m}(g)| X^{\ga} \varphi(g) | dg <\infty \tag 1.2.3 $$
for each $\ga \in \N^{q}$ and $m\in \N$.  Here $X^{\ga}\psi$
is defined by formula (1.2.1) with
$\be_{g}(\varphi)(h) = \varphi(g^{-1}h)$.
We call
$\schwsig$ the {\it $\s$-rapidly vanishing ($L^{1}$)  Schwartz
functions on $G$}.
Since we shall always use the $L^{1}$-norm
until \S 6, we shall simply write $\schwsi$ for $\schwsig$.
We topologize
$\schwsi$ by the seminorms (1.2.3).
It is easily checked that if $\s_{1} \thicksim \s_{2}$, then
$\Cal S^{\s_{1}}(G)$ is isomorphic to $\Cal S^{\s_{2}}(G)$.
\subheading{Question 1.2.4} If
$\Cal S^{\s_{1}}(G)$ is isomorphic to $\Cal S^{\s_{2}}(G)$,
then is  $\s_{1} \thicksim \s_{2}$ ?
\vskip\baselineskip
\par
We find a condition on $\s$ for which the action of $G$ by left translation
on $\schwsi$ will be  a well defined automorphism.
We use the notation $\varphi_{g}(h) = \varphi(g^{-1}h)$.  Then
$$\pa \varphi_{g} \pa_{d, \ga}^{\s_{g}} = \int_{G} \s^{d}_{g}(h) \, |\,
(X^{\ga}\varphi_{g})(h) \, | \, dh  $$
By the chain rule (see (2.2.3)), we may write
$$ (X^{\ga}\varphi_{g})(h) = \sum_{\be\leq \ga}
p_{\be}((Ad_{g^{-1}})_{ij})  (X^{\be}\varphi)_{g} (h), \tag 1.2.5$$
where $(Ad_{g^{-1}})_{ij}$ is the $ij$th matrix entry of $Ad_{g^{-1}}$,
and $p_{\be}$ is some polynomial.  Hence
$$  \aligned  \pa \varphi_{g} \pa_{d, \ga}^{\s_{g}} & \leq C_{g}
\sum_{\be \leq \ga}\int_{G} \s_{g}^{d} (h) \, |(X^{\be}\varphi)_{g}
(h) | dh \\
& =
C_{g}\sum_{\be\leq \ga}
 \int_{G} \s^{d} (h) \, |(X^{\be}\varphi)
(h) | dh \\
& = C_{g}
\sum_{\be\leq \ga}
\pa \varphi \pa_{d, \be}^{\s}.  \endaligned \tag 1.2.6$$
So for fixed $g$, the map $\varphi \mapsto \varphi_{g}$ is a
continuous linear map $\Cal S^{\s}(G) \rightarrow \Cal S^{\s_{g}}(G)$.
If $\s_{g}$ is equivalent to $\s$,
then  $\schwsi \cong \Cal S^{\s_{g}}(G)$,
so $\varphi \mapsto
\varphi_{g}$ will be a continuous isomorphism of the locally convex space
$\Cal S^{\s}(G)$, with inverse $\varphi \mapsto \varphi_{g^{-1}}$.
\subheading{Definition 1.2.7}
We say that a scale $\s$ is {\it translationally
equivalent} if $\s_{g} \thicksim \s$ for every $g\in G$.  Note that
any weight or gauge is
translationally equivalent, and that if $\s_{1}\thicksim
\s_{2}$, then $\s_{1}$ is translationally equivalent iff $\s_{2}$ is.
\subheading{Remark 1.2.8} The set of translationally equivalent scales on a
group $G$
(which we will denote by $\s$ here)
is closed under pointwise addition and multiplication.  The zero gauge acts
as the identity for addition, and the weight $\w\equiv 1$ is  the
identity for multiplication.  The set of gauges is an additive submonoid
of $\s$,
and the set of weights is a multiplicative submonoid of $\s$.  If we
denote by $\s/{\thicksim}$  the set $\s$ moded out by the  equivalence
relation on scales, then addition and multiplication are still defined
on $\s/{\thicksim}$.  One could do the same thing with the set of gauges
(for addition) and weights (for multiplication).
\vskip\baselineskip
\par From now on, assume that $\s$
is translationally equivalent.  Let $\alpha_{g}$
denote the automorphism $\varphi \mapsto \varphi_{g}$ of $\Cal S^{\s}(G)$.
We seek a condition on $\s$ so that the action $\alpha$ of $G$ on
$\Cal S^{\s}(G)$ is strongly continuous and infinitely differentiable.
By writing elements of $\schwsi$ as the sum of their real and imaginary
parts, it suffices to consider real valued functions.
For small $t$, we have
$$\split \pa \varphi_{exptX} - \varphi \pa_{d, \ga}& =
\int_{G}   \s^{d}(h)   | X^{\ga}( \varphi_{exptX} - \varphi )(h) | dh \\
& =   \int_{G}   \s^{d}(h)   | tX^{\tilde \ga}
( \varphi_{expt_{h}X})(h) | dh \leq |t|   \int_{G}
 \s^{d}((expt_{h}X)h)   | X^{\tilde \ga} \varphi(h) | dh,  \endsplit
\tag 1.2.9$$
where $t_{h}$ is a number between $-t$ and $t$ from the mean value theorem,
and $\tilde \ga$ is $\ga$ with a derivative in the $X$ direction added on.
If $\s((expt_{h}X)h)$ is bounded by a polynomial in
$\s(h)$ {\it uniformly in $t_{h}$},
then clearly 1.2.9 tends to zero as $t\rightarrow 0$.
Accordingly, we make the following definition.
\par
We say that $\s$ is {\it  uniformly translationally equivalent}
if for every compact subset $K$ of $G$ there exists $C, D>0$, $d\in \N$
depending on $K$,
such that
$$ \s_{g}(h) \leq C\s^{d}(h) + D,\qquad h \in G, \, g \in K.  \tag 1.2.10 $$
By (1.2.9), we see that if $\s$ is  uniformly translationally
equivalent,
then the action $\alpha$ on $\schwsi$ is
strongly continuous.  Similar calculations
show that the action is also differentiable.
\proclaim{Theorem 1.2.11} Let $G$ be a locally
compact group. If a scale $\s$ is translationally equivalent,
it is  uniformly translationally equivalent. (So weights and
gauges are uniformly translationally equivalent.)   Hence if $\s$ is
translationally equivalent and $G$ is a Lie group,
by our preceding remarks
the action $\alpha$ of $G$ on $\schwsi$ by
left translation is strongly continuous and differentiable.
Also, translationally equivalent scales are bounded an compacts sets.
(Hence weights and gauges are bounded on compact sets.)
\endproclaim
\demo{Proof}  (Compare
\cite{Dz 1, Prop 2.1}.)  Let $\s$ be a translationally
equivalent scale.  By replacing $\s$ with $\max(1, \s)$, we may assume
$\s \geq 1$.  Then for each $g\in G$ we have $C_{g} > 0$ and
$d_{g} \in \N$ such that
$$  \s_{g}(h) \leq C_{g} \s^{d_{g}}(h), \qquad h \in G.  \tag 1.2.12 $$
For each $g$, we may take $d_{g}$ to be the least positive integer that
satisfies (1.2.12) for some $C_{g}>0$.   Then, for this $g$ and $d_{g}$,
we let $C_{g}$ be the least positive integer satisfying (1.2.12).
\par
Note
$$ \s_{g_{1}g_{2}}(h) \leq C_{g_{2}} \s_{g_{1}}^{d_{g_{2}}}(h)
\leq C_{g_{2}} C_{g_{1}}^{d_{g_{2}}} \s^{d_{g_{1}}d_{g_{2}}} (h),
\qquad h \in G, \tag 1.2.13 $$
so
$$ d_{g_{1}g_{2}}\leq d_{g_{1}}d_{g_{2}}. \tag 1.2.14 $$
\par
Define
$$U_{m}= \{\, g \, | \quad d_{g} \leq	 m \quad \}  $$
$$V_{m}= \{\, g \, | \quad C_{g} \leq	 m \quad \}  $$
for $m \in \N^{+}$.  One easily checks from the measurability of
$\s$ that  $U_{m}, V_{m}$ are measurable subsets of $G$.
Also $e \in U_{m} \cap V_{m}$ and
$$ \cup_{m=0}^{\infty} (U_{m}\cap V_{m}) = G. \tag 1.2.15 $$
Since Haar measure is countably additive,
the Haar measure $|U_{m}\cap V_{m} |$ is
greater than zero
for some $m$. (Here we have used our assumption that $G$ is locally compact
for the existence of Haar measure.)
\par
For this $m$, the interior of $(U_{m}\cap V_{m})^{2}$ contains an open
neighborhood $V$ of the identity in $G$ \cite{Dz 2, p. 17 and 18}.
Then for $v \in V$, we have $v=v_{1}v_{2}$ with $v_{1}, v_{2} \in U_{m}\cap
V_{m}$.  Hence
$d_{v} \leq d_{v_{1}}d_{v_{2}} \leq m^{2}$ by 1.2.14, and
by (1.2.13),
$$ \s_{v}(h) \leq C_{v_{2}}C_{v_{1}}^{d_{v_{2}}} \s^{m^{2}}(h)
\leq m^{m+1} \s^{m^{2}}(h), \qquad h \in G. \tag 1.2.16 $$
(As far as I know, $d_{v}$ could be stricly less than
$m^{2}$, so (1.2.13) gives no bound on $C_{v}$.  However,
this does not interfere with the proof, since we have (1.2.16).)
Let $L$ be any compact subset of $G$, and let
$g_{1}, \dots g_{n} \in G$ be such that
$$\cup_{i=1}^{n} g_{i}V \supseteq L. $$
Say $g\in L$.  Then $g= g_{i}v$ for some $v\in V$ so
$$ \s_{g}(h) = \s_{v}(g_{i}^{-1}h)
\leq m^{m+1}\s^{m^{2}}(g_{i}^{-1}h), \qquad h \in G.
\tag
1.2.17$$
By the translational equivalence of $\s$, the quantity on the right hand
side of
(1.2.17) is bounded by a constant times a power of $\s(h)$.
Hence $\s$ is uniformly translationally equivalent.
\par
To see that translationally equivalent scales are bounded on compact sets,
let $L$ be any compact subset of $G$.
Let $C>0$ and $d \in \N$ be such that
$$ \s_{g}(h) \leq C\s^{d}(h), \qquad h \in G, \, g^{-1} \in L. \tag 1.2.18$$
Set $g = e$ in (1.2.18) and we see that $\s$ is bounded on $L$.
\qed
\enddemo
\subheading{Definition 1.2.19} Let $\s$ be any scale.  Define the
{\it $\s$-
rapidly vanishing $L^{1}$ functions $L_{1}^{\s}(G)$} on $G$ to be the
space of Borel measurable functions $f\colon G \rightarrow \C$
such that
$$ \pa f \pa_{d} = \int_{G} \s^{d}(g) \, | \, f(g) \, | \, dg \tag 1.2.20$$
is finite for all $d\in \N$.
Then $L_{1}^{\s}(G)$ is complete for the topology given
by the seminorms $\pa \quad \pa_{d}$ \cite{Schw, \S 5}.
Until \S 6, we shall often refer to $L_{1}^{\s}(G)$ as $L^{\s}(G)$.
The following characterization
of $\schwsi$ is useful (see for example Proposition 2.3.1), and it
also gives the completeness of the space $\schwsi$.
\proclaim{Theorem 1.2.21} Let $\s$ be a translationally equivalent scale
(for example a weight or a gauge) on a Lie group $G$.
Then the action $\alpha$ of $G$
by left translation is an isomorphism of
the locally convex spaces $\schwsi$ and $L^{\s}(G)$,
and this action of $G$ is
differentiable on $\schwsi$ and strongly continuous
on $L^{\s}(G)$.  The Fr\'echet space of $C^{\infty}$-vectors
$L^{\s}(G)^{\infty}$
is naturally isomorphic to $\schwsi$.  Hence $\schwsi$ is complete
and a Fr\'echet
space.
\endproclaim
\demo{Proof}  We saw that $\alpha$ acts differentiably on $\schwsi$.
By the estimate (1.2.6) with $\ga=0$, we see that each $\alpha_{g}$
is an automorphism of $L^{\s}(G)$.  Since $\alpha$ is strongly
continuous on $\schwsi$, and the inclusion map $\schwsi \hookrightarrow
L^{\s}(G)$
is continuous with dense image, it is easily checked that $\alpha$ acts
strongly continuously on $L^{\s}(G)$.
Since $\alpha$ acts differentiably on $\schwsi$, we  have $\schwsi
\subseteq L^{\s}(G)^{\infty}$.  To prove the theorem, it suffices to
show that $ L^{\s}(G)^{\infty}\subseteq \schwsi$.
See Theorem 2.1.5 for the proof of this.
\qed
\enddemo
\heading \S 1.3 Group Schwartz Algebras  \endheading
\par
We give appropriate conditions on the scale $\s$
so that $L^{\s}(G)$
and $\schwsi$ are  Fr\'echet *-algebras (not necessarily $m$-convex).
For $L^{\s}(G)$ and $\schwsi$ to be Fr\'echet algebras,
it suffices that $\s$ be a weight or a gauge, or more generally that
$\s$ be sub-polynomial (see below).
For $L^{\s}(G)$ to be a *-algebra, we require that
$\s$ be equivalent to the inverse scale $\s(g^{-1})$.  For $\schwsi$
to be a *-algebra, we
need the additional condition that $\s$ bounds $Ad$.
\vskip\baselineskip
\par
By a {\it Fr\'echet algebra}, we mean a Fr\'echet space
with an algebra structure for which the operation
of multiplication is continuous.
Multiplication is separately continuous if and only if
it is jointly continuous \cite{Wa}.
We do not require that Fr\'echet algebras have
submultiplicative seminorms (ie that they be $m$-convex) - see
\S 3.2.
A {\it Fr\'echet *-algebra}
is a Fr\'echet algebra with a continuous involution.
\par
We define an involution on functions on $G$ by $\varphi^{*}(g)=
\Delta(g){\overline \varphi}(g^{-1})$, where $\Delta$ is
the modular function on $G$.  If $\s$ is a scale, we define the
{\it inverse scale $\s_{-}$} by $\s_{-}(g) = \s(g^{-1})$.  We say that
$\s$ is {\it sub-polynomial} if there exists $C> 0$ and
$d \in\N$ such that
$$ \s(gh) \leq C(1+\s(g))^{d}(1+\s(h))^{d}, \qquad g, h \in G.\tag 1.3.1$$
If $\s \thicksim \tau$, then $\s$ is sub-polynomial if and only if
$\tau$ is.  Sub-polynomial scales are  translationally equivalent.
\proclaim{Theorem 1.3.2} If $\s$ is a
sub-polynomial scale,
then
$L^{\s}(G)$ is a
Fr\'echet algebra under convolution. If $\s$ is any
scale which is equivalent to its own inverse, then involution is
a well defined  continuous conjugate linear isomorphism of the
Fr\'echet space $L^{\s}(G)$.
In particular,
$L^{\s}(G)$ is a Fr\'echet *-algebra
if $\s$ is
a sub-polynomial scale such that $\s_{-} \thicksim \s$
(for example, if $\s$ is  a weight or a gauge).
It is in fact a dense *-subalgebra of $L^{1}(G)$,
and hence of the group C*-algebra $C^{*}(G)$ and the
reduced group C*-algebra $C^{*}_{r}(G)$.
 \endproclaim
\demo{Proof}
Without loss of generality, assume $\s \geq 1$.
Since $\s$ is sub-polynomial,
there exists $C>0$ and $d \in \N$ such that
$\s(gh) \leq C\s^{d}(g) \s^{d}(h)$.
By a change of variables, we have
$$\pa \varphi * \xi \pa_{m} \leq \int\int \s^{m}(g)
|\varphi(h)  \xi(h^{-1}g)| dhdg \leq C\pa \varphi \pa_{dm}
\pa \xi \pa_{dm}. $$
\par
For the statement about involution,
first replace $\s$ with the equivalent scale $\max(\s, \s_{-})$.
Then
$$\pa \varphi^{*} \pa^{\s}_{d} =
\int_{G} \s^{d}(g)\, |\varphi^{*}(g)|\, dg=
 \int_{G} \s^{d}(g^{-1})\, |\varphi(g)|\,
dg=\pa \varphi \pa_{d}^{\s_{-}}.$$
\qed \enddemo
\vskip\baselineskip
\par
\subheading{Example 1.3.3}
If $\s$ is not a translationally equivalent scale,
 $L^{\s}(G)$ is still a Fr\'echet space \cite{Schw, \S 5}.
We show that it is not in general an algebra
under convolution.   Let $G= \Z$.  Define  $\s(n) = e^{|n|^{|n|}}$.
We show that $L^{\s}(G)$ is not an algebra.
Let $\varphi(n) = e^{- |2n|^{|n|}}$.  Then $\varphi \in L^{\s}(G)$
since $\s^{d}(n) \varphi(n) = e^{d|n|^{|n|} - 2^{|n|} |n|^{|n|}} <
e^{d - 2^{|n|}}$ if $|n|$ is bigger than $1$.
We show that $\varphi * \varphi \notin L^{\s}(G)$.
$$\split \pa \varphi * \varphi \pa_{1} =& \sum_{n, m \in \Z}
\s(n) \varphi(m) \varphi(n-m) >
\sum_{n = 2m, m>0} \s(n) \varphi(m) \varphi(n-m)
\\&= \sum_{m>0} e^{|2m|^{|2m|}} e^{- |2m|^{|m|}} e^{- |2m|^{|m|}} =
\sum_{m>0} e^{|2m|^{|m|} (|2m|^{|m|} - 2)} > \sum_{m>0} e^{|2m|^{|m|}-2}.
\endsplit $$
The last sum diverges, so $L^{\s}(G)$ is not an algebra.   Note
that since $G=\Z$ is discrete, $\schwsi = L^{\s}(G)$ so $\schwsi $ is
not an algebra either.
\subheading{Question 1.3.4}If $L^{\s}(G)$ is a Fr\'echet algebra,
does $\s$ have to be translationally equivalent?
Does $\s$ have to be sub-polynomial?
%
\subheading{Example 1.3.5} If $\s$ is sub-polynomial, we shall see
below in Theorem 1.3.13 that $\schwsi$ is a Fr\'echet algebra.
However, we give an example showing that, unlike $L^{\s}(G)$, the algebra
 $\schwsi$ will {\it not} in general be a *-algebra if $\s \thicksim
\s_{-}$.
Let $G$ be the $ax+b$ group, namely
$2\times 2$ matrices over $\R^{2}$ of the form
$$ g= \pmatrix e^{a} & b \\ 0 & 1 \endpmatrix,
\qquad a, b\in \R.\tag 1.3.6 $$
Haar measure on $G$ is given by $e^{-a} da db$, where $dadb$ is Lebesgue
measure on $\R^{2}$.
The modular function $\Delta$ is given by $\Delta(g) = e^{a}$. (This is
consistent with the convention $ \int \varphi (g^{-1}) dg =
\int \Delta(g) \varphi(g) dg$, so that $\pa \varphi^{*} \pa_{1}
= \pa \varphi \pa_{1}$.)
 We let $\s$ be the constant weight $\s\equiv 1$.
 We show that $\schwsi$ is not closed
under the * operation $\varphi^{*}(g) =
\Delta(g) {\overline\varphi}(g^{-1})$.
\par
Define
$$\varphi (g) = {e^{a}\over{(1+a^{2})(1+b^{2})}}, $$
with $g$ as in (1.3.6).
We may think of the Lie algebra as
$2\times 2$ real matrices with the
second row zero.  Choose as a basis for the
Lie algebra the two matrices
$$ X_{1} = \pmatrix 1 & 0 \\ 0 & 0 \endpmatrix, \qquad
X_{2} = \pmatrix 0 & 1 \\ 0 & 0 \endpmatrix \tag 1.3.7$$
Then
$$\varphi(exp(-t_{2}X_{2}) exp(-t_{1}X_{1}) g) =
{e^{a+t_{1}} \over{(1+(a+t_{1})^{2})(1+(e^{t_{1}}(b-1)  +1 +t_{2})^{2})}},
$$
so since $X^{\ga}\varphi$ is given by differentiating in
$t_{1}$ and $t_{2}$, and then setting $t_{1}=t_{2}=0$ (see (1.2.1)
above with $\be_{g}(\varphi)(h)= \varphi(g^{-1}h)$), we see that
$\varphi$ and
all of its derivatives are integrable against the Haar measure.  Hence
$\varphi \in \schwsi$.
Since $\Delta(exp(-tX_{2})g) =
\Delta(g) = e^{a}$ is constant in $t$, we have
$$\split |\,X^{(0, 1)} \varphi^{*} & (g) \,|\,= \, |
\,\,\,\,{d\over{dt}}
\Delta(exp(-tX_{2})g) \varphi(g^{-1}exp(tX_{2})) \,|_{t=0}
\,\,\,\,|\\
& =\, |\, \Delta(g) (X^{(0, 1)}\varphi)(g^{-1}) e^{-a}\,|\,=
\,|\,(X^{(0, 1)} \varphi) (g^{-1})\,|,
\endsplit
\tag 1.3.8 $$
where we used $$g^{-1}exp(tX_{2}) = \pmatrix e^{-a} & e^{-a}(t-b) \\
0 & 1 \endpmatrix$$ and the chain rule for differentiation in the
second step.  By (1.3.8), we have
$$ \split \pa X^{(0, 1)}\varphi^{*} \pa_{1}
= &\int_{G} |(X^{(0,1)}\varphi) (g^{-1}) |dg \\=
 \int_{G} |(X^{(0,1)}\varphi)& (g) |\Delta(g) dg =
\int_{\R^{2}}\,
|\,\biggl[  {e^{a} \over{1+a^{2}}} {-2b\over{1+b^{2}}} \biggr]\,|\,dadb
= C \int_{\R^{2}}
{e^{a}\over{1+ a^{2}}}da,
\endsplit $$
which does not converge.  Hence $\varphi^{*} \notin \schwsi$ and
$\schwsi$ is not a *-algebra.
\vskip\baselineskip
\par
In order that $\schwsi$ be a *-algebra, we make
the following definition for a scale $\s$.  This definition
will also be important for making the smooth crossed product
$\sonas$ defined in \S 2 a Fr\'echet algebra (see for example
\cite{Sc 3, Thm 6.29}).
\subheading{Definition 1.3.9} Let $\s$ be any
 scale on the Lie group $G$.  We say that $\s$
{\it bounds $Ad$ } if
there exists $p\in \N$ and  positive constants $C, D$ such that
$$ \pa Ad_{g} \pa \leq C \s^{p}(g) + D,\qquad g\in G,
\tag 1.3.10 $$
where $\pa \quad \pa $  denotes
the operator norm
on the
space of linear operators on the Lie algebra of $G$.
If $G$ is any Lie group, then $\w(g) =  \max(\pa Ad_{g} \pa,
\pa Ad_{g^{-1}} \pa )$ is a weight on $G$ which bounds $Ad$.
If $G$ is compactly generated, then the exponentiated word length
dominates every weight on $G$ (see Theorem 1.1.21) and so bounds $Ad$.
Note that if $Ad_{g}$ always acts trivially (for example,
if $G$ is Abelian or discrete), then any scale $\s$
automatically bounds $Ad$.
If $\s$ is a scale on $G$ which bounds $Ad$, and $H$ is a
closed subgroup of
the identity component of $G$ which is normal in $G$, then
$\overline{\s}(g)=\inf_{h\in H}\s(gh)$ is a scale on
$G/H$ which bounds $Ad$.
Moreover, if $\s$ is a sub-polynomial scale, a weight or a gauge,
then $\overline{\s}$ is also. For example, if
$\s$ is submultiplicative, then
$$ \split {\overline{\s}}(g_{1}g_{2}) = \inf_{h_{1},h_{2}\in H}
\s(g_{1}g_{2}h_{1}h_{2})  =
\inf_{h_{1}, h_{2}\in H} &\s(g_{1}(g_{2}h_{1}g_{2}^{-1})g_{2}h_{2}) \\ &=
\inf_{h_{1},h_{2}\in H}
\s(g_{1}h_{1}g_{2}h_{2}) \leq \overline{\s}(g_{1})
\overline{\s}(g_{2}),\endsplit $$
so $\overline \s$ is submultiplicative.
See also Examples 1.3.14-16 below.
\par
If $\s$ is a weight or a gauge that
bounds $Ad$ (or more generally if $\s \thicksim \s_{-}$),
then clearly $\s_{-}$ bounds $Ad$ also.  We shall frequently use this,
since
one condition for the smooth crossed product
$\sonas$ in \S 2 to be an algebra will be that $\s_{-}$
bounds $Ad$.
\subheading{Example 1.3.11} We
show that in the in Example 1.3.5 above in which
$\schwsi$ is not a *-algebra,  the weight $\s$ does not bound $Ad$.
An easy calculation shows that for $g$ as in (1.3.6),
$$ Ad_{g^{-1}} =
 \pmatrix 1 & 0 \\ -b & e^{a} \endpmatrix \tag 1.3.12$$
where the matrix is in the  basis (1.3.7).
Clearly $\s\equiv 1$ does not bound $Ad$.
\proclaim {Theorem 1.3.13} Let $G$ be a Lie group, and
let $\s$ be a sub-polynomial
scale on $G$.
Then the Fr\'echet space
 $\schwsi$ is a Fr\'echet algebra under convolution.
If in addition $\s $ is equivalent to $\s_{-}$
and $\s$ bounds $Ad$ (for example, if $\s$ is
a weight or a gauge that bounds $Ad$), then $\schwsi$ is
a Fr\'echet *-algebra.
It is in fact a dense *-subalgebra of $L^{1}(G)$,
and hence of the group C*-algebra $C^{*}(G)$ and the
reduced group C*-algebra $C^{*}_{r}(G)$.
\endproclaim
\demo{Proof} See proof of Theorem 2.2.6 and
Corollary 4.9
setting $A=\C$.
\qed
\enddemo
\vskip\baselineskip
\par
We give some examples.
First, we note that if we take $G=\R$ and let $\w(r) = 1+|r|$,
then the Schwartz algebra $\schwL$ is precisely the standard
convolution algebra of Schwartz functions on $\R$.
\subheading{Example 1.3.14} Let $G=\R$.  Define $\w(r)=e^{|r|}$.
This is a  weight
on $G$ which bounds $Ad$.
The Schwartz algebra $\schwL$ is smaller than the
standard Schwartz algebra $\Cal S (\R)$.  It is the same
set of functions as $\Cal S_{exp}(\R)$ in \cite{DuC 3, \S 2.1.5}.
(In this case, it is irrelevant whether we use the $L^{1}$ norm
or the sup norm which  du Cloux uses - see Theorem 6.8.)
\subheading{Example 1.3.15}
If $H$ is any  subgroup of $G$,
such that the inclusion map $H\hookrightarrow G$
is differentiable,
then any scale, weight or gauge on $G$ which bounds $Ad$
restricts to a scale, weight or gauge on $H$ which bounds
$Ad$.  This is because for $h\in H$, the linear operator $Ad_{h}$
on the Lie algebra of $H$ is just
the restriction of $Ad_{h}$
as a linear operator on the Lie algebra of $G$.
This
allows us to define various gauges or weights (which bound $Ad$)
on  groups which are not
compactly generated, since such groups may occur
as subgroups of compactly
generated groups.  For example, the free group on
infinitely many generators is a subgroup of the compactly
generated group $SL(2, \R)$, or of the free group on two
generators.
\subheading{Example 1.3.16}
Let $G$ be a compact Lie group.  Then every scale
$\s$ bounds $Ad$ and
is equivalent to the zero scale, and $\schwsi$ is the Fr\'echet
space $C^{\infty}(G)$
of differentiable functions on $G$.   This is always a Fr\'echet
*-algebra under convolution.
\vskip\baselineskip
\heading \S 1.4 Gauges that bound $Ad$ \endheading
\par
As we have noted, there always exists a weight on $G$ which
bounds $Ad$.  In a later paper \cite{Sc 3}, it will be important to
know when there exists a gauge on $G$ which bounds $Ad$, since these will
yield Schwartz algebras which are spectral invariant \cite{Sc 2}
\cite{Sc 3}
in $L^{1}(G)$, and in $C^{*}(G)$.
We have noted many examples of groups which do possess gauges
which bound $Ad$ - for example the zero gauge bounds $Ad$
if $G$ is discrete, Abelian, or compact.
\par
We  note that gauges which bound $Ad$ do not always exist.  Assume
 there is some
$h \in G$ such that $Ad_{h}$ has an eigenvalue $\lambda$ not lying on
the unit circle.
Let $Y$ be the corresponding eigenvector.  Then
$$ \pa Ad_{h^{n}}Y\pa = \pa e^{n\lambda} Y\pa
= e^{nRe\lambda} \pa Y \pa. $$
Let $\epsilon = Re \lambda $.  By replacing $h$ with $h^{-1}$,
we may assume that
$\epsilon >0$. Assume that $\tau$ is a gauge on
$G$ which bounds $Ad$, and let $C, D>0$, $d\in \N$ be such that
$\pa Ad_{g} \pa \leq C\tau^{d}(g) + D$.  Then
$$ C\tau^{d}(h^{n}) + D
\geq \pa Ad_{h^{n}}\pa \geq e^{n \epsilon},  \tag 1.4.1$$
so $\tau(h^{n}) \geq
C^{1/d}(e^{n \epsilon} - D)^{1/d}$.  But $\tau(h^{n}) \leq n \tau(h)$
by subadditivity.  This is a contradiction.  Hence $G$ has no
gauge (or even a subadditive scale) which bounds $Ad$.
An example of such a group is
given by the $ax+b$ group above.  By (1.3.12),
it is clear that the eigenvalues of $Ad_{g}$ are not always on the
unit circle.
\proclaim{Definition 1.4.2} We say the a  Lie group
$G$ is {\it Type R } if $Ad_{g}$ has eigenvalues of complex modulus
one for all $g\in G$.  \endproclaim
\par
We shall show below (see Corollary 1.4.9) that
this generalizes the notion of a Type R connected Lie group \cite{Je}
\cite{Pa},
which says that the eigenvalues of $ad$ are imaginary.
In the connected case, a Lie
group is Type R if and only if it has polynomial growth \cite{Je}\cite{Pa}
(see \S 1.5 below).
Nilpotent Lie groups, the Mautner group,  motion groups \cite{Pa, p. 229},
and any discrete group are all examples of Type R
Lie groups.
If $H$ is any  subgroup of $G$,
such that the inclusion map $H\hookrightarrow G$
is differentiable, then $H$ is Type R if $G$ is Type R.
If $G$ is a Type R group, and $H$ is a closed subgroup of
the identity component of $G$ which is normal in $G$, then $G/H$ is Type R.
\proclaim{Theorem 1.4.3}  Let $G$ be a compactly generated
Lie group.  Assume that the quotient of $G$ by the kernel of $Ad$
has a cocompact solvable subgroup (in particular, it suffices that
$G$ have a cocompact solvable subgroup),  or assume
that $G$ is connected.
Then  the following are equivalent.
\roster
\item The group $G$ is  Type R.
\item The  word gauge bounds $Ad$.
\item There exists a gauge on $G$ which bounds $Ad$.
\endroster
The implication $
(3) \Rightarrow (1)$ is true for any Lie group.
\endproclaim
\demo{Proof} We have already seen that $ (3) \Rightarrow  (1)$ above.
Also
$(2) \Leftrightarrow (3)$ is trivial, since the word gauge dominates
every gauge (see Theorem 1.1.21).
For
$(1) \Rightarrow (2)$, we first prove two lemmas.
\proclaim{Lemma 1.4.4 (Compare \cite{Lu 1, Lemma 1.5})}
For $q \in \N$, let $G$ be the group of upper triangular matrices
$T_{1}(q, \C)$  with $1$'s on the diagonal.
Then $G$ has a gauge (the word gauge)
$\tau$ which is equivalent to the weight given
by the norm on matrices. \endproclaim
\demo{Proof}
For $g\in G$, let $\pa g \pa_{\infty}$ denote the maximum
absolute value of the {\it off diagonal} matrix entries $|g_{ij}|$.
Let $U$ be any symmetric open neighborhood of the identity
of $G$, such that $U\subseteq \{ \, g \, | \quad \pa g \pa_{\infty}
\leq 1 \quad \}$.
Let $\tau$ be the word gauge with respect to $U$.  We show that $\tau$
bounds the norm $\pa \quad \pa_{\infty}$.
\par
First by induction on $n$, we have for $1<i<k<q$,
$$  (g_{1}\dots g_{n})_{ik} = \sum_{i=i_{1}\leq\dots i_{n+1}=k}
(g_{1})_{i_{1}i_{2}}\dots (g_{n})_{i_{n}i_{n+1}},\tag 1.4.5 $$
where $g_{1},\dots g_{n} \in G$.
\par
Assume that $\tau(g)=n$ for some $g= g_{1}\dots g_{n}$,
with $g_{i}\in U$.  Then we have
$$ \aligned \pa g \pa_{\infty}&
\leq \max_{1<i<k<q} |(g_{1}\dots g_{n})_{ik}|
\leq \max_{1<i<k<q}
\biggl( \sum_{i=i_{1}\leq\dots i_{n+1}=k}(1) \biggr)\leq
\sum_{i=1_{1}\leq\dots i_{n+1}=q}(1)\\
& \leq \text{\# of possible ways to arrange
steps times \# of possible step sizes}\\
& =
{\biggl({\binom{n-1} {q-1} }
+ \dots {\binom{n-1} {1}} \biggr)}
 q = P_{q}(n) = P_{q}(\tau(g)), \endaligned$$
where $P_{q}$ is a polynomial, of degree $q$, depending only on
$q$.  Thus $\tau$ dominates the sup norm $\pa \quad \pa_{\infty}$.
\vskip\baselineskip
\par
We show that $\tau$ is dominated by the operator norm $\pa \quad \pa$
on matrices.  Let $C>1$ be such that $
 \{ \, g \quad | \quad \pa g \pa
\leq C\, \}\subseteq U$.  Let $T_{0}(n, \C)$ denote the
upper triangular complex matrices with $0$'s on the diagonal.
Recall that
$$ exp\colon T_{0}(q, \C) \rightarrow G $$
is a polynomial map with polynomial inverse.  In fact,
$$\pa exph \pa \leq e^{\pa h \pa}, \qquad  h \in T_{0}(q, \C),  $$
and
$$\pa log(g) \pa \leq P(\pa g \pa ), \qquad  g\in G,  $$
where $P$ is a polynomial.  Let $g\in G$ and $h=\log (g)$.
Let $d$ be the least integer greater than or equal to $\pa h \pa/\log(C)$.
Then $\pa e^{h/d}\pa \leq e^{\pa h \pa /d} \leq e^{\log(C)} =C$ so
$\tau(e^{h/d})\leq 1$.  Thus
$$ \tau (g) = \tau((e^{h/d})^{d}) \leq d \leq {\pa h \pa/{\log(C)} }  + 1
\leq {P(\pa g \pa)/{\log(C)}} + 1. $$
This proves Lemma 1.4.4.
\qed \enddemo
\par We say that a topological group $H$ is a {\it subgroup} of
a topological group $G$ if $H\subseteq G$ with continuous inclusion.
We say that a scale $\s$ is a {\it one-sided gauge} iff $\s$ satisfies
the inequality
$$ \s(gh) \leq C + D\s^{d}(g) +\s(h), \qquad g, h \in G, $$
for some $C, D>0$ and $d \in \N$.  Note that one-sided gauges are
translationally equivalent and therefore bounded on compact
subsets by Theorem 1.2.11.
\proclaim{Lemma 1.4.6} Let $G$ be any separable subgroup of $GL(q, \C)$
and assume that
$G$ has a closed cocompact subgroup $H$ with a
one-sided gauge which is equivalent to the matrix norm inherited from
$GL(q, \C)$.  Then $G$ has a one-sided gauge which is equivalent to
the matrix norm inherited from $GL(q, \C)$.
\endproclaim
\demo{Proof} Let $B\subseteq G$ be a rel. comp. measurable
cross section for the cosets $G/H$ \cite{Pa, App. C}.  (Here
we have used  $H$ closed (to get $G/H$ Hausdorff) and $G$ separable.)
 Let $\tau$
be a one-sided gauge on $H$ which is equivalent to the norm on
matrices.  Let $C, D>0$ and $d\in \N$ be such that
$$\tau(h_{1}h_{2}) \leq C + D\tau^{d}(h_{1}) + \tau(h_{2}), \qquad
h_{1}, h_{2}\in H. $$
Let $g\in G$, and write $g=bh$ with $b\in B$, $h\in H$.
Define
$$\ga(g) \equiv \tau(h). $$
Since $B$ is measurable, $\ga$ is measurable and hence a scale on $G$.
\par
We show that $\ga$ is a one sided gauge.
Since $B^{-1}B$ and $B$ are rel. comp. in $G$, they are both
rel. comp. and hence norm bounded sets in $GL(q, \C)$.
Thus, since $\tau$ is equivalent to the norm on matrices,
there exists a polynomial $P$ such that
$$ \tau(b_{1}hb_{2})\leq P(\tau(h)) $$
for all $h\in H$ and all $b_{1}\in B^{-1}B$ and $b_{2}\in B$
for which the product $b_{1}hb_{2}$ is in $H$.
Let $g_{1} , g_{2} \in G$ and write $g_{1}=b_{1}h_{1}$, $g_{2}=b_{2}h_{2}$,
$g_{1}g_{2} = b_{3}h_{3}$.  Then
$$\aligned
\ga(&g_{1}g_{2})= \tau(h_{3}) = \tau(b_{3}^{-1}b_{1}h_{1}b_{2}h_{2})
\leq C + D\tau^{d}(b_{3}^{-1}b_{1} h_{1} b_{2})   + \tau (h_{2}) \\
& \leq C + D (P(\tau(h_{1}))^{d} + \tau(h_{2}) =
C + D(P(\ga(g_{1}))^{d} + \ga(g_{2})
\leq C^{\prime} + D^{\prime}\ga^{k}(g_{1})   + \ga(g_{2})  \endaligned$$
for some  $k \in \N$ and $C^{\prime}, D^{\prime}>0$.
So $\ga$ is a one-sided gauge on $G$.
\par
A similar argument using the polynomial $P$ and the fact that
$\tau$ is equivalent to the norm on matrices, shows that $\ga$ is
equivalent to the norm on matrices.
\qed
\enddemo
\par
Finally, we prove $(1)\Rightarrow (2)$ of Theorem 1.4.3
using Lemmas 1.4.4 and 1.4.6.  Let $G$ be any
compactly generated Lie group.  Let $G^{\prime}$ denote
the image of $G$ in $GL(q, \C)$ via the $Ad$ map, where
$Ad_{g}$ acts on the complexification of the Lie algebra.
Then $G^{\prime}$ has a cocompact (closed since
we may take the closure of $G''$ in $G'$)
solvable subgroup $G^{\prime\prime}$ by assumption.
(We shall do the connected case later.)  Then $G^{\prime\prime}$
has a subgroup $G^{\prime\prime\prime}$ of finite index in
$G^{\prime\prime}$ which is triangulable in some basis over $\C^{q}$
\cite{KM, Thm 21.1.5}.
Since we are assuming $G$ is Type R, each triangular matrix in
$G^{\prime\prime\prime}$ has entries of modulus one along the diagonal.
By Lemmas 1.4.4 and 1.4.6, there is a one-sided
gauge equivalent to the norm on matrices on the group of
upper triangular matrices in $GL(q, \C)$, with entries
of modulus one on the diagonal.
By restriction, we obtain a one-sided gauge on
$G^{\prime\prime\prime}$ which is equivalent to the norm on matrices.
Since $G^{\prime\prime\prime}$ is cocompact and closed in $G^{\prime}$,
and since $G^{\prime}$ is separable (because $G$ is a compactly
generated Lie group, and so separable), we may apply Lemma 1.4.6 to
get a one-sided gauge on $G^{\prime}$ which
is equivalent to the norm on
matrices.  We pull this gauge back to $G$ to get a one-sided gauge on
$G$ which bounds $Ad$.  It remains to show that the word gauge dominates
this one-sided gauge.  Let $\ga$ be the one-sided gauge on $G$, and
let $\tau $ be the word gauge with respect to a generating set
$U$.  Let $K=\sup \{\ga(g)\, | \, g \in U\}$.  Then $K<\infty$ by
Theorem 1.2.11, and we have
$$\split \ga(g_{1} \dots g_{n} )& \leq (n-1)C + D\ga^{d}(g_{1}) + \dots
D\ga^{d}(g_{n-1}) + \ga(g_{n}) \\
&
\leq
nC +DK^{d}(n-1) + K = \tau(g) C + D K^{d} (\tau(g)-1) +K\endsplit $$
if $\tau(g)=n$.  Hence $\tau$ strongly dominates $\ga$, and $\tau$
bounds $Ad$.
\par
Finally, assume $G$ is
connected and Type R.  Then by \cite{Pa,  Prop 6.29},
$G$ has cocompact radical.  Hence $G$ has a cocompact solvable subgroup,
and so does the image of $G$ via $Ad$.  This proves Theorem 1.4.3
\qed \enddemo
\subheading{Question 1.4.7} Is it true for a  general compactly
generated Lie group $G$ that $G$ has a gauge
which bounds $Ad$ iff $G$ is Type R ?
\vskip\baselineskip
\proclaim{Corollary 1.4.8}  Assume that $G$ is a compactly generated
Lie group.  If $G$ has a closed cocompact
solvable subgroup $H$, such that for every $h \in H$ the eigenvalues of
$Ad_{h}$, {\it acting an operator on the Lie algebra of $G$}, lie on the
unit circle,
then $G$ is Type R and has a gauge that bounds $Ad$.
\endproclaim
\demo{Proof}
Let $G_{1}$ be the image of $G$ in $GL(q, \C )$
via the $Ad$ representation.  Let $H_{1}$ denote
the image of the cocompact subgroup $H$ via this same
representation.  Replace $H_{1}$ with its closure in
$G_{1}$ if necessary.  Then $H_{1}$ is a subgroup of $GL(q, \C)$
which is solvable and all of whose eigenvalues lie on the unit
circle.  By Lemmas 1.4.4 and 1.4.6 above, $H_{1}$ has
a one-sided gauge which is equivalent to the norm on matrices.
Hence by Lemma 1.4.6, $G_{1}$ has a one-sided gauge that
is equivalent to the norm on matrices.  We pull this back to $G$
to get a one-sided gauge which bounds  $Ad$.
Then, since $G$ is compactly generated, we saw above
(at the end of the proof of Theorem 1.4.3)
that the
word gauge dominates this one-sided gauge.  Therefore the word gauge
bounds $Ad$.  Since the existence of a gauge that bounds $Ad$ implies that
$G$ is Type R, we are done.
\qed
\enddemo
We say that a connected Lie group $G$ is $ad$ Type R
if all the eigenvalues of $ad$ are imaginary.
\proclaim{Corollary 1.4.9} A connected Lie group is Type R if
and only if it is $ad$ Type R.
\endproclaim
\demo{Proof}   Clearly Type R
implies $ad$ Type R if $G$ is connected, since the spectrum
of $Ad_{expX} = e^{adX}$ is just the set $e^{spec(ad(X))}$.
Assume that $G$ is $ad$ Type R.
Let $H$ denote the radical of $G$, and
let $\frak H$ denote the Lie algebra of $H$.
Then
the set $\{ expX \, | \, X \in $H$ \}$ is dense in $H$ by \cite{Pa, App E}.
 Since $G$ is $ad$ Type R, $Ad_{expX}$, as an operator on
the Lie algebra $\frak G$ of $G$, has all its eigenvalues lying
on the unit circle.  Since for $h \in H$, $Ad_{h}$
is the norm limit of linear operators
on $\frak G$ having eigenvalues on the unit circle,  $Ad_{h}$ must also have
eigenvalues on the unit circle.
By \cite{Pa, Prop 6.29}, $H$ is cocompact in $G$.
Hence by Corollary 1.4.8,
$G$ is Type R.
\qed
\enddemo
\vskip\baselineskip
\heading \S 1.5 Polynomial Growth Groups \endheading
\par
We exhibit a large class of Type R groups, all of which
have gauges that bound $Ad$.
We say that a locally compact group has {\it polynomial growth} if
for every relatively compact
neighborhood $U$ of the identity, the Haar measure
$|U^{n}|$ is bounded by a polynomial in $n$\cite{Pa}.  Examples of
polynomial growth groups are given by
finitely generated polynomial growth discrete groups \cite{Ji},
closed subgroups of
nilpotent Lie groups, connected Type R Lie groups,
motion groups, the Mautner group, and compact groups.
\proclaim{Proposition  1.5.1}  Let $G$ be a locally compact group.
If $G$ has a gauge $\tau$ that satisfies the integrability
condition
$$ (\exists p \in \N) \quad \int_{G} {1\over{(1+\tau(g))^{p}}} dg<\infty,
\tag 1.5.2$$
then
$G$ has polynomial growth.  If $G$ is compactly generated, then $G$ has
polynomial growth if and only if $G$ has a gauge which satisfies (1.5.2).
If $G$ has a gauge $\tau$ that
satisfies (1.5.2) and $H$ is a closed subgroup
of $G$, then the restriction of $\tau$ to $H$ also satisfies (1.5.2)
for some possibly larger $p \in \N$.
\endproclaim
\demo{Proof}
Let $\tau$ be a gauge on $G$ which satisfies (1.5.2).  Let $U$ be any
rel. comp. neighborhood of $e$.  Then for $n \in \N$,
$$\int_{U^{n}} dg \leq \sup_{g \in U^{n}}(1+\tau(g))^{p}
\int_{G}{1\over{(1+\tau(g))^{p}}} dg \leq C\bigl(1+n\sup_{g \in U}
\tau(g)\bigr)^{p}
\leq C(1+nD)^{p},$$
where we have used the fact that gauges are bounded on compact sets (see
Theorem 1.2.11).
So $G$ has polynomial growth.
\par
Assume $G$ is compactly generated and has polynomial growth.  Let  $U$ be a
generating set, and let $\tau$ be the word gauge. Let $p$ be such that
$|U^{n}|\leq C(1+n)^{p}$.  Then
$$\split \int_{G}{1\over{(1+\tau(g))^{p+2}}} dg
&\leq \sum_{n=0}^{\infty}
\int_{U^{n}-U^{n-1}} {1\over{(1+\tau(g))^{p+2}}} dg\\&
\leq \sum_{n=0}^{\infty} {|U^{n}|\over{(1+n)^{p+2}}}
\leq   \sum_{n=0}^{\infty} {C\over{(1+n)^{2}}}<\infty.\endsplit $$
\par
For the last statement of the proposition, see Proposition 6.13 below.
\qed
\enddemo
\subheading{Question 1.5.3} Does a polynomial growth
group always have a gauge which satisfies the integrability
condition (1.5.2) ?
\proclaim{Lemma 1.5.4\cite{Pa, Prop 6.6,6.9}}
If $G$ has polynomial growth, then $G$
is unimodular.
\endproclaim
The following proposition is similar to
\cite{Pa, Prop 6.20(ii)}.
\proclaim{Proposition 1.5.5}
Let $G$ be a locally compact compactly generated
polynomial growth group with closed normal subgroup $H$.
Assume that the identity component $G_{0}$ is open in $G$. Then
the quotient group $G/H$ has polynomial growth. \endproclaim
\demo{Proof}
By \cite{Bou, chap VII, \S 2, thm 2},
$$\int_{G} f(g) dg = \int_{G/H} \int_{H} f(gh)  dh dg, \tag 1.5.6 $$
for $f \in L^{1}(G)$.
If $S \subseteq H$, define
$$ \lambda_{H}(S) = \int_{S}  dh.  $$
Let $V=V^{-1}$ be a generating set for $G$.  Let $U= \pi(V)$, where
$\pi \colon G \rightarrow G/H$ is the canonical map.
For $\delta >0$, define
$$ U_{\delta} = \{ \pi(g) \in G/H \, | \, \lambda_{H}(g^{-1}V\cap H) \geq
\delta
\quad \text{and} \quad
\lambda_{H}(gV\cap H) \geq \delta,\,\,
g \in G \,\}. $$
If $\pi(g) \notin U$, then $g^{-1}V\cap H= \varnothing $.
So $U_{\delta} \subseteq U$ for each $\delta>0$.
Also, if $\pi(g) \in U$, then
$g^{-1}V$ and $gV$ are open sets whose intersection with $H$
is nontrivial.
Since $H$ is a closed subgroup of $G$ (with
relative topology), the sets $g^{-1}V\cap H$ and
$gV \cap H$ are nonempty
open sets in $H$.  Hence $\lambda_{H}(g^{-1}V\cap H) >0$
and $\lambda_{H}(gV\cap H) >0$.  It follows that
$$ \bigcup_{\delta>0} U_{\delta} = U. \tag 1.5.7$$
For sufficiently small $\delta$, $U_{\delta}$ has nonzero measure in
$G/H$ (by (1.5.7) and since $U$ is a generating set for $G/H$
by Theorem 1.1.21).
It follows that $U_{\delta}^{2}$ contains a neighborhood of zero
\cite{Dz 2, p.17-18}.
Since $U_{\delta} = U_{\delta}^{-1}$,
the set $\bigcup_{n=0 }^{\infty} U_{\delta}^{n}$
is a nonempty open subgroup of $G/H$, and hence must
contain the connected component
of the identity of $G/H$.   We show that for $\delta$ sufficiently
small, $U_{\delta}$ also generates all of $G/H$.  Let $q$ be the
canonical projection
$$ q \colon G/H \rightarrow (G/H)/(G/H)_{0}.$$
By (1.5.7) and since
$(G/H)/(G/H)_{0}$ is discrete,
for $\delta$ sufficiently small we must have
$q(U_{\delta}) = q(U)$.
(To obtain the discreteness, we used that $(G/H)_{0}$ is open,
which follows from our assumption that $G_{0}$ is open.)
Then $q(U_{\delta})$ generates the  discrete group
$(G/H)/(G/H)_{0}$.  It follows that
$U_{\delta}$ generates $G/H$.  From now on, we fix $\delta$
sufficiently small so that $U_{\delta }$ generates $G/H$.
\par
To see that $G/H$ has polynomial growth, it suffices to show that
the Haar measure of $U_{\delta}^{n}$ is bounded by a polynomial in $n$
for $n \in \N$.
To do this, we first
show that if $\pi(g) \in U_{\delta}^{n}$, then
$\lambda_{H}(g^{-1}V^{n} \cap H) \geq \delta$.  Let
$g = g_{1} \dots g_{n}$ where $\pi(g_{i}) \in U_{\delta}$.
Then $g_{i} = {\tilde g}_{i} h_{i}$ where ${\tilde g}_{i} \in V$
and $h_{i} \in H$.  So, using the fact that $H$ is a normal
subgroup,  $g^{-1}V^{n} \cap H = {\tilde h} {\tilde g}_{n}^{-1}
\dots {\tilde g}_{1}^{-1} V^{n} \cap H \supseteq {\tilde h}
{\tilde g}_{n}^{-1} V \cap H$ for some ${\tilde h}\in H$.
By left invariance of Haar measure on $H$,
$$\lambda_{H}(g^{-1} V^{n}\cap H) \geq \lambda_{H}
({\tilde g}_{n}^{-1}V\cap H) \geq \delta, \tag 1.5.8$$
since $\pi({\tilde g}_{n}) = \pi(g_{n}) \in U_{\delta}$.
\par
If $S$ is a set, we let $\chi_{S}$ denote the characteristic
function of $S$.  We estimate
$$ \aligned \int_{V^{n}} dg& =\int_{G} \chi_{V^{n}}(g) dg = \int_{G/H}
\int_{H} \chi_{V^{n}}(gh) dh dg \\
&
= \int_{G/H} \lambda_{H}(g^{-1}V^{n} \cap H) dg \geq
\int_{U_{\delta}^{n}} \lambda_{H}(g^{-1}V^{n} \cap H) dg \geq
\int_{U_{\delta}^{n}} \delta dg, \endaligned\tag 1.5.9 $$
where we used (1.5.8) in the last step.
The last expression in (1.5.9) is the
Haar measure of $U_{\delta}^{n}$ in $G/H$ times $\delta$.
Since $U_{\delta}$ is
a generating set, $G/H$ has polynomial growth.
\qed
\enddemo
\proclaim{Theorem 1.5.10(\cite{Lo, Cor to Thm 2})}
Let $G$ be a compactly generated
polynomial growth
Lie group. Then there exists a normal series of
closed subgroups $G_{1}  \vartriangleleft G_{2}\vartriangleleft G$
such that $G_{1}$ is a  connected solvable polynomial growth
Lie group, $G_{2}/G_{1}$ is a nilpotent finitely generated discrete
group, and $G/G_{2}$ is compact.
\endproclaim
\demo{Proof}
See \cite{Lo, Cor to Thm 2} and its proof (last paragraph of \cite{Lo}).
\qed
\enddemo
\proclaim{Corollary 1.5.11} Let $G$ be a compactly generated polynomial
growth Lie group. Then any closed subgroup of $G$ is compactly generated.
In particular, closed subgroups of compactly generated polynomial growth
Lie groups are compactly generated polynomial growth Lie groups.
\endproclaim
\demo{Proof}
The second statement follows from the first
and Proposition 1.5.1.  By \cite{BWY,
Thm 2.11, Cor 2} it suffices to show that
every closed subgroup of the group
$G_{2}$ from Theorem 1.5.10 is compactly generated.
\par
Note that since $G_{2}/G_{1}$
is a finitely
generated nilpotent discrete group, any closed subgroup of it is finitely
generated.
Also, by \cite{BWY, Prop 2.1, Cor 3} and
since $G_{1}$ is solvable and connected,
every closed subgroup of $G_{1}$ is compactly generated.  Let $F$ be any
closed subgroup of $G_{2}$.  We apply \cite{BWY, Prop 2.1} to see that
if $FG_{1}$ is closed in $G_{2}$, then $F$ is compactly generated.
\par
It remains to show that $FG_{1}$ is
closed.  If $f_{n}g_{n} \rightarrow g_{0}\in
G_{2}$, then $f_{n}g_{n}$ must eventually lie in the same $G_{1}$
coset $c$ of $G_{2}$. (Recall that $G_{1}$
is the connected component of the
identity of $G_{2}$.)  But $c \cap FG_{1} = c$, so $g_{0} \in FG_{1}$.
\qed
\enddemo
\proclaim{Corollary 1.5.12} If $G$  is a compactly
generated Type R polynomial growth Lie group,
then $G$  has a gauge which bounds $Ad$. \endproclaim
\demo{Proof} By Theorem 1.5.10, we know $G$ has a cocompact closed solvable
subgroup $H$.  Since $G$ is Type R, the conclusion follows from
Theorem 1.4.3.
\qed
\enddemo
\par
Every closed subgroup of a connected nilpotent Lie group
satisfies the hypotheses of Corollary 1.5.12.  So do all
connected polynomial growth Lie groups, and hence the Mautner
group and all motion groups.
\par
It is not true in general that a  Lie group
is Type R if and only if it has polynomial growth.  For example,
any discrete group is Type R.
An example of a polynomial growth group that is not Type R
is given in \cite{Lo, Ex 1}.
However, we have the following theorem.
\proclaim{Theorem 1.5.13} Let $G$ be a compactly generated
polynomial growth Lie group, and assume that the
connected component of the identity $G_{0}$ is
solvable and simply connected.  Then $G$ is Type R and has a gauge that
bounds $Ad$.
\endproclaim
\demo{Proof}
The maximal compact subgroup of the center of $G_{0}$ is trivial, so
$G$ is Type R by \cite{Lo, Thm 1 $a)\Rightarrow b)$}.
\qed
\enddemo
\subheading{Example 1.5.14}
If $G$ is nilpotent and connected,
then the Schwartz algebra $\schw$ we obtain using the word
gauge agrees with the usual notion of the Schwartz algebra of $G$
\cite{Lu 1}, \cite{Ho, p. 346}.  In this case, $G$ is easily seen to
be Type R, and has polynomial growth
by \cite{Pa, Thm 6.17, 6.39}.
\subheading{Example 1.5.15}
Let $G$ be the three dimensional Heisenberg Lie group, that is the
set of matrices of the form
$$g = \pmatrix 1 & a & b \\ 0 & 1 & c \\ 0 & 0 & 1 \endpmatrix, $$
where $a,b,c \in \R$.  Then $G$ satisfies the hypotheses of
Theorem 1.5.13.
The word gauge is equivalent to the scale
$\tau(g) =|a|+|c| + |b|^{1/2}+ |b - ca|^{1/2}.$
We omit the proof of this, since the argument is the
same as provided for the $ax+b$ group below
in Example 1.6.1(or see \cite{Lu 1,
Lemma 1.5}).  A quick calculation shows that $\tau(e) = 0$,
$\tau(g^{-1})= \tau (g)$ and $\tau(gh )\leq 3(\tau(g) + \tau(h))$,
so $\tau$ is very much like a gauge.
In an appropriate basis for the Lie algebra, we have
$$ Ad_{g} = \pmatrix 1 & 0 & 0 \\ a & 1 & -c \\ 0 & 0 & 1 \endpmatrix. $$
One can then see directly that $\tau$ bounds $Ad$.
\vskip\baselineskip
\heading \S 1.6 Examples of Weights that Bound $Ad$ \endheading
\par
Now we look at several groups which are not Type R.
Our primary concern
will be to give explicit formulas for the exponentiated word weight,
which we will refer to as the word weight for brevity.
\subheading{Example 1.6.1} Let $G$ be
the $ax+b$ group (see Example 1.3.5 above).
Define $\w(g) = e^{|a|} + |e^{-a}b| +|b| + 1$.
Then $\w$ is a continuous weight on $G$, which clearly bounds $Ad$
by (1.3.12).
We may form
the Schwartz algebra $\schwL$.  This is the same
set of functions  as the $\Cal S$ of \cite{DuC 3, \S 5.2.5}.
(Again, du Cloux uses the sup norm, and we use the $L^{1}$ norm,
but we shall see that this makes no difference by Theorem 6.8 and
Proposition 6.13 (1).)
We have seen that   since $G$ is not Type R,
 $\w$ is not equivalent to any gauge.
\par
We show that $\w$ is equivalent to the  word weight.
Let $\tau$ be the word gauge.  We get an upper bound on
$\tau$.
Let $U$ be a generating set containing all matrices of the form
$$ \pmatrix e^{a} & b \\ 0 & 1 \endpmatrix,
\qquad
|a|, \quad |b | \leq 1. $$
By the definition of $\tau$, we have
$$ \tau(g) =\tau\pmatrix e^{a} & b \\ 0 & 1 \endpmatrix
\leq \tau\pmatrix e^{a} &0 \\ 0 & 1 \endpmatrix
+\tau\pmatrix 1 &e^{-a} b \\ 0 & 1 \endpmatrix
\leq |a| +\tau\pmatrix 1 &e^{-a} b \\ 0 & 1 \endpmatrix. \tag 1.6.2 $$
Choose any  $n\in \N$ for which there exists a $\ga$ satisfying
$1/e<|\ga|\leq 1$ and
$$\pmatrix 1 &(e^{n-1}+\dots 1)\ga  \\ 0 & 1 \endpmatrix
= \pmatrix 1 &e^{-a} b \\ 0 & 1 \endpmatrix. \tag 1.6.3 $$
Then $\pmatrix 1 &e^{-a} b \\ 0 & 1 \endpmatrix\in U^{2n}$ since
$$\pmatrix 1 &(e^{n-1}+\dots 1)\ga  \\ 0 & 1 \endpmatrix
=
{\undersetbrace \text{$n$ times}
\to{\pmatrix e &\ga  \\ 0 & 1 \endpmatrix\dots
\pmatrix e&\ga  \\ 0 & 1 \endpmatrix}}
{\undersetbrace \text{$n$ times}
\to{\pmatrix e^{-1} &0 \\ 0 & 1 \endpmatrix
\pmatrix e^{-1} &0 \\ 0 & 1 \endpmatrix}}.$$
So by (1.6.2), we have $\tau(g) \leq |a| +2n$.  By
(1.6.3), we see that
$$\biggl({e^{n}-1\over{e-1}}\biggr) \ga = e^{-a}b$$
so $e^{n} \leq e(e-1) (e^{-a}b) +1$.
We bound the word weight:
$$e^{\tau(g)}\leq e^{|a|+2n} \leq e^{|a| }(e(e-1) (e^{-a}b) +1)^{2}
\leq (e(e-1))^{2} \w^{2}(g).$$
Since the word weight dominates every weight on $G$,
it follows that  $\w$ is equivalent to
the word weight on $G$.
\vskip\baselineskip
\par
We give some examples of nonsolvable Lie groups.
\subheading{Example 1.6.4} Let $G = SL(n, \R)$.  We show that the
word gauge $\tau$ is strongly equivalent to the Riemannian symmetric
function $\s$.  It suffices to show that $\s$ strongly
dominates $\tau$.
Let $g \in G$.  Let $g = UP$ be the Cartan decomposition
of $g$, and let $V$ be a unitary matrix such that $V^{*}PV$ is
a positive diagonal matrix $D$.  Then by
\cite{HW, (2.2ab)} we have
$$\s(g) = \s(P) = \s(D) = \pa \log (D) \pa = \max_{i=1, n} |\lambda_{i} |,
\tag 1.6.5
$$
where
$$ D = \pmatrix e^{\lambda_{1}}  & \dots  & 0
\\ \vdots  & \ddots & \vdots \\ 0  & \dots
& e^{\lambda_{n}} \endpmatrix. $$
Let $A$ be a generating set for $G$.
Then if $\Cal U$ denotes the set of unitary matrices,
we have that
$$ \Cal U A \Cal U \cap G$$
is a bounded subset of $G$ and so contained in $A^{m}$ for some $m$.
Let $k$ be the smallest integer greater than or equal to $\s(g)$.
Then we have
$$ \tau(g) = \tau ( UVPV^{*}) \leq m\tau(P) \leq mk\tau(P^{1/k})
\leq mkC \leq mC (\s(g) + 1),\tag 1.6.6 $$
where $C$ is a constant which depends only on the norm of $P^{1/k} $,
which is less than or equal to $e$ by (1.6.5).
Hence (1.6.6) holds for all $g \in G$ and $\s$ strongly dominates
$\tau$.
\par
It follows that $\s\thicksim_{s} \tau$ since $\tau$ strongly dominates
every gauge on $G$.  Thus if $\w= e^{\tau}$ is the word weight,
we have
$$ \w \thicksim e^{\s}. \tag 1.6.7 $$
\par
Define another weight $\theta$ on $G$ by $\theta(g) = \max (\pa g \pa,
\pa g^{-1} \pa )$.  We show that $\theta$ is also equivalent to $\w$.
By (1.6.7), it suffices to show that $\theta$ dominates $e^{\s}$.
We have
$$ e^{\s(g)} = e^{\s(D)} = e^{\max |\lambda_{i} |}
\leq \max (e^{\lambda_{i}}, e^{-\lambda_{i}})
= \max(\pa D \pa, \pa D^{-1} \pa)
= \theta(D) \leq K \theta(g),\tag 1.6.8  $$
where $K$ is any constant
greater than or equal to  the square of $\sup_{h \in \Cal U} \theta(h)$.
Since (1.6.8) holds for all $g \in G$, it
follows that $\theta \thicksim \w$.
\par
Since $\w$ dominates every weight on $G$, we know $\w$ bounds $Ad$.
We see this directly for the case $n = 2$.
A calculation shows that
if
$$g = \pmatrix E & F \\ G & H \endpmatrix,  $$
then
$$ Ad_{g} = \pmatrix EH+FG & HG & -FE \\
2HF & H^{2} & -F^{2} \\ -2GE & -G^{2} & E^{2} \endpmatrix \tag 1.6.9 $$
in some fixed basis for the Lie algebra.
Then from the fact that $\w \thicksim \theta$, it is easy to
see that $\w$ bounds $Ad$.
\subheading{Example 1.6.10}
Let $G=GL(n, \R)$, and let $\w$ be the word weight  on $G$.
Define $\theta(g) = \max(\pa g\pa,
\pa g^{-1}\pa)$.  We show that $\w \thicksim \theta$.
Let $H$ be the subgroup of $G$
of matrices with determinant $\pm 1$.
We have a natural isomorphism of groups
$$ \R \times H \cong G, \tag 1.6.11 $$
given by $(r, [a_{ij}]) \mapsto [ e^{r}a_{ij}]$.
Since $H/SL(n, \R) \cong \Z_{2}$,
an argument similar to the one used in
Example 1.6.4 shows that
$\theta \restriction_{H} $ is equivalent to the word weight $e^{\tau_{H}}$
on $H$.
Similarly, $\theta \restriction_{\R}$  gives
the weight $e^{|r|}$ on $\R$ (in the decomposition (1.6.11)),
which is easily seen to be equivalent to the
word weight $e^{\tau_{\R}}$ on $\R$.
Note that
$$ \theta(e^{r}[a_{ij}]) = \theta\restriction_{\R}
(r) \quad\theta\restriction_{H}([a_{ij}]). $$
It follows that $\theta$ dominates the product $e^{\tau_{\R} + \tau_{H}}$.
So, to see that $\theta$ dominates $\w$,
it suffices to show that the gauge on $G$ defined
by
$$ \tau(e^{r} [a_{ij}]) = \tau_{\R}(r) + \tau_{H}([a_{ij}]).\tag 1.6.12 $$
strongly dominates the word gauge.
\par
Let $U$ and $V$ be respective generating sets for $\R $ and $H$.
Then $W = U \times V$ is a generating set for $G$.  If
$(e^{r} [a_{ij}])\in W^{k}$, then $r \in U^{k} $ and $[a_{ij}] \in
V^{k}$.  Hence $\tau_{G}(e^{r} [a_{ij}]) \leq
\max(\tau_{\R}(r), \tau_{H}([a_{ij}]))
\leq \tau_{\R}(r) + \tau_{H}([a_{ij}])$.
So we have shown that $\tau$ defined by (1.6.12) strongly
dominates the word gauge $\tau_{G}$.  It follows that $e^{\tau}$
dominates $e^{\tau_{G}}$, and hence that $\theta$ dominates
the word weight $\w=e^{\tau_{G}}$.  Since $G$ is compactly
generated, the word weight dominates $\theta$, and we have
$\theta \thicksim \w$.
\vskip\baselineskip
\heading  \S 2 Smooth Crossed Products
\endheading
\par
In this section, $A$ will be a Fr\'echet algebra with
a  continuous action of a Lie group $G$.
We ignore any *-structure on $A$, and we do not assume
that the action of $G$ on $A$ is differentiable or strongly continuous.
If $\s$ is a scale on $G$,
then we shall define the $(L^{1})$ Schwartz functions
$\sonaco$ from $G$
to $A$.
If the action of $G$ on $A$
is tempered in an appropriate sense, and $\s$ satisfies certain
simple properties,
we will see that  $\sonaco$ is a Fr\'echet algebra under convolution,
which we will denote by $\sonas$.
\par
To simplify notation, we shall make the assumption that
$\s \geq 1$ throughout \S 2.  We lose no generality, since any
scale $\s$ is equivalent to either of the scales $\max(1,\s)$
or $1 + \s$.
\heading \S 2.1 Fr\'echet Space Valued Schwartz Functions
\endheading
\subheading{Definition 2.1.0} Let
$\s\geq 1$ be any scale on a Lie group $G$,
and let $E$ be any Fr\'echet space.  Assume $\s$
is bounded on compact subsets of $G$. Let $\pa \quad \pa_{m}$
be a family of increasing seminorms giving the topology of $E$.
 We define the $\s$-rapidly vanishing
$L^{1}$ Schwartz functions $\soneco$ from
$G$ to $E$ to be the set of differentiable functions
$\psi \colon G \rightarrow E$ such that the seminorms
$$ \pa \psi \pa_{d,\ga, m} =
\int_{G} \s^{d}(g) \pa X^{\ga} \psi(g) \pa_{m} dg  \tag 2.1.1
$$
are finite for each $d, m \in \N$ and $\ga \in \N^{q}$.
(Here $X^{\ga}\psi(g)$
is defined by formula (1.2.1) with action $\be_{h}(\psi)(g) =
\psi(h^{-1}g)$.)
We topologize $\soneco$ by these seminorms.
To simplify notation, we shall refer to $\soneco$ as $\sonec$. Note that
$C_{c}^{\infty}(G, E) \subseteq \sonec$.
\par
Before giving conditions for $\sonec$ to be a Fr\'echet space,
we define another set of functions from $G$ to $E$, the set of
$C^{\infty}$-vectors of which will be precisely $\sonec$.
(For this discussion, $G$ can be any locally compact group.)
 Let $\F_{d}$ denote the space of all functions
$\psi \colon G \rightarrow E$ for which
$$ \pa \psi(g) \pa_{d, m}=
\int_{G}\s^{d}(g) \pa \psi(g) \pa_{m}  dg
< \infty \tag 2.1.2 $$
for all $m$, where $\int$ denotes the upper integral.
(The upper integral of a positive function $f\colon G \rightarrow
\R $ is the infimum
of the integrals of the countably infinite linear combinations $\sum c_{n}
\xi_{n}$,
where $c_{n}\geq 0$, the functions $\xi_{n}$ are  characteristic functions
of integrable subsets of $G$
(where integrable means Borel measurable with
finite Haar measure), and $\sum c_{n} \xi_{n}$ dominates $f$
pointwise almost everywhere \cite{Tr, p. 468, p. 99}.)
For each $m\in \N$, the expression (2.1.2)
then defines a seminorm on $\F_{d}$.
Define a step function $s$ from $G$ to
$E$ to be a finite linear combination
of the form $\sum e_{n} \xi_{n}$, where $e_{n}\in E$ and $\xi_{n}$
is the characteristic function of an integrable subset of $G$.
Let $\Cal L_{d}$ denote the closure of step functions
in $\F_{d}$ in the topology given by the seminorms (2.1.2).
Since $\s^{d+1} \geq \s^{d}$,
the seminorm $ \pa \quad \pa_{d+1, m} $ dominates
$\pa \quad \pa_{d, m}$ for each $m$, and
we have natural canonical
continuous injections $\Cal L_{d+1} \hookrightarrow \Cal L_{d}$.
\par
Let $N_{d}$ denote the set of functions on which the seminorms
(2.1.2) vanish for all $m$.
This is clearly a closed subspace of $\Cal L_{d}$.
If $f\in N_{d}$,
then $\int \pa f(g) \pa_{m} \s^{d} (g) dg =0 $
for every $m$.  Hence, since $\s\geq 1$, the function $f$ vanishes
almost everywhere.  So for any $k\in \N$, we have
$\int \pa f(g) \pa_{m} \s^{k} (g) dg =0 $
for every $m$, and $f\in N_{k}$.
It follows that the sets $N_{d}$ are independent of $d$.
\par
We let $L_{d}$
denote the quotient $\Cal L_{d}/N_{d}$.
Every element
of $L_{d}$ may then be thought of as an equivalence class of measurable
functions.
By the Fisher-Riesz theorem for Fr\'echet spaces, the
space $L_{d}$ is complete for the topology given by the seminorms
(2.1.2)
\cite{Tr, p. 468}.
Since $N_{d+1}=N_{d}$, the map $\Cal L_{d+1} \hookrightarrow \Cal
L_{d}$ gives a continuous
injection
$$ {L}_{d+1} \hookrightarrow {L}_{d}. $$
\par
We define the {\it $\s$-rapidly vanishing $L^{1}$ functions from
$G$ to $E$} to be the intersection (or projective limit)
$$L_{1}^{\s}(G, E) = \cap_{d=0}^{\infty}L_{d}.\tag 2.1.3$$
The topology on $L_{1}^{\s}(G, E)$ is easily seen to be given
by the seminorms (2.1.2), where $d$ ranges over all natural numbers.
It is also easy to see from the definition (2.1.3) and the completeness of
each $L_{d}$  that $L_{1}^{\s}(G, E)$
is complete for this topology.
To simplify notation, we will denote  $L_{1}^{\s}(G, E)$ by $L^{\s}(G, E)$.
 We could have defined $L^{\s}(G,E)$
to be the closure of integrable step functions
or of continuous functions with compact support in an appropriate space
$\F_{\infty}$, and then moded out by the set of functions which lie
in the kernels of all
the seminorms (2.1.2) \cite{Schw, \S 5}.  An easy argument
shows that we would have gotten the same space.  Note that
$L^{\s}(G, E)\hookrightarrow L^{\tau}(G, F)$ with
continuous inclusion map if $E\hookrightarrow F$ with continuous
inclusion, and $\s$ dominates $\tau$.
\par
The action
$$(gF)(h) = F(g^{-1}h), \qquad g, h \in G\tag 2.1.4$$
of $G$ on $L^{\s}(G, E)$ is  strongly
continuous if $\s$ is uniformly translationally equivalent (see (1.2.10)).
To see this, let $C_{c}(G, E)$ denote
the continuous functions with compact support from $G$ to $E$.
Then (2.1.4) gives an action of $G$ on $C_{c}(G, E)$ which is easily
seen to be strongly continuous for the inductive limit
topology.  Using the uniform translational equivalence of $\s$,
the density of $C_{c}(G, E)$
in $L^{\s}(G, E)$,  and the fact that the seminorms (2.1.2) are continuous
for the inductive limit topology on $C_{c}(G, E)$, a quick calculation
shows that $G$ acts strongly continuously on $L^{\s}(G, E)$.
\proclaim {Theorem 2.1.5} Let $\s$ be a translationally equivalent
scale on a Lie group $G$ (see Definition 1.2.7).
Then the space $\Cal S_{1}^{\s}(G, E)$
is the set of $C^{\infty}$-vectors
for the action (2.1.4) of $G$ on $L_{1}^{\s}(G, E)$.  Hence
$\Cal S_{1}^{\s}(G, E)$
is complete and a Fr\'echet space.   The
set of compactly supported
functions $C_{c}^{\infty}(G, E)$ is dense in $\soneco$.
If $\s_{1} \thicksim \s_{2}$, then
$\Cal S_{1}^{\s_{1}}(G, E) = \Cal S_{1}^{\s_{2}} (G, E)$.
\endproclaim
\demo{Proof} Recall from Theorem 1.2.11 that
$\s$ is uniformly translationally equivalent
since it is translationally equivalent.
Let $\Cal S$ denote $\sonec$,  let $L$ denote
$L_{1}^{\s}(G, E)$, and let $L^{\infty}$
denote the set of $C^{\infty}$-vectors
for the action of $G$ on $L$.
A calculation similar to (1.2.9) shows that $G$ acts strongly
continuously and differentiably on $\Cal S$.
\par
We show
that $L^{\infty}\subseteq S$.
By convention, $\int \varphi(h) dh = \int \Delta(h) \varphi(h^{-1}) dh$,
where $\Delta$ is the modular function.
By \cite{DM, Thm 3.3}, we may write any element of $L^{\infty}$
as a finite sum of functions
$$ \split { f} *\psi(g) = \int_{G} { f}(h)\psi (h^{-1}g)  dh &
=\int_{G} f(gh) \psi(h^{-1})dh   \\
=\int_{G} \Delta(h) f(g h^{-1})& \psi(h) dh = \Delta(g) \int_{G}
\Delta^{-1} f(g h^{-1}) \psi(h) dh,
\endsplit
\tag 2.1.6$$
where $f\in \cc$, $\psi \in L^{\infty}$.
It suffices to prove that functions of the form
(2.1.6) lie in $\Cal S$.  We first show that such functions
are continuous.  Let $g_{n}\longrightarrow g_{0}$ in $G$,
and let ${\tilde f} = \Delta^{-1} f$,
${\tilde f}_{g}(h)={\tilde f}(g^{-1}h)$.
Since $\Delta$ is continuous, it suffices to show that $g\mapsto
\Delta^{-1}(g) (f * \psi)(g)$ is continuous.
For any seminorm $\pa \quad \pa_{d}$
on $E$,
$$ \aligned \pa
\Delta^{-1}({f}*\psi)(g_{n}) - \Delta^{-1}({f}*\psi) (g_{0}) \pa_{d}
&\leq
\int \pa {\tilde f}(g_{n}h^{-1} ) - {\tilde f}(g_{0}h^{-1})
\psi(h) \pa_{d}
 dh \\
& \leq \int \pa \psi (h) \pa_{d} dh
\pa  {\tilde f}_{g_{n}^{-1}} - {\tilde f}_{g_{0}^{-1}}
\pa_{\infty} \\
& \leq C_{f}\pa \psi \pa_{0, d}
\pa {\tilde f}_{g_{n}^{-1}} -
{\tilde  f}_{g_{0}^{-1}}\pa_{\infty}. \endaligned $$
The last expression tends to zero
as $g_{n}$ tends to $g_{0}$, since ${\tilde f}\in \cc$.  Thus ${ f}*\psi$
is a continuous function.
Similar arguments, bringing derivatives inside the integral (2.1.6)
to act on $f$, show that ${ f}*\psi$ is  a differentiable
function from $G$ to $E$.
Thus $L^{\infty}$ consists of differentiable functions from $G$ to $E$.
(An argument similar to \cite{Co-Gr, Lemma 2.3.2} may also work for this.)
Since the seminorms on $L^{\infty}$ are precisely those on $\Cal S$,
we have $L^{\infty}\subseteq \Cal S$.
\par
To see that $C_{c}^{\infty}(G, E)$ is
dense in $\Cal S$, we apply the following
lemma and the fact that $C_{c}^{\infty}(G, E)$ is dense in $L$.
(Alternatively, one could probably use approximate units.)
\proclaim{Lemma 2.1.7}  Let $E$ be a Fr\'echet space on which a Lie
group $G$ acts
strongly continuously.
Assume that $F$ is a dense $G$-invariant subspace of $E$.
Then the algebraic span $C_{c}^{\infty}(G)F$ is dense in $E^{\infty}$.
\endproclaim
\demo{Proof}
By \cite{DM, Thm 3.3}, every element of $E^{\infty}$ is a finite
sum of elements of the form
$\alpha_{\varphi} (e) =\int \varphi(h) \alpha_{h}(e) dh$, where
$e \in E$ and $\varphi \in C_{c}^{\infty}(G)$.
Let $f_{n} \rightarrow e$ in $E$.  Then $\alpha_{\varphi}(f_{n})$
tends to $\alpha_{\varphi}(e)$ in $E$.
\qed
\enddemo
\par
The last statement about equivalent scales is a simple calculation.
This proves Theorem 2.1.5.
\qed
\enddemo
\heading \S 2.2 Conditions for an Algebra
\endheading
We give  conditions on the scale $\s$ and
the action of $G$ on a Fr\'echet algebra
$A$ which makes
$\sonac$ an algebra under convolution.
We begin with a general lemma for actions of
$G$ on Fr\'echet spaces.
\proclaim{Lemma 2.2.1} Assume that the action $\alpha$ of
$G$ on a Fr\'echet space $E$ is differentiable. Let $e\in E$.
Let $\ga \in \N^{q}$,  where $q$ is the dimension
of the Lie algebra $\frak G$.
Assume that $\s\geq 1$ and $\s_{-}$
bounds $Ad$ (recall $\s_{-}(g)= \s(g^{-1})$).
Then there exists a constant $D>0$ and $p\in \N$ such that
$$\pa {T}X^{\ga}(\alpha_{h}(e))\pa_{d}\leq
D\s^{p}(h)  \sum_{\be\leq \ga} \pa
{T}\alpha_{h}(X^{\be}e)\pa_{d},
\tag 2.2.2
$$
for all $d\in \N$, $ h \in G$, $e\in E$,
and  all continuous linear maps $T\colon E \rightarrow E$.
\endproclaim
\par
Here $X^{\ga}\psi$ is defined by  formula (1.2.1)
 with $\be= \alpha$.
\demo{Proof} Consider the case $\ga_{i}=(0,\dots 1,\dots 0)$,
where the $1$ is in the $i$th spot.
By the chain rule,
$$ \aligned
X^{\ga_{i}} \alpha_{h}(e) & = {d\over dt}\alpha_{exp(tX_{i})h}
(e)\restriction_{t=0}
\quad = {d\over dt} \alpha_{h exp(Ad_{h^{-1}}(tX_{i}))}
(e)\restriction_{t=0}\\
&= \sum_{j=1}^{q} \alpha_{h}(X^{\ga_{j}}e) {d\over dt}
(Ad_{h^{-1}}(tX_{i}))^{(j)}|_{t=0}
= \sum_{j=1}^{q} \alpha_{h}(X^{\ga_{j}}e)
(Ad_{h^{-1}})_{ji}.\\
\endaligned  \tag 2.2.3
$$
(Here if $Y\in \frak G$, then the $Y^{(j)}$'s
are the coordinates of $Y$ in the basis $\basisq$.)
Since $\s_{-}$ bounds $Ad$, the matrix elements $(Ad_{h^{-1}})_{ji}$
have their absolute values bounded by $D\s^{p}(h)$ for some constant
$D$ and $p\in \N$.
Thus
$$
\pa TX^{\ga_{i}}(\alpha_{h}(e))\pa_{d}\leq
D\s^{p}(h)  \sum_{j=1}^{q} \pa
T\alpha_{h}(X^{\ga_{j}}e)\pa_{d}.
$$
The case of general $\ga \in \N^{q}$
follows by repeated application of $X^{\ga_{i}}$'s
to Eqn 2.2.3,
and again using the bound on $|(Ad_{h^{-1}} )_{ji}|$
by $D\s^{p}(h)$.
This proves Lemma 2.2.1.\qed \enddemo
\vskip\baselineskip
\subheading{Definition 2.2.4}
Let $\pa \quad \pa_{m}$ be an increasing
sequence of seminorms  on
a Fr\'echet algebra $A$ giving the topology of $A$.
We let $\text{Aut}(A)$ denote the group of continuous
algebra automorphisms of $A$.  We say that a group
homomorphism $\alpha \colon G \rightarrow \text{Aut}(A)$
is a {\it continuous action} of $G$ on $A$
if
 $g \mapsto \pa \alpha_{g}(a)\pa_{m}$
is continuous for each $m \in \N$ and $a \in A$. From
now on, we shall only consider continuous actions of $G$.
Let $\s\geq 1$ be a scale on $G$.
We say that the action of $G$ on $A$ is
{\it $\s$-tempered} if
for every $m$ there are $k, l$, and $C$  so that
$$ \pa \alpha_{g}(a) \pa_{m} \leq
C\s^{k}(g) \pa a\pa_{l}, \qquad a \in A, \, g \in G. \tag 2.2.5 $$
If the scale is understood, we simply say that the action
is tempered.
This definition is analogous to the definition of
smoothness in \cite{ENN} for an action of
$\R$ (without the
differentiability condition), and the notion of a tempered
$G$-module for a nilpotent Lie group in \cite{DuC 1, \S 4}, \cite{DuC 3}
\cite{DuC 2}.
\par
We remark that tempered does not imply strongly continuous.
For example, let $\R$ act by translation on the bounded continuous
functions $C_{b}(\R)$ with sup norm topology.
\proclaim {Theorem 2.2.6}
Let $\s$ be a sub-polynomial scale on a locally compact group $G$.
Let $\alpha$ be a $\s$-tempered
action of $G$ on a Fr\'echet algebra $A$.
Then the Fr\'echet space
$L_{1}^{\s}(G, A)$ is a Fr\'echet algebra under convolution.
\par
Assume in addition that $G$ is a Lie group, and that
either $\s_{-}$ bounds $Ad$ or that $G$ acts differentiably on $A$.
Then
the Fr\'echet space
$\sonaco$ is a  Fr\'echet algebra
under convolution.
\endproclaim
\par
When the conditions of Theorem 2.2.6 are satisfied, we define
the {\it smooth crossed product of $A$ by $(G, \s )$}
to be the Fr\'echet space $\sonaco$ with
convolution multiplication.
We denote it by $\sonas$, or $\sona$ if the scale $\s$ is understood.
\demo{Proof}
We check that  $\sonac$ is closed
under multiplication.  Without loss of generality, we assume $\s\geq 1$.
Assume that $\s_{-}$ bounds $Ad$, so that Lemma 2.2.1 applies.
We compute the seminorm of a product in $\sonac$.
$$\aligned & \pa \psi *  \varphi \pa_{m, \ga , d} =
\int_{G} \s^{m}(g) \pa X^{\ga} (\psi *\varphi )(g) \pa_{d} dg
\\
& =
 \int  \s^{m}(g) \pa\int \psi(h)
\alpha_{h}(( X^{\ga} \varphi_{h} )( g))dh \pa_{d}
dg,\quad {\text{def of conv }}\\
&\leq \int \int \s^{m}(g) \pa \psi(h) \pa_{d^{\prime}}
\pa \alpha_{h}((
X^{\ga} \varphi_{h}) (g)) \pa_{d^{\prime}}
dh dg ,\quad {\text{$A_{}$ Fr\'echet }}\\
&\leq \int \int \s^{m}(g) \s^{k}(h)C \pa \psi(h) \pa_{d^{\prime}}
\pa (X^{\ga} \varphi_{h}) (g) \pa_{l}
dh dg ,\quad {\text{tempered action   }}\\
& \leq \sum_{\be\leq \ga} \int \int \s^{m}(g)
\pa \s^{k} \psi(h) \pa_{d^{\prime}}
CD\s^{p}(h) \pa (X^{\be} \varphi )_{h}(g)\pa_{l}
dhdg , \quad {\text{  (2.2.2)}}\\
& \leq {\tilde C} \sum_{\be\leq \ga}
\int \int
\pa \s^{k+p+mr}\psi(h) \pa_{d^{\prime}}
\pa \s^{mr}(X^{\be} \varphi )(h^{-1}g)\pa_{l}
dgdh ,  {\text{   sub-poly}}\\
& = {\tilde C} \sum_{\be\leq \ga}
\int \int
\pa \s^{k+p+mr}\psi(h) \pa_{d^{\prime}}
 \pa \s^{mr}(X^{\be} \varphi )(g)\pa_{l}
dgdh , \quad {\text{ }}\\
& = {\tilde C} \sum_{\be\leq \ga}
 \pa \psi \pa_{k+p+mr, 0, d^{\prime}}
\pa \varphi \pa_{mr, \be, l} \qquad
{\text{ def of norms.}} \endaligned
\tag 2.2.7 $$
So  $ \sonac$  is a Fr\'echet
algebra.
The computation (2.2.7) without the derivatives shows that $L^{\s}(G, A)$
is a Fr\'echet algebra (without assuming that $\s_{-}$
bounds $Ad$).  If $G$ acts
differentiably on $A$, the estimates \cite{Sc 3, (6.23),(6.24)}
show that $\sonac$ is a Fr\'echet algebra
under convolution. (Alternatively,
one could do an appropriate change of variables for $h$ in the
early stages of (2.2.7).)
This proves Theorem 2.2.6.
\qed
\enddemo
\par
See \S 5 for examples.
We prove some elementary results on ideals and quotients,
which we will be using in a later paper.
\proclaim{Proposition 2.2.8} Let $G$, $\s$ and $A$
be as in Theorem 2.2.6.   Let $I$ be any closed $G$-invariant
two-sided ideal of $A$.
Let $J$ denote the closed subspace of $\sona$ of functions which
take their values in $I$.  Then $J$ is a two-sided ideal in $\sona$
and $\sona/J \cong G\rtimes (A/I)$.
\endproclaim
\demo{Proof} First note that if the action of $G$ on $A$ is tempered,
then the action of $G$ on $A/I$ is tempered.  Hence we may form
the crossed product $G\rtimes (A/I)$.
Let $\pi \colon \sona \rightarrow G\rtimes(A/I)$ be the canonical map.
We show that
$\pi $ is onto.  The map $L^{\s}(G, A) \rightarrow L^{\s}(G, A/I)$
is onto since the projective tensor product of surjective maps
is surjective\cite{Schw, \S 5}\cite{Tr, Prop 43.9}.
By Theorem A.8, the corresponding map of $C^{\infty}$-vectors
for the action of $G$ be left translation is surjective, so by
Theorem 2.1.5, $\pi$ is surjective.
\par
Clearly the kernel of $\pi$ is $J$.  \qed \enddemo
\vskip\baselineskip
\heading  \S 3 $m$-convexity  \endheading
\heading \S 3.1 Conditions for
$m$-convexity of the Crossed Product \endheading
\par
Assume that $A$ is an $m$-convex Fr\'echet algebra.
That is, assume that the seminorms topologizing
$A$ can be taken to be submultiplicative.
We give conditions so that the convolution algebra
 $\sona$ is $m$-convex.
\subheading{Definition 3.1.1} Let $\s\geq 1$ be
a scale on a topological group $G$.
We say that the action of $G$ on $A$ is
{\it $m$-$\s$-tempered } if there exists a family of submultiplicative
seminorms $\bigl\{ \normm\bigr\}$ topologizing $A$
such that for every $m \in \N$, there exists $C>0$ and
$d \in \N$, for which
$$
\,  \pa \alpha_{g}(a) \pa_{m} \leq
C\s^{d}(g) \pa a \pa_{m}, \qquad a\in A,
g\in G.
\tag 3.1.2 $$
Note that this is similar to the
definition of a $\s$-tempered action (2.2.4),
except that we have submultiplicative seminorms, and the
same seminorm appears both on
the left and the right of the inequality (3.1.2).
We call the seminorms $\normm$   {\it ($G$, $\s$)-stable
submultiplicative } seminorms for $A$.
If the scale $\s$ is understood,
we simply say that the action is $m$-tempered.
If the action of $G$ is $m$-tempered, it is clearly tempered.
This more restrictive notion of an $m$-tempered
action
seems to be necessary in the proof of Theorem 3.1.7 below, to obtain the
$m$-convexity of the crossed product.
\subheading{Question 3.1.3} Is a
tempered action of a group on an $m$-convex
Fr\'echet algebra always $m$-tempered ?
\par
The reverse direction of the following theorem follows from
\cite{MRZ, Lemma 2}.  The forward direction is an easy exercise.
See \cite{Pe, Thm 3.2} for a similar result for matrix algebras.
\proclaim {Theorem 3.1.4} A Fr\'echet algebra $\Cal A$ is $m$-convex
if and only if for every  family (or equivalently for any one family)
 $\bigl\{ \normm\bigr\}$ of increasing seminorms giving
the topology of $\Cal A$,
and for every $p \in \N$, there exists $C>0$, $r\geq p$
such that
$$
\pa a_{1}\dots
a_{n}\pa_{p} \leq C^{n} \pa a_{1} \pa_{r}
\dots \pa a_{n} \pa_{r},  \tag 3.1.5 $$
for all $n$-tuples $a_{1}\dots a_{n}$ of elements of $A$,
and all $n\in \N$.
\qed
\endproclaim
\vskip\baselineskip
\par
Let $\Cal A_{m}$ be a family of $m$-convex
Fr\'echet algebras with continuous algebra
homomorphisms
$\pi_{mk}\colon \Cal A_{m} \rightarrow \Cal A_{k}$
for $m \geq k$.
Assume that $\pi_{mk} \circ \pi_{lm} = \pi_{lk}$
for $l\leq m \leq k$, $\pi_{mm}=\text{id}$, and that
$\Cal A_{m+1}$ has dense image in $\Cal A_{m}$.
Let $\prod_{m=0}^{\infty} \Cal A_{m}$
be the topological cartesian product
of the spaces $\Cal A_{m}$.  We call the
subset of $\prod_{m=0}^{\infty}
\Cal A_{m}$ of tuples $(a_{0}, a_{1}, \dots)$
which satisfy $\pi_{m+1m}(a_{m+1})= a_{m}$ for all $m$
the {\it projective limit} of the
$\Cal A_{m}$'s, denoted by $\Cal A=\varprojlim \Cal A_{m}$.
The image of $\Cal A$ is dense in
each $\Cal A_{m}$, and is complete for the
strongest topology which
dominates the topology on each $\Cal A_{m}$
(or, in other words
$\Cal A$ is a Fr\'echet space for the topology given
by all the seminorms on
all the $\Cal A_{m}$'s) \cite{Mi} \cite{Ph, 1.4-5}.
Since we may take every seminorm submultiplicative, $\Cal A$ is
also an $m$-convex Fr\'echet algebra.
Similarly, if $\Cal A$ is
any $m$-convex Fr\'echet algebra, and  $\pa \quad \pa_{m}$
is a family of increasing submultiplicative seminorms giving the
topology for $\Cal A$, then  $\Cal A$ is the projective limit
$\varprojlim \Cal A_{m}$ of the Banach algebra
completions
$\Cal A_{m}$ in $\pa \quad \pa_{m}$\cite{Mi}.
We say that a scale $\s\geq 1$ is {\it $m$-sub-polynomial} if there
exists $C>1$ and $l \in \N$ such that
$$
\s(g_{1} \dots g_{n}) \leq C^{n} \s^{l}(g_{1}) \dots \s^{l}(g_{n}),
\tag 3.1.6
 $$
for all $g_{1}, \dots g_{n} \in G$ and $n \in \N$.
If $\s_{1} \thicksim \s_{2}$, then $\s_{1}$ is $m$-sub-polynomial iff
$\s_{2}$ is.  Any scale which is equivalent to a weight or a gauge
is $m$-sub-polynomial.
Recall that $\s_{-}(g)= \s(g^{-1})$.
\proclaim{Theorem 3.1.7}
Let $\s$ be an $m$-sub-polynomial scale
on a Lie group $G$ such that
either $\s_{-}$ bounds $Ad$ (for example, $\s$ could be
any weight or  gauge that bounds $Ad$) or $G$ acts differentiably
on $A$,
and let the action of $G$ on an $m$-convex
Fr\'echet algebra $A$ be $m$-$\s$-tempered. Then
the convolution algebra $\sonas$ is an $m$-convex Fr\'echet
algebra.  In fact $\sonas$ is the projective limit of
Fr\'echet algebras $G\rtimes^{\s}A_{m}$,
where $A_{m}$ is the completion of $A$ in the $m$th ($G$, $\s$)-stable
submultiplicative seminorm
$\pa \quad \pa_{m}$.
\par
Similarly, if $G$ is any locally compact group,
then the corresponding conclusions are true for
the Fr\'echet algebra $L_{1}^{\s}(G, A)$, using only the assumptions
that the action of $G$ on $A$ is $m$-$\s$-tempered,
and  that $\s$ is $m$-sub-polynomial.
\endproclaim
\demo{Proof}
First consider the case where $A$ is a Banach algebra.
Without loss of generality, assume $\s \geq 1$. Let $C_{1}>1$ and
$l \in \N$ be such that
$$ \s(g_{1} \dots g_{n}) \leq C_{1}^{n} \s^{l}(g_{1}) \dots \s^{l}(g_{n}),
\qquad g_{1}, \dots g_{n} \in G. $$
Let $\pa \quad \pa$ be a $(G, \s)$-stable submultiplicative
seminorm giving the topology for $A$,
and let $C>1$ and $m \in \N$ be such that $\pa \alpha_{g} (a) \pa \leq
C \s^{m}(g) \pa a \pa$.
We verify the $m$-convexity by showing that the
family of increasing seminorms
$$\pa \psi \pa^{\prime}_{p}
=\max_{|\ga|\leq p}
\pa \psi \pa_{p,\ga}= \max_{|\ga|\leq  p}
\int \pa \s^{p}X^{\ga}\psi(g) \pa dg
\tag 3.1.8$$
satisfies (3.1.5).
\par
To prepare to
estimate  $\pa \psi_{1}*\dots \psi_{n} \pa_{d, \ga}$,
we write $\psi_{1} *\dots \psi_{n} (g)$ as
$$ \int \dots \int \alpha_{\eta_{1}}(\psi_{1}(h_{1}))
\dots \alpha_{\eta_{n-1}}(\psi_{n-1}(h_{n-1}))
\alpha_{\eta_{n}}(\psi_{n}(h_{n})) dh_{1}\dots dh_{n-1} \tag 3.1.9 $$
where $h_{1}, \dots, h_{n-1}$ are the variables of integration,
$h_{n} = h^{-1}_{n-1} \dots h_{1}^{-1} g$, $\eta_{1}= e$,
$\eta_{k}= h_{1}\dots h_{k-1}$.
We proceed to estimate.  Using (3.1.9) and the left invariance of
Haar measure,
$$\aligned
& \pa \psi_{1}*\dots \psi_{n} \pa_{d, \ga} =
\int \s^{d}(g)
\pa(X^{\ga}( \psi_{1}\dots \psi_{n}))(g) \pa dg\\
&\leq
 \int \dots \int \s^{d}(g) \pa \alpha_{\eta_{1}}(\psi_{1}(h_{1}))
\\ & \qquad
\dots         \alpha_{\eta_{n-1}} (\psi_{n-1}(h_{n-1}))
\alpha_{\eta_{n}}
(X^{\ga}
(\psi_{n})_{\eta_{n}}(g))\pa dh_{1}\dots dh_{n-1}
dh_{n}.
\endaligned \tag 3.1.10 $$
By Lemma 2.2.1, for some $j,C_{2}>1$ depending only on our $q$-tuple
$\ga$,
the integrand in (3.1.10) is bounded by
$$C_{2} \s^{d}(g)\s^{j}(\eta_{n}) \sum_{\be\leq \ga}
\pa \alpha_{\eta_{1}}(\psi_{1}(h_{1}))
\dots
\alpha_{\eta_{n-1}}(\psi_{n-1}(h_{n-1}))
\alpha_{\eta_{n}}(X^{\be}\psi_{n})(h_{n})\pa.
\tag 3.1.11 $$
Since the action of $G$ is  $m$-$\s$-tempered,
we may bound the summand in
(3.1.11).
$$\aligned
\pa \alpha_{\eta_{1}}( & \psi_{1}(h_{1}))
\dots
\alpha_{\eta_{n}}(X^{\be}\psi_{n})(h_{n})\pa \\
&
\leq
\pa \psi_{1}(h_{1}) \pa \,\pa \alpha_{h_{1}}(\psi_{2})(h_{2})\dots
\alpha_{\eta_{n}}
(X^{\be}\psi_{n})(h_{n}) \pa \quad \text{since $\eta_{1}=e$,
$\eta_{2}=h_{1}$} \\
&\leq
\pa \psi_{1}(h_{1}) \pa\, C\s^{m}(h_{1})\,\pa \psi_{2}(h_{2})\pa
\,\pa \alpha_{h_{1}^{-1}\eta_{2}}(\psi_{3})(h_{3}) \dots
\alpha_{h_{1}^{-1}\eta_{n}}
(X^{\be}\psi_{n})(h_{n}) \pa \\
& \dots \qquad \dots \qquad \dots \\
&\leq C^{n-1}\s^{m}(h_{1})
\dots \s^{m}(h_{n-1})
\,\pa \psi_{1}(h_{1}) \pa\,\pa \psi_{2}(h_{2})\pa
\dots
\pa
X^{\be}\psi_{n}(h_{n}) \pa \endaligned
\tag 3.1.12 $$
where we have repeatedly used the $(G, \s)$-stable submultiplicative
property of our seminorm $\pa \quad \pa$.
Since $\s$ is $m$-sub-polynomial,
$$ \s^{j}(\eta_{n}) \leq C_{1}^{jn}\s^{jl}(h_{1})\dots \s^{jl}(h_{n-1})
\tag 3.1.13
$$
and
$$\s^{d}(g)\leq C_{1}^{dn}\s^{dl}(h_{1})\dots \s^{dl}(h_{n}).\tag 3.1.14
$$
Now we are ready to plug everything back in and make our
estimate.
Plugging
(3.1.12) and (3.1.13) into (3.1.11), we see that
the integrand of (3.1.10) is bounded by
$$\aligned \s^{d}(g)
\sum_{\be\leq \ga}
C_{2}(CC_{1}^{j} )^{n}
& \pa \s^{jl+m}\psi_{1}(h_{1})\pa \\
& \dots
\pa \s^{jl+m} \psi_{n-1}(h_{n-1})
\pa
\pa (X^{\be}\psi_{n})(h_{n})\pa
.\endaligned \tag 3.1.15
$$
Using (3.1.14), we see that
(3.1.15) is bounded by
$$
\sum_{\be\leq \ga} C_{2}(CC_{1}^{j+d} )^{n}
\pa \s^{jl+m+dl} \psi_{1}(h_{1})\pa
\dots
\pa \s^{dl}(X^{\be}\psi_{n})(h_{n})\pa. \tag 3.1.16 $$
Let $t=jl+m+dl$.
Then plugging (3.1.16) back into the integral (3.1.10), we see that
$$\aligned
& \pa \psi_{1}\dots \psi_{n} \pa_{d, \ga}\\
&
\leq \sum_{\be\leq \ga} C_{2}(CC_{1}^{j+d} )^{n}
\pa \psi_{1}\pa_{t,0}
\dots
\pa \psi_{n-1} \pa_{t,0}
\pa \psi_{n}\pa_{t,\be}\\
&  \leq C_{3}^{n}
\pa \psi_{1}\pa_{p}^{\prime}
\dots
\pa \psi_{n-1} \pa_{p}^{\prime}
\pa \psi_{n}\pa_{p}^{\prime}
\endaligned\tag 3.1.17$$
where $p\geq t,  |\ga|$, and $C_{3}$ is a sufficiently
large constant.  By Theorem 3.1.4, this
proves the $m$-convexity of $G\rtimes^{\s}A$.
\par
Now let $A$ be any $m$-convex Fr\'echet algebra, and
let $\pa \quad \pa_{m}$ be a family of increasing $(G, \s)$-stable
submultiplicative seminorms for the topology of $A$.  Let $A_{m}$
denote the completion of
$A$ in $\pa \quad \pa_{m}$.  Then $G\rtimes^{\s}A_{m}$
is $m$-convex by what we have just shown.  Since $\sonas$ is precisely
the set of differentiable functions for which each of the seminorms for
$G\rtimes^{\s}A_{m}$ is finite for each $m$ (see  (2.1.1)),
$\sonas$ can be identified with the projective limit
$\varprojlim G\rtimes^{\s} A_{m}$, and so $\sonas$ is an $m$-convex
Fr\'echet algebra.
\par
The same proof without derivatives gives the corresponding results
for $L_{1}^{\s}(G, A)$.
\qed
\enddemo
\par
Just as the notion of $m$-convexity has an equivalent condition
in terms of an arbitrary increasing family of
seminorms on $A$ (see Theorem 3.1.4),
so does the notion of an $m$-tempered action. We prove this for the case
when $\s$ is a weight.
The case of a gauge $\s$ may be obtained by applying
the theorem to the weight $1 + \s$.
\proclaim{Theorem 3.1.18} An action of a group $G$ with weight $\s$
on a Fr\'echet algebra
$A$ is $m$-$\s$-tempered if and only if for every
family (or equivalently for any one family)
$\bigl\{ \normm\bigr\}$ of increasing seminorms for $A$ we have
that for every $m \in \N$, there exists $C>0$, $d,r\geq m$ such that
$$ \pa
\alpha_{g_{1}}(a_{1}) \dots  \alpha_{g_{n}}(a_{n})
\pa_{m}
\leq C^{n}\biggl(\s(g_{1})
\s(g_{1}^{-1} g_{2})\dots \s(g_{n-1}^{-1}g_{n}) \s(g_{n})\biggr)^{d}
\pa a_{1} \pa_{r} \dots \pa a_{n} \pa_{r},
\tag 3.1.19 $$
for all $n$-tuples $a_{1}\dots a_{n}\in A$, $g_{1}\dots g_{n}\in G$ and
all $n \in \N$.
\endproclaim
\demo{Proof}
First assume that the action of $G$ on
$A$ is $m$-$\s$-tempered.  Then we have
$$ \aligned \pa \alpha_{g_{1}}&(a_{1}) \dots  \alpha_{g_{n}}(a_{n})
\pa_{m}
\leq
C\s^{d}(g_{1}) \pa a_{1} \pa_{m} \pa \alpha_{g_{1}^{-1}g_{2}}
(a_{2}) \dots  \alpha_{g_{1}^{-1}g_{n}}(a_{n})
\pa_{m}\\
& \leq C^{2} \s^{d}(g_{1}) \s^{d}(g_{1}^{-1}g_{2})
\pa a_{1} \pa_{m} \pa a_{2} \pa_{m} \pa \alpha_{g_{2}^{-1}g_{3}}
(a_{3}) \dots  \alpha_{g_{2}^{-1}g_{n}}(a_{n})
\pa_{m}\\
& \dots \qquad \dots \\
& \leq C^{n} \s^{d}(g_{1})
\s^{d}(g_{1}^{-1}g_{2})\dots \s^{d}(g^{-1}_{n-1}g_{n})
\pa a_{1} \pa_{m} \dots \pa a_{n}
\pa_{m}\endaligned
\tag 3.1.20 $$
for any ($G$, $\s$)-stable
submultiplicative seminorms for $A$.  It follows
easily that if $\normm^{\prime}$ is any equivalent family of seminorms,
then (3.1.19) is satisfied for $\normm^{\prime}$ with  possibly different
constants $C$ and $d$ than were used in (3.1.2). (In fact, it is
easy to check that
if one family of seminorms satisfies (3.1.19), then every increasing
one does.)
\par
To show the converse, let $\normm$ be any increasing family  of
seminorms which satisfy (3.1.19).  For $r \in \N$,
let $m$, $C$ and $d$ be as in (3.1.19).
Let $U_{m}=\, \{\, a \in A \, | \, \pa a \pa_{m}\leq 1 \, \}$.
Then by (3.1.19),
$$ \alpha_{g_{1}}(U_{r}) \dots \alpha_{g_{n}}(U_{r})
\subseteq
C^{n}\biggl(
\s(g_{1})\s(g_{1}^{-1}g_{2})\dots \s^{d}(g^{-1}_{n-1}g_{n})
\s (g_{n}) \biggr)^{d}
U_{m} \tag 3.1.21 $$
for all $n$-tuples $g_{1}, \dots g_{n}\in G$.
For $g_{1},\dots g_{n} \in G$, define the balanced neighborhoods of zero
$$ W_{g_{1},\dots g_{n}}= {1\over
{C^{n}(\s(g_{1})
\s(g_{1}^{-1}g_{2})\dots \s(g^{-1}_{n-1}g_{n})\s(g_{n})
)^{d}}}
\biggl( \alpha_{g_{1}}(U_{r})\dots
\alpha_{g_{n}} (U_{r})\biggr).\tag 3.1.22 $$
Let
$$W_{n}=
\bigcup_{g_{1}, \dots g_{k} \in G, \quad k\leq n}
\biggl( W_{g_{1}, \dots g_{k}} \biggr).
$$
Then $W_{n+1} \supseteq W_{n}$.
Define decreasing convex neighborhoods of zero $V_{r}$ by
$$V_{r} = \bigcup_{n=1}^{\infty} \, conv(W_{n})
\tag 3.1.23
$$
where $conv(S)$ denotes the smallest convex subset of $A$ which contains
a set $S$.
Then $V_{r}\subseteq U_{m}$ so the $V_{r}$
form a basis of zero neighborhoods for the topology of $A$.
Also, the $V_{r}$ are convex and balanced. ($V_{r}$ is convex since
it is the union of increasing convex sets.)
\par
We show that $W_{n}W_{m} \subseteq W_{n+m}$.
Let $k \leq n$, $l \leq m$.  It suffices
to show that
$W_{g_{1},\dots g_{k}}
W_{h_{1}, \dots h_{l}} \subseteq W_{g_{1},
\dots g_{k}, h_{1}, \dots h_{l}}$.
Let $w_{1} =\alpha_{g_{1}}(u_{1})\dots
\alpha_{g_{k}}(u_{k}) /C^{k}(\s(g_{1})
\s(g^{-1}_{1}g_{2})\dots \s(g_{k}) )^{d}\in
W_{g_{1}, \dots g_{k}}$ and
$w_{2}=\alpha_{h_{1}}({\tilde u}_{1})
\dots \alpha_{h_{l}}({\tilde u}_{l})/
C^{p}(\s(h_{1})
\s(h^{-1}_{1}h_{2})\dots \s(h_{l}) )^{d}
 \in
W_{h_{1}, \dots h_{l}}$.
Then
$$ \split &{1\over{C^{k+l}(\s(g_{1})
\s(g_{1}^{-1}g_{2}) \dots \s(g_{k})\s(h_{1})
\s(h^{-1}_{1}h_{2}) \dots \s(h_{l}))^{d}}}\\ &\qquad \qquad \leq
 {1\over{C^{k+l}(\s(g_{1})\s(g_{1}^{-1}g_{2}) \dots \s(g_{k}^{-1}h_{1})
\s(h^{-1}_{1}h_{2})
\dots \s(h_{l}))^{d}}}, \endsplit
 \tag 3.1.24 $$
since $\s(g_{k}^{-1}h_{1}) \leq \s(g_{k}) \s(h_{1})$.
So $w_{1}w_{2} \in  W_{g_{1}, \dots g_{k}, h_{1}, \dots h_{l}}$.
It follows that $V_{r}^{2}\subseteq V_{r}$.
\par
Also, from the submultiplicativity of $\s$, we
see that $\alpha_{g}(W_{g_{1}, \dots g_{n}}) \subseteq \s^{2d}(g)
W_{gg_{1}, \dots gg_{n}}$.  It follows that
$\alpha_{g}(W_{n}) \subseteq \s^{2d}(g) W_{n}$ and
$\alpha_{g}(V_{r})\subseteq
\s^{2d}(g) V_{r}$.
\par
Define new seminorms by $\pa a \pa^{\prime}_{r} = \inf\{ \, c>0 \, |
\, a \in cV_{r} \}$.  Then   $\pa ab \pa^{\prime}_{r}\leq
\pa a\pa^{\prime}_{r} \pa b \pa^{\prime}_{r}$
and $\pa \alpha_{g}(a) \pa^{\prime}_{r}
\leq \s^{2d}\pa a \pa^{\prime}_{r}$.
Hence condition (3.1.19) implies that
the action of $G$ is $m$-$\s$-tempered.
\qed\enddemo
\subheading{Question 3.1.25} Let $\s$ be a weight on $G$.
Are there examples of $m$-convex Fr\'echet algebras
$A$ with $\s$-tempered (not $m$-$\s$-tempered)
actions of $G$ for which the crossed product
$\sonas$ is not $m$-convex ?
\subheading{Example 3.1.26} In a preliminary version of this paper
\cite{Sc 1, chap 1, \S 3}, I had
a different condition than $m$-temperedness to prove the $m$-convexity
of the crossed product. Let $\s\geq 1$ be a
weight on $G$. We say that an action of $G$ on $A$ is {\it
strongly $m$-$\s$-tempered} if for every family   (or equivalently for
any one family) $\bigl\{ \normm\bigr\}$
of increasing
seminorms for $A$ we have that
for every $m \in \N$, there exists $C>0$, $r\geq m$, $d \in \N$ such that
$$\pa
\alpha_{g_{1}}(a_{1}) \dots  \alpha_{g_{n}}(a_{n})
\pa_{m}
\leq C^{n}\biggl( \max_{\tau_{1}+\dots \tau_{n}\leq d}
\s^{\tau_{1}}(g_{1})\dots \s^{\tau_{n}}(g_{n}) \biggr)
\pa a_{1} \pa_{r} \dots \pa a_{n} \pa_{r},
\tag 3.1.27 $$
for all $n$-tuples $a_{1}\dots a_{n}\in A$, $g_{1}\dots g_{n}\in G$ and
all $n \in \N$.
(Note  that $r$ and $d$ do not depend on $n$.)
\par
This condition implies $m$-temperedness since if $\tau_{i_{1}}, \dots
\tau_{i_{d}}$ are the only nonzero $\tau$'s on the right of the
inequality (3.1.27), we have
$$\aligned &\s^{\tau_{1}}(g_{1})\dots \s^{\tau_{n}}(g_{n})
\leq
 \s^{d}(g_{i_{1}}) \dots \s^{d}(g_{i_{d}})\\
&\qquad \leq
\biggl(
(\s(g_{i_{1}}^{-1}g_{i_{1}+1}) \dots \s(g_{n-1}^{-1}g_{n})
\s(g_{n})) \dots
(\s(g_{i_{d}}^{-1}g_{i_{d}+1}) \dots \s(g_{n-1}^{-1}g_{n})
\s(g_{n}))
\biggr)^{d}
\\
&\qquad \leq
\biggl(\s(g_{1}^{-1}g_{2}) \dots \s(g_{n-1}^{-1}g_{n})
\s(g_{n}) \biggr)^{d^{2}}\endaligned$$
by the submultiplicative property of the weight $\s$.
\par
We give an example of an $m$-tempered action which is not strongly
$m$-tempered.
Let $\Q^{\infty}$ be the direct sum
$$  \bigoplus_{n=0}^{\infty} \Q$$
of countably many copies of $\Q$ with pointwise addition.
Let $\w$ be the weight on $\Q^{\infty}$ defined by
$$\w({\vec r}) = \prod_{i=1}^{\infty} ( 1 + |r_{i}|).  $$
(Note that if $\vr \in \Q^{\infty}$, only finitely many $r_{i}$'s are
nonzero, so the definition of $\w$ makes sense.)
Let $A$ be the  commutative Banach algebra
$l^{1}(\Q^{\infty}, \w)$ of functions from $\Q^{\infty}$ to
$\C$ which satisfy
$$ \pa \varphi \pa_{1} = \sum_{\Q^{\infty}} \w({\vr})
|\varphi({\vr})| < \infty,\tag 3.1.28$$
with convolution multilplication.
The norm $\pa \quad \pa_{1}$ is easily seen to be submultiplicative.
\par
Let $G$ be the discrete group consisting of tuples
$(q_{1}, \dots ,q_{n}, 1, 1,\dots)$ of infinite length,
with only finitely many entries not equal to one, and each
$q_{i}$ a positive rational number.
With pointwise multiplication,  $G$
is easily seen to be a group.
(Note that $G$ is not finitely generated.)
Define a weight $\ga$ on $G$ by
$$\ga(\vq) = \prod_{i=1}^{\infty} \max(q_{i}, 1/q_{i}). $$
Define an action $\alpha$ of $G$ on $A$ by
$$\alpha_{\vq}(\varphi )(\vr ) = \varphi({\vq}\, . \,{\vr} ), $$
where $\vq \,.\, \vr =
(q_{1}, \dots q_{n}, 1,\dots)\,.\, (r_{1}, \dots r_{n},0,\dots)$
is given by pointwise multiplication
$ (q_{1}r_{1}, \dots q_{n}r_{n}, 0, \dots)$.
Then $\alpha_{1/{\vq}} = \alpha_{\vq}^{-1}$, and it is not too
hard to show that $\alpha$ gives an action of $G$ on $A$
by algebra automorphisms.  We have
$$ \pa
\alpha_{\vq}(\varphi) \pa_{1} = \sum_{\vr \in \Q^{\infty}}
\w(\vr) |\varphi({\vq}\,.\, {\vr} )| =
 \sum_{\vr \in \Q^{\infty}}
\w(1/{\vq}\,.\, {\vr}) |\varphi( \vr )| \leq \ga(\vq) \pa \varphi \pa_{1},
\tag 3.1.29$$
since $(1+|r_{i}/q_{i}|) \leq \max(q_{i}, 1/q_{i}) (1+|r_{i}|)$.
Hence $\alpha_{\vq}$ is a continuous automorphism of $A$.
In fact, by (3.1.29), the submultiplicative norm $\pa \quad \pa_{1}$ on $A$
is ($G$, $\ga$)-stable.  Hence the action of $G$ on $A$
is $m$-$\ga$-tempered,
and the smooth crossed product $G\rtimes^{\ga}A$ is $m$-convex
by Theorem 3.1.7.
\par
We show that this action is {\it not} strongly $m$-$\ga$-tempered.
Let $e_{i}$ be the element $(0, \dots0, 1,0,\dots) \in \Q^{\infty}$ where
the 1 is in the $i$th spot.  Let $\delta_{\vq} \in A$ be the delta
function at $\vq\in \Q^{\infty}$.  Let $q\in\Q$ with $q>1$ and let $\vq_{i}$
denote $(1, \dots1, q, 1, \dots)\in G$,
where the $q$ is in the $i$th spot.
Then
$$ \alpha_{\vq_{1}}(\delta_{e_{1}}) *\dots \alpha_{\vq_{n}}(\delta_{e_{n}})
=\delta_{(q,0,\dots)} *\dots \delta_{(0,\dots, q,0,\dots)}
=
\delta_{(q, \dots q,0,\dots)},
$$
where $n$ $q$'s occur.
The norm in $A$ of $\delta_{(q, \dots q,0,\dots)}$ is precisely
$\w(q,\dots q,0,\dots)= (1+q)^{n}$.
If the action were strongly $m$-tempered, we would
have this norm bounded by
$$ \aligned \max_{\tau_{1} +\dots \tau_{n}\leq d}& \ga^{\tau_{1}}(\vq_{1})
\dots \ga^{\tau_{n}}(\vq_{n})\pa \delta_{e_{1}} \pa_{1}
\dots \pa \delta_{e_{n}} \pa_{1} \\
& \leq \max_{\tau_{1} +\dots \tau_{n}\leq d} q^{\tau_{1}}
\dots q^{\tau_{n}}\w(e_{1})
\dots \w(e_{n})
 \leq q^{d} 2^{n}\endaligned $$
times some fixed constant $C$ to the $n$th power.  This would
have to hold for any $q$, which cannot be since $(1+q)^{n} $
is not bounded by $q^{d}2^{n}C^{n}$
as $q$ tends to infinity. Hence the action $\alpha$ of $G$ on
$A$ is $m$-$\ga$-tempered but not strongly $m$-$\ga$-tempered.
\vskip\baselineskip
\subheading{Question 3.1.30} If $G$ is a compactly generated
group, do the notions of $m$-tempered and strongly $m$-tempered
coincide ?
\vskip\baselineskip
\heading \S 3.2 Non $m$-convex Group Algebras \endheading
\par
We construct several non $m$-convex group algebras,  using
sub-polynomial scales which are not $m$-sub-polynomial.
\par
Let $G$ be a locally compact compactly generated group, and let $\tau$
be the word gauge on $G$.  Define scales $\s_{k}$ on $G$
by
$$ \s_{k}(g) = e^{\tau (g)^{k}} \tag 3.2.1 $$
\proclaim{Proposition 3.2.2}  If $G$ is compact, then
every $\s_{k}$ is equivalent to $1$.
Assume $G$ is not compact.  Then the  scales $\s_{k}$
form a family of increasing inequivalent sub-polynomial scales.
If $k$ is zero,
$\s_{k}\equiv 1$, and if $k=1$, $\s_{k}$ is the exponentiated
word gauge.  If $k\geq 2$, then $\s_{k}$ is not $m$-sub-polynomial.
\endproclaim
\demo{Proof}
Assume $G$ is not compact. Clearly $\s_{k+1}\geq \s_{k}$
since $\tau $ is integer valued.  We show that $\s_{k}$ cannot
dominate $\s_{k+1}$. Assume  that
$$\s_{k+1}(g)\leq C\s_{k}^{d}(g), \qquad g\in G. $$
Since $G$  is not compact, $\tau$ can take on arbitrarily large
integer values.  So we have  $e^{m^{k+1}} \leq C e^{d{m}^{k}}$
for arbitrarily large values of $m$.  This is impossible since
$C$ and $d$ are fixed.
\par
 We show that if $k\geq 2$, then $\s_{k}$ is not $m$-sub-polynomial.
Assume that
$$\s_{k}(g_{1}\dots g_{n}) \leq C^{n}\s_{k}^{d}(g_{1})
\dots \s_{k}^{d}(g_{n}). \tag 3.2.3$$
For $n \in \N$,
let $g_{1}, \dots g_{n}$ be elements of the generating set $U$,
such that $\tau(g_{1}\dots g_{n})= n$.  Then
$\s_{k}(g_{i})=e$ and
by (3.2.3) we have
$$ e^{n^{k}} \leq C^{n} e^{dn}. \tag 3.2.4 $$
But since $G$ is not compact, $n$ can be taken arbitrarily large with
$d$ and $C$ fixed.  So (3.2.4) will be violated if $k\geq 2$.  Hence
$\s_{k}$ is not $m$-sub-polynomial if $k\geq 2$.
\par
We show that each $\s_{k}$ is sub-polynomial.  To do this, we show that
$$  \tau^{k}(gh) \leq  2^{k} (\tau^{k}(g) + \tau^{k}(h) ).  \tag 3.2.5$$
By the subadditivity of $\tau$, it suffices to show that $(a+b)^{k}
\leq 2^{k}(a^{k} + b^{k} )$ for nonnegative integers $a$ and $b$.
For this, it suffices to show that the real valued function
$f\colon [0, \infty) \times [0, \infty )  \rightarrow \R^{+}$
given by
$$f(x, y)= {(x+y)^{k}\over {x^{k} + y^{k}}} $$
is bounded by $2^{k}$.
The function $f$ is globally bounded since it is continuous
and $\lim_{z\rightarrow \infty} (1+z)^{k}/(1+z^{k})$ is one.
Setting the first partial
derivatives of $f$
to zero, we see that $f$ has a local maximum only if $x= \pm y$.
If $x=y$, then $f(x, y)=2^{k}$.
So it follows that $f\leq 2^{k}$, and in fact that $2^{k}$ is the
best bound.
\qed
\enddemo
\subheading{Definition 3.2.6} If $G$ is a compactly generated Lie group
and $k\geq 1$, we let $\Cal S^{k}(G)$
be the Schwartz algebra corresponding
to the scale $\s_{k}$.
We let $L^{k}(G)$ denote the Fr\'echet *-algebra $L_{1}^{\s_{k}}(G)$.
\par
Note that $\s_{k}$ dominates the exponentiated
word weight $\s_{1}$, and
hence dominates every weight on $G$ - so in particular $\s_{k}$
bounds $Ad$  and $\Cal S^{k} (G)$ is a Fr\'echet *-algebra.
\par
Let $G$ be any finitely generated discrete group which is not finite.
For example, $G$ could be the discrete subgroup
$$G= \pmatrix 1 & \Z & \Z \\ 0 & 1 & \Z \\ 0 & 0 & 1 \endpmatrix
\tag 3.2.7$$
of the three dimensional Heisenberg Lie group.
Let $k\geq 2$.
Then $\Cal S^{k}(G)$ in not $m$-convex by the following
proposition and Proposition 3.2.2.
Note that for the group (3.2.7) this gives an example
of a noncommutative non $m$-convex Fr\'echet *-algebra.
\proclaim {Theorem 3.2.8} Assume that $G$ is a discrete
group.  Let $\s$ be a sub-polynomial
scale on $G$.  Then the convolution algebra
$\schwsi$ is $m$-convex if and only if $\s$ is
$m$-sub-polynomial. \endproclaim
\demo{Proof}
By Theorem 3.1.7, it suffices to show that if $\schwsi$ is
$m$-convex, then $\s$ is $m$-sub-polynomial.
Let $e_{g}$ be the function in $\schwsi$
defined by
$$e_{g}(h) =\cases 1 & {\text {if $h=g$}} \\
0 &  {\text {if $h\not= g$}} \endcases
\tag 3.2.9  $$
For convenience, we replace $\s$ with the equivalent scale $1+ \s$ so
that $\s \geq 1$.  A family of increasing seminorms for $\schwsi$
is given by
$$ \pa \varphi \pa_{m} = \sum_{g \in G} \,\s^{m}(g) \, | \, \varphi(g) \, |
.$$
Then
$$
\aligned \biggl(
{ \pa e_{g_{1}}\dots
e_{g_{n}}\pa_{m}\over
{\pa e_{g_{1}} \pa_{k}
\dots \pa e_{g_{n}} \pa_{k} }}
\biggr)^{1/n}
& =
\biggl(
{ \pa e_{{g_{1}}\dots
{g_{n}}}\pa_{m}\over
{\pa e_{g_{1}} \pa_{k}
\dots \pa e_{g_{n}} \pa_{k} }}
\biggr)^{1/n}
\\
&=
\biggl(
{\s^{m}(g_{1}\dots g_{n})\over
{\s^{k}(g_{1})
\dots \s^{k}(g_{n}) }}
\biggr)^{1/n}.
\endaligned \tag 3.2.10
$$
By Theorem 3.1.4, the $m$-convexity of $\schwsi$ implies that for all $m$
there is some $k$ such that (3.2.10) is bounded by a constant $C$.  In
particular, taking $m=1$ we see that
$$
\s(g_{1}\dots g_{n})
\leq C^{n}\s^{k}(g_{1}) \dots \s^{k}(g_{n})
\tag 3.2.11
$$
and $\s$ is $m$-sub-polynomial.
\qed
\enddemo
\subheading{Question 3.2.12} Is it true for a general Lie group
$G$ that $\schwsi$ is
$m$-convex if and only if $\s$ is $m$-sub-polynomial ?
\vskip\baselineskip
\par
We give some examples of non $m$-convex
Schwartz algebras for connected Lie groups.
\proclaim{Proposition 3.2.13} Let $G$ be the Abelian group  $\R^{N}$.
Let $\s_{k}$ be
the non $m$-sub-polynomial scale defined above for any  $k\geq 2$.
Then $\s_{k}$ is equivalent to the sub-polynomial scale
$$\delta_{k}(x)= e^{|x|^{k}}$$
on $\R^{N}$.
The
convolution algebras $\Cal S^{k} (G)$ and
$L^{{k}}(G)$ are {\it not} $m$-convex.
\endproclaim
\demo{Proof} To see that $\delta_{k}\thicksim \s_{k}$,
take as generating set for $G$ the open unit ball $U$.
It suffices to show that
$|x|\leq \tau(x) \leq |x| + 1$,
where $\tau$ is the word gauge corresponding to $U$.
If $\tau (x) = n$, then $x= x_{1}+\dots x_{n}$ with
$x_{i}\in U$.  Hence $|x| \leq n$ and we have
$|x| \leq \tau(x)$ for all $x$.
If $|x|=c$, let $n$ be the smallest integer greater than
or equal to $c$.  Then $|x/n|<1$ and $\tau(x) \leq n \leq c + 1$.
\par
To show that the Schwartz algebra is not $m$-convex,
we use Theorem 3.1.4.
To simplify notation, we denote $\delta_{k}$ by
$\delta$, where $k$ is assumed greater than or equal to $2$.
We topologize $\Cal S^{k} (G)$
by the increasing family of seminorms
$$ \pa \psi \pa_{m}^{\prime} =
\max_{|\ga|\leq m} \int_{G} \delta^{m}
(x) |X^{\ga}\psi(x)|dx, \tag 3.2.14 $$
and topologize $L^{k}(G)$ by the increasing family
$\pa \quad \pa_{m}$ of seminorms
$$ \pa \psi \pa_{m} = \int_{G} \delta^{m}(x) \,|\psi(x)|\, dx.
\tag 3.2.15 $$
Clearly $\pa \psi \pa_{m}\leq
\pa \psi \pa_{m}^{\prime}$ for $\psi \in \Cal S^{k} (G)$.
\par
Let $\psi$ be any positive $C^{\infty}$-function satisfying
$$ \psi(x) = \cases 0 & |x|>2 \\
1 & |x|\leq 1 \endcases \tag 3.2.16 $$
Let $ u$ be a unit vector.
Define $\psi_{n}(x) = \psi(x-nu)$.
If $k\geq 2$, we show that for any $l\in \N$ and
for any $C>0$ , there is some $n$ sufficiently large
such that
$$ \biggl(
\undersetbrace \text {$n$-times} \to { \pa \psi_{1}\dots
\psi_{1}\pa_{1}\over
{\pa \psi_{1} \pa_{l}^{\prime}
\dots \pa \psi_{1} \pa_{l}^{\prime} }}
\biggr)^{1/n}
>C \tag 3.2.17 $$
By Theorem 3.1.4, this will establish that $\Cal S^{k} (G)$
is not $m$-convex.  And since the prime norms dominate the
unprimed norms, it will also establish that
$L^{k}(G)$ is not $m$-convex.
\par
Since the denominator of (3.2.17) is just
$\pa \psi_{1} \pa_{l}^{\prime}$,
it suffices to show that $\pa \psi_{1}^{*n}\pa_{1}^{1/n}$
tends to infinity as $n\longrightarrow \infty$.
We estimate $\pa \psi_{1}^{*n}\pa_{1}$.
$$ \aligned \pa \psi_{1}^{*n}
\pa_{1} & =\int_{G} \delta(x)
|\psi_{1}^{*n}(x)|dx \\
& = \int_{G} \delta(x+nu) |\psi^{*n}(x)|dx. \endaligned
\tag 3.2.18
$$
We write
$$
|\psi^{*n}(x)|=
\int_{G} \dots \int_{G}
\psi(x_{1}) \dots \psi(x_{n-1}) \psi(x-x_{1}-\dots
x_{n-1})dx_{1}\dots dx_{n-1} \tag 3.2.19 $$
since $\psi$ is positive.
For definiteness, we choose $|x|$ to
be the maximum of the components of the
vector $x$.
Let $W$ be a  neighborhood of $0$ such that $W^{n} $
is contained in the unit ball with respect to
$|\quad |$, and the Lebesgue measure  of $W$ is
greater than or equal to $1/n^{N}$ (For example $W$ could
be the unit ball of $\R^{N}$
to the $1/n$th power.)
If, in the integral (3.2.19),
we restrict  $x_{1}, \dots x_{n-1}$ to be within the
set $W$, then $\psi(x_{i})=1$
always.  Thus
$$ | \psi^{*n}(x)|\geq
\int_{W^{n-1}}\psi (x-x_{1}-\dots x_{n-1})dx_{1}\dots dx_{n-1}.
\tag 3.2.20 $$
If we  make the further restriction that $x \in W$,
then  the integrand of (3.2.20) is greater than or equal to one and
$$ | {\psi^{*n}(x)}|   \geq  \int_{W^{n-1}}dx_{1}\dots dx_{n-1}
\geq 1/n^{(n-1)N}. \tag 3.2.21 $$
Plugging into (3.2.18), we see that
$$ \aligned \pa \psi_{1}^{*n}\pa_{1} & \geq
\int_{W} \delta(x+nu)  | {\psi*\dots \psi(x)}|dx
\\
& \geq \inf_{x\in W} \delta(x+nu)
\int_{W}1/n^{(n-1)N} dx
\geq \delta((n-1)u) /n^{nN} =
e^{(n-1)^{k}} /n^{nN}. \endaligned \tag 3.2.22 $$
Thus
$$
\pa \psi_{1}^{*n}\pa_{1}^{1/n}  \geq
{e^{(n-1)^{k-1}}\over {(n)^{N}}}.
\tag 3.2.23 $$
Since $k\geq 2$,  $e^{(n-1)^{k-1}}$ grows much
faster than $(n)^{N}$.  So given any $C>0$,
it is always possible to find an $n$ sufficiently large so that (3.2.17)
holds.  This proves that
$\Cal S^{k} (G)$ and $L^{k}(G)$
are not $m$-convex by Theorem 3.1.4,
and completes the proof of Proposition
3.2.13
\qed \enddemo
\par
For other constructions of non $m$-convex Fr\'echet algebras,
see \cite{Ar}, \cite{RZ, \S 7}, \cite{Ze, Theorem 26},
and \cite{Pe}.
\heading  \S 4 Conditions for a *-algebra  \endheading
\par
Let $\s$ be a weight or a gauge on $G$ that bounds $Ad$.
The main result of this section is that if $A_{1}$ is a Fr\'echet
*-algebra,
with a tempered and strongly continuous action of
 $G$ by *-automorphisms,
 and $A$ is
the set of $C^{\infty}$-vectors for the
action of  $G$ on $A_{1}$,
then $\sonas$ is a Fr\'echet *-algebra
(see Corollary 4.9 below).
We first show that a tempered  action of $G$
on $A_{1}$ gives a tempered  action of $G$
on $A$.  This fact will prove useful in examples, since the
temperedness of an action on $A_{1}$ is often easier to verify
than the temperedness on $A$.  We also prove the appropriate $m$-versions
of our results.
\par
See the beginning of \S 1.3 for the definition
of a Fr\'echet *-algebra.
We remark that
for a Fr\'echet *-algebra $A$, there always exists a set of
*-isometric seminorms for $A$.  Simply define
$\pa a \pa_{m}^{\prime} = \max{ ( \pa  a \pa_{m}, \pa a^{*} \pa_{m} )}$.
({\it *-isometric} means $\pa a^{*} \pa^{\prime}_{m}
= \pa a \pa'_{m}$ for all $a \in A$ and $m \in \N$.)
Similarly, an $m$-convex Fr\'echet *-algebra is topologized by a
family of submultiplicative *-isometric seminorms.
We say that $G$ acts by {\it *-automorphisms} on $A$
if $\alpha_{g}(a^{*})=
\alpha_{g}(a)^{*}$ for $a \in A$ and  $g\in G$.
\par
Let $A$ be a Fr\'echet algebra with tempered and strongly continuous
action of $G$.
Even if $A$ is a Fr\'echet *-algebra and $G$ acts by *-automorphisms,
the convolution
algebra $\sona$ may not be closed under the involution $\varphi^{*}(g)=
\Delta (g) \alpha_{g}(\varphi(g^{-1})^{*})$.
For example, let $G=\R$ with $\s$ given by the absolute
value, and let $A_{}=C_{0}(\R)$ where $\R$ acts by
translation.  Then $\sona= \Cal S (\R, C_{0}(\R))$,
the Schwartz functions on $\R$ taking values in $C_{0}(\R)$.
The function
$$ \psi(r, s) = e^{-r^{2}} {1\over {1+ |s|}} \tag 4.1 $$
is in $\Cal S (\R, C_{0}(\R))$, but
$$\psi^{*}(r, s) = e^{-r^{2}}{1\over {1 +|s-r|}} $$
is not differentiable in the variable $r$.  In this
case, note that $\sona$ would be closed under the involution
if the action
of $G$ on $A$ were differentiable.
\par
Let $H$ be a Lie group, and assume $G\subseteq H$ with differentiable
inclusion map.
Let $A$ be the set of
$C^{\infty}$-vectors for a strongly continuous action
of $H$ on a  Fr\'echet algebra $A_{1}$.
(We shall always assume
that the action of $G$ on $A_{1}$ is then given by the
restriction of this action of $H$.)
Let $\basisp$ be a
basis for the Lie algebra $\frak H$ of $H$, where $p$ is the
dimension of $\frak H$.
If $a\in A$ and $\ga \in \N^{p}$, let $X^{\ga}a$
be given by formula (1.2.1) with $\be = \alpha$ and $q=p$,
where $\alpha$ denotes
the action of $H$ on $A$.
The algebra $A$ is topologized by the seminorms
$$\pa a \pa_{l,m} = \max_{|\ga|\leq l}
\pa X^{\ga}a \pa_{m}. \tag 4.2 $$
In this topology, $A$ is a Fr\'echet algebra
($m$-convex if $A_{1}$ is), and is a dense
subalgebra of $A_{1}$ with continuous inclusion.
Also, the Lie group $H$ leaves $A$ invariant and acts
pointwise differentiably on $A$ by continuous automorphisms.
See Theorem A.2 in the appendix  for proofs of these  facts.
\subheading{Definition 4.3} We say that a scale $\s$ on $G$
{\it bounds $Ad$ on $H$} if there exists $C, D \geq 0$
and $d\in \N$ such that
$$ \pa Ad_{g} \pa_{\frak H} \leq C \s^{d}(g) + D, \tag 4.4 $$
for each $g\in G$.  Here the norm is taken as an operator on the
Lie algebra $\frak H$ of $H$.
\subheading{Example 4.5} We give an example of
a scale $\s$ on $G$ which bounds
$Ad$ on $G$ but not on $H$.  Let
$G$ be the subgroup of integer valued matrices of the
three dimensional Heisenberg Lie group
$$ H = \pmatrix 1 & \R & \R \\ 0 & 1 & \R \\ 0 & 0 & 1 \endpmatrix. $$
If  $\s$ is the trivial scale, then
$\s$ of course bounds $Ad$ on $G$ since $G$ is discrete.
However, the norm of $Ad_{g}$ as an operator on $\frak H$ is
equivalent to the sum of the  two matrix entries just above the
diagonal (see Example 1.5.15), so $\pa Ad_{g} \pa_{\frak H}$
is unbounded for $g \in G$ and $\s$ does not bound $Ad$ on $H$.
\par
We prove the following useful theorem.
Recall that $\s_{-}(g) = \s(g^{-1})$.
\proclaim {Theorem 4.6} Let $A_{1}$ be an ($m$-convex)
Fr\'echet algebra with strongly continuous action
of a Lie group $H$.  Let $G$ be a Lie group and assume $G\subseteq H$
with differentiable inclusion map.  Let $\s$ be a
scale on $G$ such that $\s_{-}$ bounds $Ad$ on $H$.
Assume that $A_{1}$ is an ($m$-convex)
Fr\'echet algebra on which the action of
$G$ is $\s$-tempered ($m$-$\s$-tempered).
Then the algebra $A$ of $C^{\infty}$-vectors
for the action of $H$ on $A_{1}$ is an ($m$-convex) Fr\'echet
algebra, and the action of $G$ on $A$ is $\s$-tempered
($m$-$\s$-tempered).
\par
If $A_{1}$ is a Fr\'echet *-algebra and $H$ acts by *-automorphisms
on $A_{1}$, then $A$ is a Fr\'echet *-algebra for which the above
statements hold.
\endproclaim
\demo{Proof}
For the first and last statements, see Theorem A.2.
We prove the second statement.
To simplify notation, assume that $\s \geq 1$, by replacing
$\s$ with $1+ \s$ or $\max (\s, 1)$.
Since  $\s_{-}$ bounds $Ad$ on $H$, by Lemma 2.2.1 we have
$$
\pa X^{\ga}\alpha_{g}(a)\pa_{m}\leq
C\s^{d}(g)  \sum_{\be\leq \ga} \pa
\alpha_{g}(X^{\be}(a))\pa_{m},
\quad  a\in { A}, \quad  g\in G.
\tag 4.7
$$
Therefore
$$
\pa \alpha_{g}(a) \pa_{l, m} \leq K\s^{d}(g) \max_{|\be|\leq l}
\pa \alpha_{g} (X^{\be}(a)) \pa_{m}. \tag 4.8 $$
If the action of $G$ on $A_{1}$ is tempered, (4.8) tells us that
${\pa \alpha_{g}(a) \pa_{l, m}} \leq
{\tilde K} \s^{d+k}(g) {\pa a \pa_{l,t}}$
for some $t\geq m$, so the action of $G$ on $A$ is tempered.
If the action of $G$ on $A_{1}$ is $m$-tempered,
and $\pa \quad \pa_{m}$ are $(G, \s)$-stable submultiplicative
seminorms for $A_{1}$, then we have
$\pa \alpha_{g}(a) \pa_{l, m} \leq  {\tilde K} \s^{d+k}(g) \pa a \pa_{l,m}$,
so the seminorms $\pa \quad \pa_{l, m}$ on $A$ are $(G, \s)$-stable.
Also, a constant times $\pa \quad \pa_{l, m}$ gives
$(G, \s)$-stable submultiplicative
seminorms for $A$ (see (A.6) and
following remarks), so the action of $G$ on $A$
is $m$-tempered.
\qed
\enddemo
\proclaim{Corollary 4.9} Let $G$, $H$, and $A_{1}$  be as in Theorem 4.6
above, and
assume that $A_{1}$ is an ($m$-convex)
Fr\'echet *-algebra on which $H$ acts by *-automorphisms.  Let
$A$ be the ($m$-convex) Fr\'echet *-algebra
of $C^{\infty}$-vectors for the action of $H$
on $A_{1}$.
If the action of $G$ on $A_{1}$ is
$\s$-tempered ($m$-$\s$-tempered), and $\s$ is a sub-polynomial
($m$-sub-polynomial) scale
on $G$ which bounds $Ad$ on $H$, and satisfies $\s\thicksim\s_{-}$
(for example,  $\s$ could be a weight or a gauge on $G$ which bounds
$Ad$ on $H$),
then the convolution algebra
$\sonas$ is an ($m$-convex) Fr\'echet *-algebra.
\endproclaim
\demo{Proof}
For everything except the *-algebra part, see Theorems 2.2.6, 3.1.7,
and 4.6 above.  It remains to show that involution is well
defined and continuous on the algebra $\sona$.
We let $\pa \quad \pa_{d} $
be seminorms (4.2)
for $A$, and let $\pa \quad \pa_{m, \ga, d}$
be seminorms (2.1.1) for $\sona$ as in Definition 2.1.0.
Just as in the proof of Theorem 4.6 above,
we may assume that $\s \geq 1$.
Recall that $\psi^{*}(g) = \alpha_{g}(\psi(g^{-1})^{*})\Delta (g)$,
where $\Delta$ is the modular function for $G$.
We must show that
$$ \pa \psi^{*} \pa_{m, \ga, d} =
\int_{G} \pa \s^{m} X^{\ga}\psi^{*} (g)
\pa_{d} dg \tag 4.10$$
is bounded by a linear combination of seminorms in $\psi$. We
do most of the work in the following lemma.
\proclaim{Lemma 4.11}
 Let $\ga \in \N^{q}$ (where $q$ is
the dimension of the Lie algebra $\frak G$
of $G$) and assume
that $\s\geq 1$ bounds $Ad$ on $G$, and $\s\thicksim \s_{-}$.
Then there is some $j\in \N$
so that
$$\pa (X^{\ga}\psi^{*})(g) \pa_{m}\leq
\sum_{\be, {\tilde \be}\leq \ga} C\s^{j}(g)
\pa \alpha_{g}(X^{\be}((X^{\tilde \be}\psi)(g^{-1})))^{*}
\pa_{m}, \tag 4.12 $$
where the operator $X^{\tilde \be}$ is from
the Lie algebra $\frak G$ of $G$
acting by left translation on $\psi\in \sona$, and $X^{\be}$ is from
$\frak G$ acting via $\alpha$ and the  inclusion $\frak G \subseteq
\frak H$.
\endproclaim
\demo{Proof} It is well known fact that the modular
function  $\Delta
\colon G\longrightarrow \R^{+}$ is differentiable
\cite{War, Vol. I}.
We show that all of the derivatives of $\Delta$
 are bounded
by some $\s^{p}$.
Let $\tilde K$ be a constant such that
 $|(X^{\be}\Delta)(e)|\leq {\tilde K} $
for all $\be \leq \ga$.
Since $\Delta $ is the absolute value of the determinant
of $Ad$ \cite{War, Vol. I, Appendix}
and $\s$ bounds $Ad$, there are $C$ and $p$ such that
$$\Delta (g) \leq C\s^{p}(g). \tag 4.13$$
So by the multiplicativity of $\Delta$, we have
$$\aligned |(X^{\be}\Delta )(g)| & \leq {\tilde K}\Delta (g) \\
& \leq K\s^{p}(g) \quad {\text {by (4.13)}}
\endaligned \tag 4.14 $$
where $K= {\tilde K}C$ and $\be\leq \ga$.  This
inequality will be very useful in our calculations.
\par
Using the product rule and the chain rule, and the formula
$\psi^{*}(g) = \alpha_{g}(\psi (g^{-1})^{*}) \Delta(g)$,
we see that
for $\ga\in \N^{q}$,
$$\split & (X^{\ga} \psi^{*})(g)  \\
& = \sum_{\be_{1},
{\be_{2}}, \be_{3}\leq \ga}
\alpha_{g} (X^{\be_{1}}((X^{ {\be_{2}}}\psi)(g^{-1})^{*}))
(X^{\be_{3}}\Delta)(g))
P_{
\be_{1},  {\be_{2}}, \be_{3}}((Ad_{g^{-1}})_{ji})
\endsplit
\tag 4.15
$$
where the $P_{\be_{1},  {\be_{2}}, \be_{3}}$
are polynomials with degree bounded by $|\ga|$.
Note that the matrix entries $(Ad_{g^{-1}})_{ij}$ in (4.15)
are from $Ad_{g^{-1}}$ as an operator
 on the Lie algebra of $G$, not just $H$ (since $X^{\ga}$ and
all the other differential operators in (4.15) come from $\frak G$).
Since $\s_{-}\thicksim \s $, and  $\s$ bounds $Ad$ on $G$,
there is some $s\in \N$ and
$C_{1}$ such that
$$ |P_{\be_{1},  {\be_{2}}, \be_{3}}
((Ad_{g^{-1}})_{ji})|
\leq C_{1}\s^{s}(g). \tag 4.16  $$
Estimating $\pa X^{\ga}\psi^{*} (g) \pa_{m}$
using (4.15), (4.16), and (4.14),
we get
$$\aligned
\pa (X^{\ga}\psi^{*})(g) \pa_{m} & \leq
\sum_{\be_{1},
 {\be_{2}}\leq \ga}
C_{1}K \s^{s+p}(g)
\pa \alpha_{g}(X^{\be_{1}}((X^{ \be_{2}}\psi)(g^{-1})))^{*}
\pa_{m}
\endaligned \tag 4.17 $$
Taking $j=p+s $ and $C= C_{1}K$ yields (4.12) and
proves Lemma 4.11.
\qed \enddemo
We have
$$\aligned
& \pa \Delta (X^{\ga}\psi^{*})(g^{-1})\pa_{m} \\
&\leq
\sum_{\be, {\tilde \be}\leq \ga} C_{1} \s^{j}(g)
\pa  \alpha_{g^{-1}}(X^{\be}((X^{\tilde \be}\psi)(g)))^{*}
\pa_{m}
 \qquad  {\text{  by Lemma 4.11, (4.14)}}
\\
&
\leq \sum_{\be, {\tilde \be}\leq \ga}{ C_{2}\s^{r}(g)
\pa  X^{\be}((X^{\tilde \be}\psi)(g))^{*}
\pa_{m}}
\quad {\text{tempered action, $\s_{-} \thicksim \s$}}\\
&\leq
\sum_{{\tilde \be}\leq \ga} {C_{3}}
\pa \s^{r} (X^{\tilde \be}\psi)(g)
\pa_{p} \quad {\text{continuity of * and of differentiation on $A$,}}
\endaligned \tag 4.18
$$
for constants $C_{1}, C_{2}, C_{3}>0$, and sufficiently large
integers $r, p$.
So we have
$$\aligned
 \pa \psi^{*} \pa_{m, \ga, d}
& = \int_{G} \pa \s^{m}(X^{\ga}\psi^{*} ) (g) \pa_{d}
dg \\
& =
\int_{G} \pa \s^{m} \Delta
 (X^{\ga}\psi^{*} ) (g^{-1}) \pa_{d}
dg \\
& \leq  \int_{G}
\biggl( \sum_{{\tilde \be}\leq \ga} {C_{3}}
\pa  \s^{m+r} (X^{\tilde \be}\psi)(g)
\pa_{p} \biggr) dg {\text {      (4.18)}}\\
& = \sum_{{\tilde \be}\leq \ga}
C_{3} \pa \psi \pa_{m+r,{\tilde  \be} , p}.
\endaligned \tag 4.19 $$
Thus $*$ is continuous.
This proves Corollary 4.9.  \qed \enddemo
\heading \S 5 Examples and Scaled $G$-Spaces  \endheading
\par
We develop the notion of a scaled $(G, \w)$-space, where $\w$ is a
sub-polynomial scale on $G$, which will give
us many examples of dense subalgebras (with tempered and
$m$-tempered action of $G$) of the commutative
C*-algebra $B= C_{0}(M)$.  (Here $C_{0}(M)$ denotes the continuous
functions on $M$, vanishing at infinity with sup norm.)
We then look at
examples of scaled $(G,\w)$-spaces, and show that we can form
the smooth crossed product $\GSno$.
We also make note of some examples of
dense subalgebras of $B$ in the case that
$B$ is a noncommutative C*-algebra.
\par
To simplify our formulas, unless specified otherwise,
we assume that $\s \geq 1$. This
can always be achieved by replacing $\s$ with one of the
equivalent scales $\max(\s, 1)$ or
$1+\s$.
\subheading{Definition 5.1}
Let $M$ be a locally compact space.  Let $
\s\geq 1 $
be a scale on $M$ which is bounded on compact subsets of $M$.
\par
Let $\schwcm$ be
the set of functions
$$\{\, f\in C_0(M)\,| \, \pa \s^{d} f\pa_{\infty} < \infty
\quad \forall d\in \N \,\},$$
which we call the continuous functions which
vanish $\s$-rapidly.  It follows that $\s^{d} f$
vanishes at infinity for every $d\in \N$, even if $\s$
is not a proper map.
We topologize $\schwcm$ by the seminorms
$$\normdp = \pa \s^{d}{f }\pa_{\infty}, \qquad d=0,1,\dots
\tag 5.2$$
We will usually denote the seminorm $\normo$ by $\normi$.
If $\s$ and $\tau$ are equivalent scales, then $\schwcm =
C^{\tau}(M)$.
If $\{ {f }_n \}$ is a Cauchy sequence in $\schwcm$, then ${f }_n
\longrightarrow {f }_{0} $ in $\normi$ for some ${f }_{0} \in {C_0}(M)$.
Note
$$\s^d (m)|{f_{0}} (m)| \leq
\s^{d}(m)|{f_{0}} (m) - {{f }_k} (m)|  +
   {\pa {{f }_k} - {{f }_n} \pa_d}   +
{\pa {{f }_n}
\pa_d}
\tag * $$
for $m\in M$, $k,n \in \Bbb N$.  Let $N$ be so large that
$$n,k\geq N \Longrightarrow
{{\pa} {{f }_k} - {{f }_n} {\pa_d}} <1.$$
Take $n=N$
in $(*)$ and fix $m\in M$.  By letting $k$ run in $(*)$, we see
that $\s^{d}(m) |{f_{0}} (m)| \leq 2 + \pa {f }_{N} \pa_d$.
Hence ${f_{0}} \in \schwcm$.  Similar arguments show that
$f_{n} \longrightarrow f_{0}$ in $\schwcm$, so $\schwcm$ is complete.
The space $\schwcm$ is an $m$-convex (see \S 3) Fr\'echet *-algebra,
since
$$ \aligned \pa f_{1} \dots f_{n} \pa_{d} & =
\pa\s^{d}f_{1} \dots f_{n}\pa_{\infty} \\
& \leq \pa f_{1} \pa_{d}
\pa f_{2} \pa_{\infty}
\dots \pa f_{n} \pa_{\infty}. \endaligned
\tag 5.3 $$
Since   $\s$ is bounded on compact sets,
the compactly supported continuous
functions
$C_{c}(M)$ are contained
 in $\schwcm$, so $\schwcm$ is dense in $C_0(M)$.
The seminorms (5.2) are continuous for the inductive limit
topology on $C_{c}(M)$.
\subheading {Remark 5.4} If  $\s$ is not
bounded on compact subsets,
then $\schwcm$ may not be dense in $C_0(M)$.  For let $M=\R$ and set
$$\s (r) = \cases {1\over |r|} + |r| & r\neq 0.\\
0 & r=0. \endcases\tag 5.5 $$
Then every element of $ \schwcm $ vanishes at $0$.
\subheading{Definition 5.6}
Now assume that $M$ is an $H$-set, with $H$ a
Lie group.
If $f\in \schwcm$, define $\alpha_{h}(f)(m) = f(h^{-1}m)$.
We impose a condition that makes $\alpha$ a strongly continuous
action of $H$ on
$\schwcm$.
We say that $\s$ is {\it  uniformly $H$-translationally equivalent} if
for every compact subset $K\subseteq H$, there exists  $ l \in \N$ and
$ C>0$ so that
$$
\s(hm)\leq C\s^{l}(m), \qquad h \in K, m \in M.
\tag 5.7 $$
Just as we saw in Theorem 1.2.11
that for scales on a group, translational equivalence
implies uniform translational equivalence,
it is probably also true that the notion of $H$-translational
equivalence (which we haven't stated)
implies uniform $H$-translational equivalence.
We shall not be needing this, however.
We show that if $\s$ is uniformly $H$-translationally equivalent,
then $H$ acts on
$\schwcm$ by continuous
automorphisms, and moreover that the action is
strongly continuous.
For $h\in K$, $f\in \schwcm$, we have
$$ \aligned \pa \alpha_{h}(f) \pa_{d} &=
\pa\s^{d} \alpha_{h}(f)\pa_{\infty}
= \sup_{m\in M} \s^{d}(m) |f(h^{-1}m)| \\
& = \sup_{m\in M} \s^{d}(hm) |f(m)|
 \leq C^d \pa \s^{ld}f \pa_{\infty}
 = C^d \pa f \pa_{ld}. \endaligned
\tag 5.8 $$
So $H$ clearly leaves
$\schwcm$ invariant and acts by bounded,
and hence continuous, automorphisms on $\schwcm$.
Since the action  of $H$ on $C_{c}(M)$ is strongly
continuous for the inductive limit topology, a
standard argument shows that the action of $H$
on $\schwcm$ is strongly continuous.
\par
If $\s$ is not uniformly
$H$-translationally equivalent, then $H$ may not leave
$\schwcm$ invariant as the following
example shows.
\subheading {Example 5.9}
Let $M= \R^{+}$, $\s (r) = 1+ |r|$, and $H= \Z / 2\Z$.
Let $H$ act on $M$ via $\alpha (r) = 1/r$.  Then
$\s$ is not uniformly $H$-translationally
equivalent, since (5.7) is not
satisfied.  The function
$f(r)= \min(r, e^{1-r^{2}})$ is in $\schwcm$, but
$\alpha(f)$ is not.
This is because $e^{1-r^{2}}$ vanishes faster at infinity than any
power of $1/r$, but $1/r$ does not.
\par
If $ \s $ is uniformly $H$-translationally equivalent, we
define the {\it $\s$-rapidly vanishing
$H$-Schwartz functions on $(M, \s)$},
denoted by $\schwmh$, to be the set of $C^{\infty}$-vectors
for the action of $H$ on
$\schwcm$.  Since what group we
are using is usually clear, we often write
$\schwm$ to abbreviate $\schwmh$.
If it is clear both what group and scale we are using,
we will often simply write $\Cal S(M)$.
If $M$ happens to be the group $H$, we
may write $\soinfh$, with
subscript $\infty$ to denote the sup norm,
to contrast the space
with the $L^{1}$ Schwartz functions $\soh$ on $H$.
\par
The space $\schwm$ has a natural Fr\'echet topology,
and is dense in $\schwcm$ and
hence dense in $C_{0}(M)$.  Also, the action of $H$ by
left translation leaves $\schwm$ invariant and the  action
of $H$ on $\schwm$ by continuous
automorphisms is pointwise differentiable.
The Fr\'echet space $\schwm$ is also
an $m$-convex Fr\'echet *-algebra.
See Theorem A.2  below for proofs of these  facts.
\par
Let $p$ be the dimension of the Lie algebra of $H$.
The topology on $\schwm$ is given by
the norms
$$\pa f\pa_{l,\ga} =
\pa\s^{l}X^{\ga}f\pa_{\infty}, \tag 5.10$$
where $\ga\in \N^{p}$, and the notation $X^{\ga}f$
is defined in equation (1.2.1) with
action $\be = \alpha$ and $q=p$.
\subheading{Example 5.11} Let $M=H=\R$, where $H$ acts
by translation,  and let $\s(r) = 1+|r|$.
Then $\s$ is uniformly $H$-translationally equivalent and $\schwm$
is just the standard set of Schwartz functions on $\R$.
\par
We would like to define a smooth crossed product $G \rtimes \schwmh$
where $G$ is an arbitrary subgroup of $H$.  For example,
we would like to be able to define smooth crossed products
$\R\rtimes \Cal S(\R)$ and
$\Z \rtimes \Cal S (\R)$.  To facilitate this,
we make the following definition.
\subheading{Definition 5.12}  Let $G$ be a Lie group, which is a
subgroup of $H$ with differentiable inclusion map.
Let $ \s$ be a
uniformly $H$-translationally equivalent scale on $M$.
Then we know that the action $\alpha$ of $H$ on $\schwmh$ by left
translation is differentiable.  This action clearly restricts to
a differentiable action of $G$.
Let $\w \geq 1$ be a sub-polynomial scale on $G$.
We impose a
condition that makes $\alpha$ an $m$-$\w$-tempered
action of $G$ on
$\schwmh$.
We say that {\it  $(M, \s, H)$ is a scaled $(G, \w)$-space} if
$\w_{-}$ bounds $Ad $ on $H$ (see Definition 4.3)
(for example, $\w$ could be a weight or a gauge on $G$
that bounds $Ad$ on $H$) and
there exists  $ l \in \N$ and
$ C>0$ so that
$$
\s(gm)\leq C\w^{l}(g)\s^{l}(m), \qquad g \in G, m \in M.
\tag 5.13 $$
If we say that $(M, \s)$ or simply $M$
is a scaled $(G, \w)$-space, with no
mention of a group $H$, then it will be implied that we are taking
$H=G$.  It is useful to note that
condition (5.13) implies that $\s$
is uniformly $G$-translationally equivalent, since $\w$ is bounded
on compact sets.
\subheading{Example 5.14} For example, let $M=H=\R$, $\s(r)=1+|r|$,
and $\Cal S_{H}^{\s}(\R) = \schwr$ as in Example 5.11. Let $G$ be the
closed subgroup $\Z$ of $H$, with $\w(n) = e^{|n|}$ or $\w(n) = 1+ |n|$.
Then for either choice of $\w$, we have
$$  \s(n+r)\leq 1 + |n| + |r| \leq \w(n) \s(r).$$
So $(M, \s, H)$ is a scaled $(G, \w)$-space.  However, if we choose the
weight $\w \equiv 1$, then
$(M, \s, H)$ is {\it not} a scaled $(G, \w)$-space.
\vskip\baselineskip
\par
Let $(M, \s, H)$ be a scaled $(G, \w)$-space.
We show that the action of $G$
on $\schwmh$ is  $m$-$\w$-tempered.
For this, by Theorem 4.6 it suffices to
show that the action of $G$ on
$\schwcm$ is  $m$-$\w$-tempered.
By the estimate (5.8) with
$C$ replaced by $C\w^{l}(g)$, we have
$$ \pa \alpha_{g}(f)
\pa_{d}\leq
C\w^{l}(g)\pa f\pa_{d}.
$$
Since each seminorm $\pa \quad \pa_{d}$ is submultiplicative
(recall we are assuming $\s\geq 1$), this shows that the action of
$G$ on $\schwcm$ is $m$-$\w$-tempered.
\par
We may form the convolution algebra $G\rtimes^{\w}\schwmh$,
which is topologized by the seminorms
$$ \pa \varphi\pa_{d,\ga,l,\be} =
\int_{G}\sup_{ m\in M} \biggl(
\w^{d}(g) \s^{l}(m) |X^{\ga} {\tilde X}^{\be}
\varphi(g,m)|\biggr) dg, \tag 5.15 $$
where $X^\ga$ acts on the first argument of $F$, and
${\tilde X}^{\be}$ on the second.  We place the usual
operations of multiplication and involution on $G\rtimes^{\w} \schwmh$:
$$\aligned \varphi*\psi(g,m) &=
\int_{G}\varphi(h,m)\psi(h^{-1}g, h^{-1}m)dh  \\
\varphi^{*}(g, m) & =
{\overline \varphi}(g^{-1}, g^{-1}m)\Delta (g).
\endaligned  \tag 5.16$$
\proclaim {Theorem 5.17} Let $(M, \s, H)$ be a
scaled $(G, \w)$-space.
Then the space $G\rtimes^{\w} \schwmh$ is a dense  Fr\'echet
subalgebra of $L^{1}(G, C_{0}(M))$, and of
the C*-crossed products $\GC$ and $G\rtimes_{r}C_{0}(M)$.
It is $m$-convex if $\w$ is $m$-sub-polynomial, and is a *-algebra if
$\w \thicksim \w_{-}$. (These last
two conditions are satisfied if $\w$ is
a weight or gauge on $G$.)
\endproclaim
\demo{Proof}
We have seen that the action of $G$
on $\schwmh$ is $m$-$\w$-tempered
and differentiable.
Hence by  Theorem 2.2.6, we have the first statement
of the theorem.  Since $G\rtimes^{\w} \schwmh$ contains
$C_{c}^{\infty}(G, \schwmh)$, we have the density statement.
For the  $m$-convexity, see Theorem 3.1.7.  See Corollary 4.9 for
the last statement.
\qed\enddemo
\vskip\baselineskip
\par
Now
we give several examples of scaled $(G, \w)$-spaces.
\subheading{Example 5.18} In Example 5.14 with $\w(n)= 1+|n|$,
the resulting crossed product
$G\rtimes \schwmns$ is $\Z\rtimes \schwr = \Cal S(\Z, \schwr)$, the
standard set of Schwartz functions from $\Z$ to $\schwr$.
\subheading{Example 5.19} We note that the scale and weights
we choose may depend on the action of $G$ on $M$.  For example,
let $G=H= \R$, $M=\R$, and let $G$ act by translation.  Then the
weight $\w(r) = 1+ |r|$ and scale $\s (m) = 1 + |m|$
make $(M, \s)$ a scaled $(G, \w)$-space.  However, if the action of $G$
on $M$ is given by $r(m) = e^{r}m$, they do not.
If we define a weight $\delta$ on $G$ by $\delta(r) = e^{|r|}$,
then $(M, \s)$ is a scaled $(G, \delta)$-space
for this new action.
\subheading{Example 5.20} Let $G=H$ be the $ax + b$ group,
and let $M=\R$.  Let $G$ act on $M$ by
$$ \pmatrix e^{g_{1}} & g_{2} \\ 0 & 1 \endpmatrix
m = e^{g_{1}}m + g_{2} \tag 5.21 $$
Let $\w$ be the  weight $\w(g)= e^{|g_{1}|}
+ |e^{-g_{1}}g_{2}| + |g_{2}| +1 $
on $G$ of Example 1.6.1.  Let $\s(r) = 1+|r|$ be the scale on $\R$.
Then
$$ \aligned \s(gm) & = (1+ |e^{g_{1}}m + g_{2}|)
\leq (1 + \w(g)|m| + \w(g)) \\
&\leq 2 \w(g) (1+ |m|)
= 2 \w(g)\s(m) , \endaligned \tag 5.22 $$
so $(M, \s)$ is a scaled $(G, \w)$-space.
The resulting crossed product  $G\rtimes \schwr$ is
isomorphic as a Fr\'echet space to $\schwL
\widehat \otimes_{\pi} \schwr$, where $\schwr$ denotes the standard
space of
Schwartz functions on $\R$,
$\schwL$ is the group convolution
algebra of Example 1.6.1, and $\widehat
\otimes_{\pi}$ denotes the completed projective tensor product.
\subheading {Example 5.23}
We let $G=H$ be any closed subgroup of $GL(n, \R)$, and let
 $\w(A) = \max(\pa A \pa, \pa A^{-1} \pa )$
be the weight  on $G$ from Example 1.6.10.
We consider two $G$-spaces.
\roster
\item The group $G$ acts on $\R^{n}$ by left multiplication.
\item The group $G$ acts on $M(n, \R)$ by conjugation.
\endroster
We define a scale $\s$ on $M(n, \R)$ by setting
$\s (S) = 1+ \pa S \pa$, where
$\pa S \pa$ denotes the operator norm of the
matrix $S$.  We define
$\s$ on $\R^{n}$ in a similar way.
In each case $(M, \s)$ is a scaled $(G, \w)$-space.  We verify this
for $M=M(n, \R)$.  Let $A\in G$ and $S\in M$.  Then
$$ \split \s(ASA^{-1}) = & 1+ \pa ASA^{-1} \pa
 \\
&\leq 1 +\pa A \pa \pa S \pa
\pa A^{-1} \pa
\leq  \pa A \pa \pa A^{-1} \pa
\s (S). \endsplit $$
 Since $\w^{2}(A) \geq \pa A \pa \pa A^{-1} \pa$,
 this shows that $(M, \s)$
is a scaled $(G, \w)$-space.
\subheading{Example 5.24}
Let $G$ be a Lie group with weight $\w$ that bounds $Ad$ on $G$.
Let $N$ be a subgroup of $G$ and let $(M, \s)$
be a scaled $(N, \w)$-space.
Define an equivalence relation on $G\times M$
by $(g, m) \sim ({\tilde g}, {\tilde m})$ if and only if
there is some $n\in N$ such that ${\tilde g} = gn^{-1}$ and
${\tilde m}=nm$.  Let $M_{G}= G\times_{N} M$ denote
the quotient space.  Then $G$ acts as a transformation group on $M_{G}$.
The $G$-space $M_{G}$ is {\it induced} from the $N$-space $M$.
\par
We define a scale $\s_{G}$ on $M_{G}$ by
$$\s_{G}([g, m)]) = \inf_{n\in N} \w(gn^{-1})\s(nm).\tag 5.25$$
Then	$\s_{G}({\tilde g}[(g, m)]) \leq
\w({\tilde g}) \s_{G}([(g, m)])$ so $(M_{G}, \s_{G})$
is a scaled $(G, \w)$-space.  We let $\Cal S(M_{G})$ denote
the corresponding $G$-Schwartz functions.
\par
For a specific example, let $G=\R$, $N= \Z$ and let $M$ be the circle
$\T$.  Let $\Z$ act on $\T$ by an irrational rotation $\theta$. Let
$\w(r)=1+|r|$ be the weight on $\R$ and $\s \equiv 1$
be the constant scale on $\T$.  Then $(M_{\R} = \R \times_{\Z} \T,
\s_{\R})$ is a scaled $(\R, \w)$-space. We may form $\Cal S(M_{\R})$.
   Note that for $[(r, t)] \in M_{\R}$ we have
$\s_{\R}([(r, t)])= \inf_{n\in \Z} \w(r-n)\s(n\theta + t)
= \inf_{n\in \Z} \w(r-n) \leq 2$.
\subheading{Example 5.26} The scale $\s_{0} \equiv 1$ on $M$
always makes $(M, \s_{0})$ a scaled $(G, \w)$-space, for any
group $G$ and weight $\w$ which bounds $Ad$ on $G$. If $M$ is a
scaled $H$-space, with $H$ properly containing $G$, then
$(M, \s_{0}, H)$ is a scaled $(G, \w)$ space as long as $\w$ bounds
$Ad$ on $ H$.
\subheading{Example 5.27}  To see some cases where $A_{}$ is
noncommutative, we could in general take $A_{}=B$, or $A_{}= B^{\infty}$,
where $B^{\infty}$ denotes the set of $C^{\infty}$-vectors for the
action of a Lie group $G$ on a C*-algebra $B$.
Let $\w$ be any weight on $G$ which bounds $Ad$.  If $A_{}=B$,
$$\pa \alpha_{g} (a) \pa = \pa a \pa,
\tag 5.28 $$
so the action of $G$ is always tempered and $m$-tempered.
We saw in Theorem 4.6 that if
$A_{} = B^{\infty}$, then the action of $G$ on $A$ is still tempered
and $m$-tempered.
So in either case, we may form the $m$-convex Fr\'echet
*-algebra $G\rtimes^{\w} A$.
We say that $G$ is an elementary Abelian group
if it is of the form $\T^{n} \times \R^{m} \times \Z^{p} \times F$,
where $F$ is a finite Abelian group.
If $G$ is an elementary Abelian group, $A=B^{\infty}$,
and $\w$ is any weight equivalent to the word gauge on $G$,
then the algebra $G\rtimes^{\w} A$ is precisely
the dense subalgebra defined in \cite{Bo, \S 2.1.4}
with $E=A$. (We shall see in Theorem 6.8 and Proposition 6.13(2) below
that although
Bost uses the sup norm and I use the $L^{1}$ norm,
we still get the same algebra.)
\par
For a specific example where $B$ is noncommutative,
we let $B$ be the two dimensional irrational rotation
algebra $A_{\theta}$.
Here we think of $A_{\theta}$ as the C*-algebra generated
by two unitaries $U$, $V$ with commutation relations
$$  UV= e^{-2\pi i\theta}VU   $$
We let $G=SL(2, \Z)$ and $\w$ be any weight on $G$.
Let
$$g= \pmatrix a & b \\ c&d\endpmatrix \in SL(2,\Z).
$$
Then $G$ acts on $A_{\theta}$ by
$$ \alpha_{g}(U^{n}V^{m})= \xi(g, n, m)U^{an+bm}V^{cn+dm}  $$
where $\xi $ is a map to $\T$ given by
$$\xi (g, n, m) = e^{(nm - (an+bm)(cn+dm))i\pi \theta/2} $$
We thus may form the dense $m$-convex Fr\'echet *-subalgebra
$SL(2, \Z) \rtimes^{\w} A_{\theta}^{\infty}$
($=\Cal S_{1}^{\w}(SL(2, \Z),
A_{\theta})$)
of the C*-algebra $SL(2, \Z) \rtimes A_{\theta}$.
\vskip\baselineskip
\heading  \S 6 Sup Norm and $L^{r}$ norm
in Place of $L^{1}$ norm, and nuclearity  \endheading
\par
We define the $L^{r}$  Schwartz functions
on $G$ for $1\leq r\leq \infty$.  We denote these by
$\solp$.
We show that $\solp$ is always
contained in $\soinf$ with continuous inclusion map.
If the reciprocal of the scale is in $L^{p}(G)$
for some $p\geq 1$, we show that
$\schwsig = \solp = \soinf$.
We show that in many of
our examples, this condition
is satisfied.
We compare our algebras with those in \cite{Jo}, \cite{Ji}, and with
Rader's algebra \cite{War, Prop 8.3.7.14} and the zero Schwartz space
for $SL(2, \R)$ \cite{Bar, \S 19}.
Finally we prove that if the above
integrability condition  on $1/\s$ is satisfied,
then $\solp$ is a nuclear Fr\'echet space.
\subheading {Definition 6.1}
We defined $\soinf$ in \S 5.  Let $1\leq r <\infty$.
Let $\s\geq 1$ be a scale on a Lie group $G$.
Let $\solp$ be the set of differentiable functions $\psi \colon
G \rightarrow \C$ satisfying
$$ \pa \psi \pa^{(r)}_{m,\ga} =
\pa \s^{m}X^{\ga}\psi \pa_{r} =
\biggl(\int_{G} |\s^{m} X^{\ga} \psi(g) |^{r} dg\biggr)^{1/r}
 <\infty \tag 6.2$$
for each $\ga \in \N^{q}$ and $m\in \N$.
We call
$\solp$ the {\it $\s$-rapidly vanishing ($L^{r}$)  Schwartz
functions on $G$}.
We topologize
$\solp$ by the seminorms (6.2).
It is easily checked that if $\s_{1} \thicksim \s_{2}$, then
$\Cal S_{r}^{\s_{1}}(G)$ is isomorphic to $\Cal S_{r}^{\s_{2}}(G)$.
\subheading{Definition 6.3} Let $\s \geq 1$
be any scale, and let $r\geq 1$.
Define the {\it $\s$-rapidly vanishing $L^{r}$ functions
$L_{r}^{\s}(G)$ on $G$} to be the
space of Borel measurable functions
$f\colon G \rightarrow \R$ such that
$$ \pa f \pa_{d}^{(r)} = \biggl(\int_{G} \,|\,\s^{d} f(g) \, |^{r} \, dg
\biggr)^{1/r}\tag 6.4$$
is finite for each $d \in \N$.
Then $L_{r}^{\s}(G)$ is complete for the topology given
by the seminorms $\pa \quad \pa_{d}^{(r)}$
\cite{Schw, \S 5}.  Just as in \S 1
for $r=1$, we have the following theorem.
\proclaim{Theorem 6.5}  Let $\s$ be a
translationally equivalent scale (for
example a weight or a gauge) on a Lie group $G$.
Then the action $\alpha_{g}(\varphi)(h) = \varphi(g^{-1}h)$ of
$G$ on $\solp$ and $L_{r}^{\s}(G)$
by left translation is well defined,
strongly continuous, and differentiable on $\solp$.  The Fr\'echet
space of $C^{\infty}$-vectors $L_{r}^{\s}(G)^{\infty}$ is naturally
isomorphic to $\solp$.  Hence
$\solp$ is complete and a Fr\'echet space.
\endproclaim
\demo{Proof}
The proof is the same as for
the case $r=1$ - see Theorems 1.2.21 and 2.1.5.
\qed
\enddemo
\par
The Schwartz functions $\solp$
may not be algebras if $1<r\leq \infty$, even if $\s$ is a weight or a
gauge.  For example, let $\s\equiv 1$.  Then $\schwsig$
is an algebra.
We have $\soinf= C_{0}^{\infty}(G)$,
where $C_{0}^{\infty}(G)$ is the set of differentiable functions
$\psi$ such that $\psi$ and all its derivatives vanish at infinity.
In this case, $\soinf$ is rarely an algebra for convolution.
For example, take $G=\Z$.
  Let $\psi \in c_{0}(\Z)$ be defined by
$$\psi (n)= {1\over{ (1+|n|)^{1/2}}}.\tag 6.6$$
Then
$$\aligned \psi *\psi (0) &= \sum_{m}\psi(m) \psi(0 -m) \\
& = \sum_{m} {1\over { (1+|m|)^{1/2}}}
{1\over { (1+|-m|)^{1/2}}}\\
& = \sum_{m} {1\over{1+|m|}}, \endaligned \tag 6.7 $$
which diverges.  (Note that in this case, we do not have $\soinf
\subseteq C^{*}(G)$ or $\subseteq L^{1}(G)$.)
\par
For an example when $\solp$ is not an algebra for $r<\infty$,
let $r=2$ and $G=\Z$.
Then $\solto= l^{2}(\Z)$ which is not a convolution
algebra. (To see this, note that $l^{2}(\Z)\cong L^{2}(\T)$
via the Fourier transform, which is not closed under
pointwise multiplication.)
\proclaim{Theorem 6.8} Let $\s\geq 1$ be a
translationally equivalent scale
on a Lie Group $G$.
Then  $\solp$ is contained in $\soinf$ for any $r\geq 1$,
with continuous inclusion map.
\par
Assume in addition that there is some $p>0$ such that
$$ \int_{G} {1\over \s^{p}(g)} dg < \infty.\tag 6.9 $$
Then $\solp =\schwsig$ as Fr\'echet spaces for all $r \in [1, \infty]$.
\endproclaim
\par
If $\s$ satisfies condition (6.9), we say that {\it $1/\s$
is in some $L^{p}$-space}.  Note that condition (1.5.2) is
equivalent to (6.9).
\demo{Proof} By Theorem 6.5 and \cite{DM, Thm 3.3}, we may write any
element of $\solp$ as a finite sum of functions of the form
$$ f*\psi(g)=  \int f(h) \psi(h^{-1}g) dh, $$
where $f \in C_{c}^{\infty}(G)$ and $\psi \in \solp$.
By the chain rule and
uniform translational equivalence of the scale (see Theorem 1.2.11),
$$ \aligned | \s^{m} X^{\ga}(f*\psi)(g) |
&\leq \sum_{\be\leq \ga}\int \, |\, \s^{m}(g)
p_{\be}((Ad_{h^{-1}})_{ij}) f(h)
 X^{\be}\psi (h^{-1}g)\, | \, dh \quad{(2.2.3)} \\
&\leq \sum_{\be\leq \ga}{C}\int\, | \,
p_{\be}((Ad_{h^{-1}})_{ij}) f(h)
\s^{mk}(h^{-1}g) X^{\be}\psi (h^{-1}g)\, | \, dh \\
& \leq D \sum_{\be\leq \ga}
\biggl(\int \, | \, \s^{mk}X^{\be} \psi (h) \, |^{r} \,dh \biggr)^{1/r}
\\
& = { D} \sum_{\be\leq \ga}
\pa \psi \pa_{mk, \be}^{(r)}, \endaligned \tag 6.10
$$
where the $p_{\be}((Ad_{h^{-1}})_{ij})$ are polynomials
in the matrix entries  of $Ad_{h^{-1}}$.
So we have
$\solp \subseteq \soinf$.
\par
To see that the inclusion maps are continuous, we use the closed
graph theorem just as in Remark 1.2.2 of \cite{Jo}.
Let $\psi_{n}\longrightarrow 0$ in $\solp$, and assume
that $\psi_{n}\longrightarrow \psi$ in $\soinf$.
Assume for a contradiction that $\psi \not= 0$.
Then for sufficiently large $n$, each $\psi_{n}$ must have
absolute value bigger than or equal to some $\delta>0$
in a fixed neighborhood $U$ in $G$.  But this contradicts
$\psi_{n} \longrightarrow 0$ in $\solp$.
So the graph of the inclusion map
$\solto \hookrightarrow \soinf$
must be closed, and
thus the inclusion map is continuous.
This proves the first statement of the theorem.
\par
For the second statement of the theorem, we  use condition
(6.9) to show that $\soinf \subseteq \solp$
if $r<\infty$, with continuous inclusion maps.
We let $\psi \in \soinf$.  Then
$$\aligned \pa \s^{m} X^{\ga}\psi \pa_{r} & =
\pa  {1\over{\s^{p}}}
\s^{m+p} X^{\ga} \psi \pa_{r} \\
&
\leq
\pa  {1\over{\s^{p}}} \pa_{r}
\pa \s^{m+p} X^{\ga}\psi\pa_{\infty} \\
& \leq C^{1/r} \pa \s^{m+p}X^{\ga}\psi \pa_{\infty},
\endaligned \tag 6.11 $$
where $C$ is the value of the integral in (6.9).
So the seminorms on $\solp$ are dominated
by the seminorms on $\soinf$.
Thus we have the desired continuous inclusions. From
the first part of the theorem, it follows
that $\solp= \so $ for any $r\in [1, \infty]$, all with equivalent
topologies.
\qed \enddemo
\vskip\baselineskip
\proclaim{Corollary 6.12} If $\s$ is any sub-polynomial scale on
a Lie Group $G$ and
the identity component $G_{0}$ of $G$
is non-trivial, then the group algebra
$\Cal S^{\s}_{1}(G)$ has no (left or right) bounded approximate unit.
\endproclaim
\demo{Proof}  By Theorem 6.8, any bounded approximate unit must
be bounded in the sup norm.  But with convolution multiplication,
any approximate unit $e_{n}$ must become unbounded near the
identity as $n\rightarrow \infty$.
\qed
\enddemo
\par
We show that for many of our examples
of gauges and weights, the integrability
condition
(6.9)
is satisfied.  Part (2) of the
following proposition is proved in \cite{Ji}
for discrete groups, and in Proposition 1.5.1 above.
\proclaim {Proposition 6.13} In each of the following cases,
the reciprocal of the scale $\s$ on $G$ is in
some $L^{p}$-space.
\roster
\item The group $G$ is  compactly generated and $\s$ is the
exponentiated word weight on $G$.
\item The group $G$ is compactly generated, of polynomial growth,
and $\s(g)= 1+\tau(g)$, where $\tau$ is the word gauge. (Moreover,
a gauge which satisfies the integrability condition exists on
a compactly generated group $G$ iff $G$ has polynomial growth.)
\item The group $G$ is $GL(n, \R)$ and $\s(g) = \max(\pa g \pa, \pa
g^{-1} \pa )$.
\item The group $G$ is a closed subgroup of a Lie group $H$,
and $\s$ is the restriction to $G$ of
an $H$-translationally equivalent scale (also denoted by $\s$)
on $H$, such that
both $\s$ and $\s_{-}$
bound $Ad$ on $H$, and $1/\s \in L^{p}(H)$.
(For example, $\s$ could be any weight or gauge on $H$ which bounds
$Ad$ and satisfies $1/\s \in L^{p}(H)$.)
\endroster
\endproclaim
\demo{Proof}
Part (2) is part of
Proposition 1.5.1.  We prove part (1).  Let $U$ be an open rel. comp.
generating set for a compactly generated group $G$.  Then $|U^{n}|$
grows at most exponentially.  We see this as follows (compare \cite{Je,
p.114-
115}).  Since $U^{2}$ is rel. comp., and  $U$ is open,
there is some $l\geq 1$ and finitely many
$u_{1}, \dots u_{k} \in U^{l}$ such that
$U^{2} \subseteq \cup_{i=1}^{k} u_{i}U$.  It follows that
$$ U^{n} \subseteq \cup_{i_{1}, \dots i_{n} =1}^{k} u_{i_{1}}\dots
u_{i_{n}} U.  $$
Hence the Haar measure of $U^{n}$ is bounded by $k^{n}$ times the
Haar measure of $U$.
\par
So let $p\in \N$ be such that the Haar measure of $U^{n}$
is bounded by $e^{n(p-1)}$.  Let $\tau$ be the word gauge, and
$\w= e^{\tau}$ be the exponentiated word weight.  We have
$$\split \int_{G} {1\over{\w^{p}}}dg & = \int_{G} e^{-p\tau(g)} dg
\leq \sum_{n=0}^{\infty}
\int_{U^{n}} e^{-pn} dg\\ &= \sum_{n=0}^{\infty}\,|U^{n}|\, e^{-pn}
\leq \sum_{n=0}^{\infty} e^{-n} <\infty. \endsplit $$
This proves (1).
\par
For (3), use (1) and Example 1.6.10.
For (4), assume that $\s$ is defined on a Lie group $H$
which contains $G$ as a closed subgroup, and assume that
$\s$ and $\s_{-}$ bound $Ad$ on $H$ and $1/\s\in L^{p}(H)$.
Then $\s$ bounds $Ad$ on $G$ (see Example 1.3.15) so $\s$
bounds the  modular function $\Delta_{G}$
(see (4.13) and preceding remarks).  Also, $\s_{-}$ bounds
the modular function $\Delta_{H}$.
Let $C>0$ and $d\in \N$ be such that
$$\Delta_{G}(g) \leq C\s^{d}(g) \qquad g \in G, $$
and
$$\Delta_{H}(h^{-1}) \leq C\s^{d}(h) \qquad h \in H. $$
Since $1/\s \in L^{p}(H)$, by \cite{Bou, chap VII, \S 2, thm 2}
we have
$$ \int_{G} {1\over{\s^{p}(hg)}} \Xi(g) dg
< \infty \tag 6.14$$
for almost every $h\in H$, where $\Xi(g) =
{\Delta_{H}(g)/{\Delta_{G}(g)}}$.
  Fix some $h\in H$ such that (6.14)
is true.  Then since $\Xi(g) \geq 1/(C^{2}\s^{2d}(g))$, we have
$$ \aligned \infty &> \int_{G}
{1\over{\s^{p}(hg)}} \Xi(g) dg
\\
& \geq C_{h}
\int_{G}
{1\over{\s^{pk}(g)}}\Xi(g)dg
\qquad {\text{trans eq of $\s$}}\\
& \geq D
\int_{G}
{1\over{\s^{pk+2d}(g)}}dg
\endaligned $$
for some positive constant $D$.
Hence $1/\s \in L^{pk+2d}(G)$.
This proves Proposition 6.13. \qed \enddemo
\subheading{Remark 6.15} Since every nilpotent Lie group has polynomial
growth \cite{Pa, Cor 6.18}, it follows from part (2)
 of Proposition 6.13 that there exists a gauge on any closed
subgroup of a connected nilpotent Lie group which satisfies
the integrability condition (6.9)(or (1.5.2)).
By parts (3) and (4) of Proposition 6.13, we see that the weight
$\theta$ on  $GL(n, \R)$
defined in Example 1.6.10 satisfies (6.9), and also  that the
restriction of this weight to $SL(n, \R)$
(see Example 1.6.4)   or to any closed subgroup of
$GL(n, \R)$ satisfies (6.9).
It also follows from Proposition 6.13, or by simple calculations, that
condition (6.9)  is satisfied in Examples 1.3.14, 1.6.1, 5.20, and 5.23.
Thus it makes no difference in these examples  which
space $\solp$ of Schwartz functions we use.
\subheading{Remark 6.16} We compare our algebras with those in \cite{Jo}.
(The algebras in \cite{Ji} are special cases of these.)
If $G$ is amenable and discrete and $\s$ is a gauge on $G$,
then an easy argument
shows that $\solto\subseteq C_{r}^{*}(G)$ if and only if
$\solto \subseteq l^{1}(G)$ (see proof of Corollary 3.1.8 of \cite{Jo}).
Also, if $\solto \subseteq l^{1}(G)$, then $\solto \subseteq \so$
since $\s^{m}\psi(g)\in l^{1}(G)$ for all $\psi \in \solto$.
So if $G$ is amenable, $\solto$ is a subalgebra of $C^{*}_{r}(G)$
if and only if $\solto = \so$.
\par
If $G$ is not amenable, however, it is not clear to me if
$\solto \subseteq C_{r}^{*}(G)$ implies
$\solto \subseteq l^{1}(G)$.  So, as far as I know,
there could be cases when $\solto$ is a dense subalgebra of
$C_{r}^{*}(G)$ and $\so \not= \solto$.
In such a case, Theorem 3.1.7 above does not directly imply the
$m$-convexity of $\solto$.  However, the $m$-convexity
of $\solto$ (and in fact all the dense subalgebras in
\cite{Jo}, including the appendix) is clear from Theorem 3.1.4
above (or \cite{Mi, Prop 4.3})
and the estimates of \cite{Jo, Lemma 1.2.4}.
\subheading{Remark 6.17}
Theorem 3.1.7 of \cite{Jo} says that for a discrete
group
where $\s$ is the word gauge, $\solto\subseteq l^{1}(G)$
if and only if $G$ is of polynomial growth.
But we saw in Remark 6.16 that
$\solto \subseteq l^{1}(G)$ if and only if
$\solto \subseteq \so$, so  $\so=\solto$ if and only if
$G$ has polynomial growth.
If we drop the hypothesis that $\s$
be a gauge, we can
have $\so=\solto$ even when $G$ does not have polynomial growth.
For an example, let $G$
be the free group on two generators, and $\s(g)= e^{|g|}$ where
$|g|$ denotes the word length of $g$.
Apply Proposition 6.13(1) and Theorem 6.8.
\subheading{Remark 6.18} We compare our group algebras
with Rader's algebra $\Cal C (G)$ \cite{War, vol. II, Prop 8.3.7.14}.
Here $G$ is a reductive Lie group - see \cite{War, vol. II, \S 8.3.7}.
A function $f$ is in $\Cal C (G)$ if and only if
$$ {}_{D_{1}}|f|_{r, D_{2}} = \sup_{x\in G}(1+\s (x))^{r}
\Xi^{-1}(x) |D_{1} fD_{2}(x)|<\infty, \tag 6.19 $$
where $D_{1}$ is a differential operator acting
on $f$ on the left, and $D_{2}$ is a differential
operator acting on $f$ on the right.  The Riemannian symmetric function
$\s$ is a gauge (in fact equivalent to the
word gauge by Example 1.6.4)  on $G$
and $\Xi$ is the zonal spherical function on $G$.
The space $\Cal C(G)$ is then topologized
by the seminorms (6.19).
\par
To contrast this with the situation in this paper,
we consider the case $G=SL(2, \R)$ (see Example 1.6.4).
Let $sl(2, \R)$ be the Lie algebra of $SL(2, \R)$\cite{War}.
Let $\pa \quad \pa_{sl}$ denote the norm on
$sl(2, \R)$ used in \cite{HW, \S 2}. (The group $SL(2, \R)$ is a
special case of the groups considered there.)
Let
$$g=\pmatrix e^{a} & 0 \\ 0 & e^{-a} \endpmatrix.$$
Then by \cite{HW, (2.2a)}, we have
$$ \s(g)= \s\biggl( exp \pmatrix a & 0 \\ 0 & -a \endpmatrix
\biggr) = \pa \pmatrix a & 0 \\ 0 & -a \endpmatrix
\pa_{sl} \propto
 | a |. $$
In  (1.6.9), we saw that
$$ Ad_{g} = \pmatrix 1 & 0 & 0 \\ 0 & e^{2a} & 0 \\
0 & 0 & e^{-2a} \endpmatrix.$$
So $\s$ does not bound $Ad$.
So in order to get a *-algebra, we would use some scale other than
$\s$ to form our Schwartz algebra.
As we saw in Example 1.6.4,
the weight $\theta(x)= \max(\pa x \pa, \pa x^{-1} \pa )$
defined on $SL(2, \R)$
does bound $Ad$. (In fact, we saw that it is equivalent
to the exponentiated word weight.)
Here $\pa \quad \pa$ denotes the norm on a matrix of $SL(2, \R)$.
Our weight $\theta$ grows like the sum of the
absolute values of the matrix entries, and thus exponentially
faster than $\s$.  Since our weight bounds $Ad$,
we would not gain anything by having differential operators
act on both sides as in (6.19).  Using Theorem 6.8,
our group algebra is the set
$$\Cal S^{\theta}_{\infty}(G) = \{\, f\,|\, \sup_{x\in G}
\theta^{r}(x)|(Df)(x)| <\infty, \quad
r\in \N, \quad D {\text { a diff op }}\}. \tag 6.20 $$
Another difference between the two sets of functions
(6.20) and $\Cal C (G)$ is that the seminorms (6.19) for
$\Cal C(G)$ make use of the zonal spherical
function $\Xi$,  whereas the definition of our space $\Cal
S^{\theta}_{\infty}(G)$
does not allow that.
\par
As we noted in Remark 6.15, the  weight $\theta$
satisfies the integrability condition (6.9).
To have this condition make sense for $\s$,
we must consider $1+\s$, since $\s(I)=0$ (see (1.5.2)).
The function $1+\s(x)$ does not satisfy (6.9),
since otherwise $SL(2, \R)$ would have
polynomial growth by Proposition 6.13,
part 2 (or 1.5.1).
And $SL(2, \R)$ does not have polynomial growth \cite{Je, Prop 1.3,
Thm 1.4}.
\par
We note that Rader's algebra is $m$-convex.  By the estimate
just preceding Prop 8.3.7.14 of \cite{War, vol II},
we have
$$ {}_{D_{1}}|f*g|_{s, D_{2}} \leq
M {}_{D_{1}}|f|_{r}|g|_{s, D_{2}} \tag 6.21 $$
for some constants $M$ and $r\in \N$.  Hence the
seminorms
$$ \pa f \pa_{l} = \max_{|D_{1}|, |D_{2}| \leq l}
{}_{D_{1}}|f|_{l, D_{2}}
\tag 6.22
$$
satisfy the condition of Theorem 3.1.4 above and $\Cal C(G)$ is
$m$-convex. (In fact, the $m$-convexity follows from
(6.22) and \cite {Mi, Prop 4.3}.)
\par
We show that our group algebra
$\Cal S^{\theta}(G)$ for $SL(2, \R)$ is in fact
the zero Schwartz space defined in \cite{Bar, \S 19}.
By \cite{Bar, (19.1)}, the zero Schwartz space $\Cal C^{0}(G)$
is topologized by the seminorms
$$ \rho^{0}_{D_{1}, D_{2}, l} (f) = \sup_{x\in G}|e^{(1+l)\s(x)}
D_{1}fD_{2}(x)|,\tag 6.23
$$
where $D_{1}$ and $D_{2}$ are differential operators on the right
and left respectively, and $l\in \N$.
Recall from Example 1.6.4 that $\theta \thicksim e^{\s}$.
Thus our Schwartz space $\Cal S^{\theta}(G)$
above is topologized by the seminorms
$$  \pa f \pa_{l, \ga} = \sup_{x \in G} |e^{(1+l)\s(x)} X^{\ga}f(x)|,$$
where $l \in \N$.
These two topologies are clearly the same, except that
(6.23) uses differential operators both on the
left and the right.  But since $e^{\s}$ bounds $Ad$,
this makes no difference.
Thus $\Cal C^{0}(G)= \Cal S^{\theta}(G)$.
\vskip\baselineskip
%Begin NUCLEARITY
\par
 We show that if the integrability condition
(6.9) holds, then $\schwsi$ has the  property
of being a nuclear Fr\'echet space\cite{Tr}.
If (6.9) fails to hold, then
$\so$ may fail to be nuclear.  For example, $l^{1}(\Z )$
is an infinite dimensional Banach space and hence not
nuclear \cite{Tr, Cor 50.2}.
\proclaim {Theorem 6.24} Let $\s$ be a
translationally equivalent scale on a Lie group $G$
which satisfies the integrability condition (6.9).
Then $\so = \solp$ for all $r\in [1, \infty]$ and
all of these spaces are nuclear Fr\'echet spaces.
\endproclaim
\demo{Proof}
The first statement is Theorem 6.8.
For the second, we imitate the proof of Theorem 6.2.5 of \cite{Pi},
which gives the nuclearity for the standard set of Schwartz
functions on $\R$.
\par
We first recall some general elementary facts and definitions
on nuclearity from \cite{Pi}.  Let $E$ be any locally
convex space.  If $U$ is a neighborhood of zero in $E$
(following \cite{Pi}, we do not assume neighborhoods are open), then
the {\it polar} $U^{0}$ of $U$ denotes the weakly compact
subset
$$  \{ \varphi \in E^{\prime} |
\quad | \varphi (u)|\leq 1, \quad u\in U\}
\tag 6.25 $$
of the topological dual $E^{\prime}$ of $E$. (The space $E^{\prime}$
consists of all linear functionals on $E$  which are continuous
for the locally convex topology on $E$.)  The seminorm $p_{U}$
associated with $U$ is defined by
$$ p_{U}(x) = \inf \{  \rho > 0 |  x\in \rho  U \} \tag 6.26 $$
\proclaim{Lemma 6.27 \cite{Pi, Prop 4.1.5}}
The locally convex space $E$ is nuclear if and only if
 some fundamental system $\{ U_{\alpha} \}$ of
zero neighborhoods for $E$ has the property
that for all $U\in \{ U_{\alpha} \}$ there is some
$V \in \{ U_{\alpha} \}$  and a positive Radon measure
$\mu$ defined on the polar $V^{0}$ such that
$$ p_{U}(x) \leq \int_{V^{0}}|\varphi^{\prime} (x)\,|\, d\mu
(\varphi^{\prime})
\tag 6.28 $$
for all $x \in E$.  \endproclaim
\par
We use the notation $\schw$ for $\so = \soinf$.  Define
seminorms $\pa \quad \pa_{l}^{1} $ and
$\pa \quad \pa_{l}^{\infty}$ on $\schw$ by
$$\pa \psi \pa_{l}^{i} = \sum_{|\be|\leq l}
\pa \s^{l}X^{\be}\psi \pa_{i}, \tag 6.29 $$
where $i = 1, \infty$.
These give equivalent topologies by assumption (6.9) and Theorem 6.8.
We show that the condition for nuclearity given in
Lemma 6.27 holds for the fundamental system
of $L^{1}$ neighborhoods of zero $\{ \psi \in \schw
| \pa \psi \pa_{l}^{1}\leq 1 \}$.
Let $U_{l}$ be the zero neighborhood
$$U_{l}= \{ \psi \in \schw
| \pa \psi \pa_{l}^{1}\leq 1 \}.  $$
(Then $p_{U_{l}}(\psi ) = \pa \psi \pa_{l}^{1}$.)
We let $V_{l+p}$ be any $L^{1}$ neighborhood
of zero contained in the sup norm neighborhood
$$  \{ \psi\in \schw
| \pa \psi \pa_{l+p}^{\infty}\leq 1 \},  \tag 6.30$$
where $p$ is from condition (6.9).
\par
For each $|\be|\leq l$ and $g\in G$, the formula
$$ \epsilon_{g}^{[\be ]} (\psi ) = \s^{l+p} (g)
(X^{\be}\psi )(g)  \tag 6.31  $$
defines a  function from $ \schw$ to $\C$,  continuous
in $\psi$.  Note that
$\epsilon_{g}^{[\be ]} $ lies in the polar set $V_{l+p}^{0}$ of
$V_{l+p}$, since it is in the polar set of (6.30).
\par
Define a positive Radon measure $\mu$ on $V_{l+p}^{0}$ by the
equation
$$  \int_{V_{l+p}^{0}} F(\varphi^{\prime}) d\mu (\varphi^{\prime}) =
\sum_{|\be|\leq l} \int_{G}
{1\over {\s^{p}(g)}} F(\epsilon_{g}^{[\be]})dg
\tag 6.32 $$
for all continuous functions $F$ in $C(V_{l+p}^{0})$.
(Note that condition
(6.9) implies that this measure $\mu$ is well defined and continuous.)
Finally,  we have
$$\aligned p_{U_{l}}(\psi)=\pa \psi \pa_{l}^{1} &
= \sum_{|\be|\leq l} \pa \s^{l}X^{\be} \psi \pa_{1}
\qquad {\text { def (6.28)}}\\
& =
\sum_{|\be|\leq l}
\pa
{1\over {\s^{p}}} \s^{l+p} X^{\be} \psi \pa_{1}
\\ & = \sum_{|\be|\leq l} \int_{G}
{1\over {\s^{p}(g)}} |\epsilon_{g}^{[\be]} (\psi )| dg
\qquad {\text { def of $\epsilon_{g}^{[\be ]}$}}\\
& = \int_{V_{l+p}^{0}} |\varphi^{\prime} (\psi )| d\mu (\varphi^{\prime})
\qquad {\text { def of $\mu$ (6.32)}}\\
\endaligned  $$
Thus by Lemma 6.27, the Fr\'echet space $\schw$ is nuclear.
This proves Theorem 6.24.  \qed \enddemo
\subheading{Question 6.33} Is $\schwsi$ nuclear if and only if
$\s$ satisfies the integrability condition (6.9)?
(See \cite{Jo, Thm 3.1.7} for the discrete case.)
\par
We prove some consequences of our nuclearity result.
\proclaim{Proposition 6.34} Let $E$ be any Fr\'echet
space. If $\so$ is nuclear, then
$\soneco$ is isomorphic as a Fr\'echet space to the projective completion
$\so \tensp E$. \endproclaim
\demo{Proof} Let $\schw$ denote $\so$, and let
$\Cal S(G, E)$ denote $\sonec$.
  Recall that $\Cal S(G, E)$ is the
set of $C^{\infty}$-vectors for the action (2.1.4) of $G$ on
$L_{1}^{\s}(G, E)$, and $\schw$ is the set of $C^{\infty}$-vectors for the
the same action  of $G$ on $L_{1}^{\s}(G)$. From this and
by Theorem A.8 of the appendix, and the fact that we
have a $G$-equivariant isomorphism $L_{1}^{\s}(G, E)
\cong L_{1}^{\s}(G)\tensp E$\cite{Schw, \S 5}, we always have a continuous
surjection
$$\theta \colon \schw \tensp E \longrightarrow \Cal S(G, E). $$
\par
We show that if $\schw$ is nuclear, then this map is injective.
Let $\tense$ denote the completed $\epsilon$ tensor product \cite{Tr}.
By nuclearity, we have $\schw \tensp E \cong \schw \tense E$.
Assume $\theta(\varphi)=0$ for some $\varphi \in \schw \tense E$.
Let $e', x' $ be continuous linear functionals on $E$ and
$\schw$ respectively.  We show that $(x' \otimes e')(\varphi)=0$.
Let ${\tilde e}' \colon \Cal S(G, E) \rightarrow \schw$
denote the continuous linear map given by ${\tilde e}' (\psi)(g)
= e'(\psi(g))$ for $\psi \in \Cal S(G, E)$.  Then $0 = {\tilde e}'
(\theta(\varphi)) = (id \otimes e')(\varphi)$.
But $(x'\otimes e')(\varphi)
= x'((id \otimes e')(\varphi))$.  Hence $(x'\otimes e')(\varphi)=0$
for all pairs $(x',e')\in \schw' \times E'$.   By the definition
of the completed epsilon tensor product $\schw \tense E$ \cite{Tr,
Def 43.5, 43.1}, this implies that $\varphi = 0$.  Hence $\theta $
is injective and we have proved Proposition 6.34.
\qed \enddemo
\proclaim{Theorem 6.35} Let $(M, \s, H)$ be a scaled $(G, \w)$-space.
 Assume that $1/\w$ is in $L^{p}(G)$ for
some $p$.
  Let
$\schwmns$ denote $\schwmh$, and let $\Cal S(G\times M)$
denote $\Cal S^{\w\times \s}_{G\times H}(G\times M)$.
We then have
$$ \schw \tensp \schwmns \cong \schwgtm. \tag 6.36 $$
By Proposition 6.34 above, we  have $\GSno\cong \schwgtm$,
so that we may view $\GSno$ as a space of $G\times H$-differentiable
$\w\times \s$-rapidly vanishing functions.
\endproclaim
\demo{Proof} We show this by an argument
analogous to \cite{Tr, Thm 51.6},
using the fact that $\schw$ is nuclear.
A simple calculation shows that the map
$$\schw \times \schwmns \longrightarrow \schwgtm $$
$$(\psi, \phi) \mapsto \psi(g)\phi(m)$$
is a continuous bilinear map.  Thus we have an inclusion
$$\schw \otimes \schwmns \hookrightarrow \schwgtm,$$
continuous if we place the projective topology on the uncompleted
tensor product.  We will be done if we show that $\schwgtm$
induces a topology on $\schw \otimes \schwmns$ which
is stronger than (and hence equivalent to)
the projective topology, and that $\schw \otimes
\schwmns$ is dense in $\schwgtm$.
\enddemo
\subheading{Dense}
For the denseness, we use the fact  that
the algebraic tensor product of continuous
functions $C_{c}(G) \otimes C_{c}(M)$
is dense in $C_{c}(G\times M)$ in the
inductive limit topology.  It follows that $C^{\w}(G)\otimes
\schwcm$ is dense in $C^{\w\times \s}(G\times M)$, so
that the canonical map $
C^{\w}(G)\tensp
\schwcm \longrightarrow C^{\w\times \s}(G\times M)$
has dense image.
Using Theorem A.8 and the definitions of our spaces, we
see that $\schw \tensp \schwmns \longrightarrow \schwgtm$
has dense image.  Thus the algebraic tensor product
$\schw \otimes \schwmns$ is dense in $\schwgtm$.
\subheading{Topologies equivalent}
Since $\schw$ is nuclear, the projective
topology on $\schw\otimes \schwmns$
agrees with the $\epsilon$-topology.  So it suffices to show that
if a sequence $\Xi_{n}$ in $\schw\otimes \schwmns$
converges to zero in $\schwgtm$, then it
converges to zero in the $\epsilon$-topology.
For this it suffices to show that if $A\subseteq \schw^{\prime}$
and $B\subseteq \schwmns^{\prime}$ are equicontinuous, then
$\Xi_{n} \longrightarrow 0$ uniformly on
$A\otimes B \subseteq \schwgtm^{\prime}$ \cite{Tr}.
(Here, $A\otimes B$ denotes the set
$\{a\otimes b |a\in A, \quad b\in B\}$.)
This is true by the following lemma.
\proclaim {Lemma 6.37} If
 $A\subseteq \schw^{\prime}$
and $B\subseteq \schwmns^{\prime}$ are equicontinuous, then
$A\otimes B \subseteq \schwgtm^{\prime} $ is an equicontinuous
subset of $\schwgtm^{\prime}$.\endproclaim
\demo{Proof}
Note that $\schwgtm$
is topologized by the seminorms
$$ \pa \w^{p}\s^{q}X^{\ga}{\tilde X}^{\be}
\Xi \pa_{\infty} =
\sup_{x\in G,  m\in M}
\w^{p}(x) \s^{q}(m) |X^{\ga} {\tilde X}^{\be}
\Xi (g,m)| \tag 6.38 $$
where $X^\ga$ acts on the first argument of $\Xi $, and
${\tilde X}^{\be}$ on the second.
By equicontinuity, let $p, q$ and $C$, $D$ be such that
$$|S(\phi )|\leq C\max_{|\ga|\leq p} \pa \w^{p}X^{\ga}\phi
\pa_{\infty}, \qquad \phi \in \schw, \tag 6.39 $$
and
$$|T(\psi )|\leq D\max_{|\be|\leq q} \pa
\s^{q}{\tilde X}^{\be}\psi
\pa_{\infty}, \qquad \psi \in \schwmns, \tag 6.40 $$
for all $S\in A$ and $T\in B$.   We
follow the argument of \cite{Hor, Lemma 3, P. 371}
to show that
$$|(S\otimes T)(\Xi )| \leq CD \max_{|\ga|\leq p, |\be|\leq q}
 \pa \w^{p}\s^{q}X^{\ga}{\tilde X}^{\be}
\Xi \pa_{\infty}, \qquad \Xi \in \schwgtm, \tag 6.41	$$ from
which the lemma follows.
\par
We introduce the integral notation
$$S(\phi ) = \int_{G} S(x) \phi (x) dx $$
$$ T(\psi ) = \int_{M} T(m) \psi (m) dm.$$
Then if $\Xi \in \schwgtm$,
$$|\int_{G}S(x) {\tilde X}^{\be}
\Xi (x, m) dx | \leq C \max_{|\ga|\leq p} \sup_{x\in G}
\w^{p}(x)|X^{\ga} {\tilde X}^{\be} \Xi (x, m)| \tag 6.42 $$
by (6.39). So we have
$$ \aligned
|(S\otimes T) \Xi | & = |\int_{M} T(m)\biggl( \int_{G}
S(x) \Xi (x, m)dx \biggr) dm | \\
& \leq D
 \max_{|\be|\leq q} \sup_{m\in M}
\s^{q}(m)|{\tilde X}^{\be}\int_{G} S(x) \Xi (x, m)dx|
\qquad (6.40)\\ & = D
 \max_{|\be|\leq q} \sup_{m\in M}
|\int_{G}\s^{q}(m) S(x) {\tilde X}^{\be}\Xi (x, m)dx|
\\ &\leq CD
 \max_{|\ga|\leq p, |\be|\leq q}
\pa \w^{p}\s^{q} X^{\ga}{\tilde X}^{\be} \Xi \pa_{\infty}
\endaligned \tag 6.43 $$
by (6.42).
This proves the lemma, and the isomorphism (6.36).
\qed
\enddemo
%
%THE APPENDIX BEGINS HERE
\heading  Appendix. Sets of $C^{\infty}$-vectors.  \endheading
\par
We show that the set of $C^{\infty}$-vectors $\Einf$ for the action
of a Lie group $G$ on a Fr\'echet space $E$ is a Fr\'echet space
for the natural topology, and that
$G$ leaves $\Einf$ invariant and acts differentiably on $\Einf$
by continuous automorphisms. If $E$ is an $m$-convex Fr\'echet *-algebra,
we show that $\Einf$ is also.
These facts are well known \cite{DM},
but the literature seems to prove them only for Banach
or Hilbert spaces \cite{Go},
\cite{Po}, \cite{Ta}. We conclude the appendix with a
technical lemma which is useful in Proposition 2.2.8,
Proposition 6.34 and Theorem 6.35.
\par
Let $\alpha$ denote the action of $G$ on $E$.
We assume that $\alpha$ gives a strongly continuous
action of $G$ on $E$ by continuous automorphisms.
We define the set of $C^{\infty}$-vectors $\Einf$
for the action of $G$  on $E$, to be the set of $e\in E$ such
that all the derivatives  $X^{\ga}e$  (see (1.2.1) with
$\be = \alpha$)
exist for the topology of
$E$.  Here $q$ is the dimension of the Lie algebra $\frak G$
of $G$ and  $\ga \in \N^{q}$.
Let $\normm$ be a family of seminorms topologizing $E$.
Then we topologize $\Einf$ by the seminorms
$$\pa e \pa_{l, m} = \max_{|\ga|\leq l}\pa
X^{\ga}e \pa_{m} \tag A.1   $$
\proclaim {Theorem A.2}
For the topology given by the seminorms
$\pa \quad \pa_{l, m}$,
the locally convex space $\Einf$ is complete.
Thus $\Einf$ is a Fr\'echet space.  In the topology of
$E$, the set $\Einf$ is dense.  The
action $\alpha$ of $G$ on $E$ leaves $\Einf$ invariant,
and each $\alpha_{g}$ restricts to a continuous automorphism
of $\Einf$.  For fixed $e\in \Einf$, the map $g\mapsto \alpha_{g}(e)$
is infinitely
differentiable from $G$ to $\Einf$ (or equivalently, all
the derivatives of elements of $\Einf$ converge in the
topology of $\Einf$).
If $E$ is an ($m$-convex) Fr\'echet [*]-algebra,
and $G$ acts by [*]-automorphisms on $E$, then $\Einf$ is
an ($m$-convex) Fr\'echet [*]-algebra.
\endproclaim
\demo{Proof}
Let $e_{m}$ be a Cauchy sequence in $\Einf$, and for each
$\ga \in \N^{q}$ let $e^{\ga} \in E$ be such that
$$  X^{\ga}e_{m} \longrightarrow e^{\ga}   \qquad {\text{ in $E$. }}
\tag A.3$$
If $\ga = (0, \dots ,0)$, let $e$ denote $e^{\ga}$.
We show $e_{m}\longrightarrow e$ in $\Einf$.  Let $\ga_{i} =
(0, \dots 1,0, \dots )$, with a 1 in the $i$th spot.  Then
$$\aligned  \int_{0}^{t} \alpha_{exptX_{i}}(e^{\ga_{i}})dt
& = \lim_{m} \int_{0}^{t} \alpha_{exptX_{i}}(X^{\ga_{i}}e_{m})dt
\qquad {\text { (A.3) and unif conv}}\\
& = \lim_{m}  \alpha_{exptX_{i}}(e_{m}) - e_{m}
\qquad {\text { Fund thm calc}}\\
& =
\alpha_{exptX_{i}}(e) - e\endaligned \tag A.4 $$
so the derivative $X^{\ga_{i}}e$ exists in the topology of $E$ and
equals  $e^{\ga_{i}}$.
  Repeated application of this argument shows that
$X^{\ga}e$ exists in $E$ for all $\ga \in \N^{q}$, and equals
 $e^{\ga}$.  By (A.3), this implies $e_{m}\longrightarrow e$
in the topology of
$\Einf$.  Thus $\Einf $ is complete, and a Fr\'echet space.
We remark that the completeness is shown in \cite{Go, Cor 1.1}
for $E$ a Hilbert space.
\par
The argument of \cite{Ta, p. 11} shows that if $G$
acts strongly continuously on $E$, then $\Einf$
is dense in $E$.  The simple argument indicated in
\cite{Po, \S 1} shows that $\Einf$ is $G$-invariant.
Also, the proof of \cite{Po, Prop 1.2} works for Fr\'echet
spaces, so $g\mapsto \alpha_{g}(e)$ is differentiable.
\par
Now assume that $E$ is a  Fr\'echet algebra, and that
$\alpha_{g}$ is an  algebra automorphism of $E$
for each $g \in G$.
Let $\ga \in \N^{q}$ and $e, f\in E$.
By the product rule,
$$ X^{\ga}(ef) = \sum_{\be, {\tilde \be}} c_{\be, {\tilde \be}}
(X^{\be}e)(X^{\tilde \be}f),\tag A.5
$$
where $\be$ and $\tilde \be$ have order less than or equal to
$|\ga|$, and $c_{\be, {\tilde \be}}$ are
appropriate constants. We have
$$
\aligned
\pa ef \pa_{l,p} &
= \max_{|\ga| \leq l} \pa X^{\ga} (ef)\pa_{p}
\quad {\text {for the seminorms (A.1)}}\\
& \leq \sum_{\be, {\tilde \be}}
c_{\be, {\tilde \be}}
\pa (X^{\be}e)(X^{\tilde \be}f) \pa_{p}, \quad
{\text{by (A.5)}}
\\ & \leq \sum_{\be, {\tilde \be}}
{\tilde c}_{\be, {\tilde \be}}
\pa (X^{\be}e)\pa_{d}
\pa (X^{\tilde \be}f) \pa_{d}
\qquad {\text {$E$ Fr\'echet alg}}
\\
& \leq C\pa e \pa_{l,d }
\pa f \pa_{l,d }, \quad {\text{by def of norms (A.1)}}
\endaligned \tag A.6
$$
for some $C>0$ depending only on $l$ and $d$. Thus
$\Einf$ is a Fr\'echet algebra.
If $E$ is $m$-convex, then we may take $d=p$ in (A.6), so
the $m$-convexity of $\Einf$
clearly follows (without using Theorem 3.1.4).
If $E$ is a Fr\'echet *-algebra and
$G$ acts by
*-automorphisms on $E$,
then $X^{\ga}e^{*} = (X^{\ga}e)^{*}$ for all $e\in \N$.
Hence the * operation will be continuous on the seminorms (A.1)
and  $\Einf$ is  a Fr\'echet *-algebra.
This proves Theorem A.2.  \qed \enddemo
\par
Let $F$ be any Fr\'echet space.
Then $G$ has a natural continuous
action on the completed projective tensor product
$E\hat \otimes_{\pi} F$
given by
$$\alpha_{g}(e\otimes f)  = \alpha_{g} (e) \otimes f\tag A.7  $$
on elementary tensors.
It is easily checked using the absolutely convergent
series expression for elements of $E\tensp F$ in
\cite{Tr, Thm 45.1} that this gives a strongly continuous
action of $G$.
The surjectivity in the following theorem
uses the theorem of Dixmier and
Malliavin \cite{DM, Thm 3.3}.
\proclaim {Theorem A.8} Let $C_{\alpha}^{\infty}$
denote the set of $C^{\infty}$-vectors for the action of
$\alpha$ on $E\hat \otimes_{\pi} F$.  Then
there is a continuous surjection
$\Einf \hat \otimes_{\pi} F\longrightarrow
C_{\alpha}^{\infty}$
 of Fr\'echet spaces.
\endproclaim
\demo {Proof} Since we shall always be dealing with the projective
tensor product in this theorem, we omit the subscript $\pi$ from
$\hat \otimes_{\pi}$.
We have a natural continuous map $\pi \colon \EinF \longrightarrow E\hat
\otimes F$.  We show that the image of $\pi$ is contained in $\Cai$.
If $x\in \EinF$, we may write
$x$ as an absolutely convergent series $\sum \lambda_{n}
e_{n}\otimes f_{n}$, where $e_{n}\longrightarrow 0$
in $\Einf$, $f_{n}\longrightarrow 0$ in $F$, and
$\sum |\lambda_{n}|<1$ \cite{Tr,  Thm 45.1}.
Let $X$ be in the Lie algebra of
$G$.  If $z \in \Cai$, then $Xz$ is the limit
$$ Xz = \lim_{t\rightarrow 0}
{\alpha_{(exptX)}(z) - z \over {t}} \tag A.9$$
We show that if we plug  the series expression for $x$
in place of $z$ into  (A.9), then the limit converges to
an element of $E \hat \otimes F$
in the topology of $E \hat \otimes F$.
This will imply that $x \in \Cai$.
Define $$ \psi(t)(e) = \alpha_{(exptX)}(e) $$
We abbreviate $\psi (t)(e_{n})$ by $\psi_{n} (t)$.
Since $e_{n} \longrightarrow 0$ in $\Einf$,
 $\psi^{\prime}_{n}(0)$ converges
to zero in $\Einf$.  Define
$$ x^{\prime} \equiv \sum_{n=0}^{\infty} \lambda_{n}\psi^{\prime}_{n}
(0) \otimes f_{n} \tag A.10 $$
The series clearly converges absolutely in $E\hat \otimes F$ and so
$x^{\prime} \in E \hat \otimes F$.
We proceed to show that the limit (A.9) with $z=x$
converges in $E \hat \otimes F$ to $x^{\prime}$.
By elementary calculus
$$ \psi_{n} (t) = \psi_{n}(0) + t\psi^{\prime}_{n}(0) +
t^{2}
\int_{0}^{1} (1-s)\psi_{n}^{\prime \prime}(ts) ds. \tag A.11$$
Subtracting (A.10) from (A.9) with the series for $x$ plugged in
for $z$, we see that $Xx-x^{\prime}$ is  the
limit (if it exists), as $t$ tends
to zero, of the infinite series
$$\split \sum_{n=0}^{\infty} \lambda_{n}
\biggl({\psi_{n} (t) -
\psi_{n}(0) \over {t}} - \psi^{\prime}_{n}(0)\biggr)
\otimes f_{n} = \\   \sum_{n=0}^{\infty} \lambda_{n} \biggl(t
\int_{0}^{1} (1-s)\psi_{n}^{\prime \prime}(ts) ds\biggr)
\otimes f_{n}. \endsplit\quad {\text {by (A.11)}} \tag A.12 $$
Let $\pa \quad \pa_{d}$, $\pa \quad \pa_{l}$
be norms for the topology
on $E$, $F$ respectively and let $\pa \quad \pa_{d,l}$
be the induced norm on the tensor product.
Let $T_{t} \colon \Einf \longrightarrow E$ be the continuous
linear operator $T_{t}(e) =
\psi^{\prime \prime} (t) (e) - X^{2}e$.  Then
$T_{t}$ converges to zero pointwise as $t\longrightarrow 0$.
By the uniform boundedness principle for Fr\'echet spaces\cite{RS, Thm
V.7},
there is some $k, p\in \N$, $C>0$ such
that
$$ \pa T_{t}e \pa_{d} \leq C\max_{|\ga|\leq p}
 \pa X^{\ga}e \pa_{k}
\tag A.13
$$
for $t \in [0, 1]$, $e\in \Einf$.  It follows that
$$
\pa \psi_{n}^{\prime \prime}
(t)  \pa_{d} \leq C\max_{|\ga|\leq p}
 \pa X^{\ga}e_{n} \pa_{k}
+ \pa X^{2} e_{n} \pa_{d}.$$
for all $t\in [0, 1]$ and $n\in \N$.
Since the terms on the right hand side tend to zero as
$n\longrightarrow \infty$, there is some $D>0$
such that
$$\sup_{t\in [0, 1], n\in \N} \pa
 \psi_{n}^{\prime \prime}(t)\pa_{d}<D.
\tag A.14$$
The right hand side of (A.12) evaluated at $\pa \quad\pa_{d,l}$
is then less than or equal to
$$|t|\sum_{n=0}^{\infty} \lambda_{n}D  \pa f_n \pa_{l}
\tag A.15 $$
if $|t|<1$.
This clearly tends to zero as $t$ tends to zero.  Thus the limit
$Xx$ converges to $x^{\prime} $ in $E\hat \otimes F$.  It follows that
$\pi$ maps
$\EinF $ continuously into $\Cai$.
\par
 We now prove that  $C^{\infty}_{\alpha}
\subset \pi (\EinF)$.
By  \cite {DM, Thm 3.3} any element of $C^{\infty}_{\alpha}$
is a finite sum of elements of the form $\alpha_{f}(y)$, where
$f\in C^{\infty}_{c}(G)$ and $y\in \EF$.  By definition, the
expression $\alpha_{f}(y)$
is the integral
$$\int_{G} f(g) \alpha_{g} (y) dg.\tag A.16 $$
If we write $y$ as an absolutely convergent series of
elementary tensors (as in \cite {Tr, Thm 45.1} )
and take the integral inside the sum, we get a series
converging absolutely in $\EinF$.  Hence
the series for
 $\alpha_{f}(y)$ converges in $ \EinF$.  Since
 $\pi$ maps this series to $\alpha_{f}(y)$ in
$E\hat \otimes F$,
 we have proved that $\Cai \subseteq \pi (\EinF)$.
This proves Theorem A.8.
\qed \enddemo
%
%THE APPENDIX ENDS HERE
%

\vfill\eject
%
%BEGIN INDEX
%
\hoffset=0.1truein
\voffset=0.4truein
\hsize=6.5truein
\vsize=8.5truein
\addto\tenpoint{\normalbaselineskip=14pt\normalbaselines}
\addto\eightpoint{\normalbaselineskip=11pt\normalbaselines}
\newdimen \fullhsize
\fullhsize=6.5in \hsize=3.2in
\def\fullline{\hbox to \fullhsize}
\def\makeheadline{\vbox to 0pt{\vskip-22.5pt
                  \fullline{\vbox to8.5pt{}\the\headline}\vss}
      \nointerlineskip}
\def\makefootline{\baselineskip=24pt \fullline{}}
\def\centerline#1{\fullline{\hss#1\hss}}
\let\lr=L \newbox\leftcolumn
\output={\if L\lr
      \global\setbox\leftcolumn=\columnbox \global\let\lr=R
   \else \doubleformat  \global\let\lr=L\fi
   \ifnum\outputpenalty>-20000 \else\dosupereject\fi}
\def\doubleformat{\shipout\vbox{\makeheadline
    \fullline{\box\leftcolumn\hfil\columnbox}
    \makefootline}
   \advancepageno}
\def\columnbox{\leftline{\pagebody}}
\document
\line{\hfill}
\centerline{\smc Index}
\line{\hfill}
\line{\hfill}
\line{action (of $G$ on a Fr\'echet algebra) \ 2.2.4. \hfill}
\line{$ax+b$ group \ 1.3.5, 1.3.11, 1.6.1, 5.20.\hfill}
\line{approximate units \ 6.12.\hfill}
\line{Arens, R. \ \ \cite{Ar}, \S 0, end \S 3.2. \hfill}
\vskip\baselineskip
\line{Bagley, R. W. \ \ \cite{BWY}, Pf. of 1.5.11. \hfill}
\line{Barker, W.H. \ \ \cite{Bar}, \S 0, Intro. to \S 6, 6.18.\hfill}
\line{Bost, J. B. \ \ \cite{Bo}, \S 0, 5.27. \hfill}
\line{Bourbaki, N. \ \ \cite{Bou}, Pf. of 1.5.5, Pf. of 6.13.\hfill}
\line{bounds $Ad$\ \S 1.3, 1.3.10-16, 1.4.3, 1.4.8, 1.5.12,\hfill}
\line{\ \ \ \ 2.2.1, 2.2.6, 3.1.7, 3.2.6, 4.11, 6.12-13, 6.18.\hfill}
\line{bounds $Ad$ on $H$ \ 4.3, 4.5-6, 4.9, 5.12, 5,26, 6.13.\hfill}
\vskip\baselineskip
\line{$C^{\infty}$-vectors \ 1.2.1, 1.2.21, 2.1.5, 4.6, 4.9, \hfill}
\line{\ \ \ \ 5.27, 6.5, appendix.\hfill}
\line{$C_{c}^{\infty}(G, E)$ (compactly supp.
$C^{\infty}$ functions) \hfill}
\line{\ \ \ \ Pf. of 2.1.5, 2.1.7, Pf. of 6.35. \hfill}
\line{closed subgroups \ 1.5.1, 1.5.10, 5.14, 5.23, 6.13. \hfill}
\line{compactly generated group \ 1.1.20, 1.1.21, 1.4.3, \hfill}
\line{\ \ \ \ 1.4.8, 1.5.1, 1.5.5, 1.5.10, 1.5.13, 3.1.30, 3.2.6, \hfill}
\line{\ \ \ \ 3.2.19, 6.13 \hfill}
\line{compactness \ 1.1.21, 1.2.11, 1.3.16.\hfill}
\line{Connes, A. \ \ \cite{BC}, \S 0, \cite{Co}, \S 0. \hfill}
\line{continuous action \ 2.2.4. \hfill}
\line{Corwin, L. \ \ \cite{Co-Gr}, Pf. of 2.1.5. \hfill}
\vskip\baselineskip
\line{differential operator from Lie algebra \ ($X^{\ga}$)\ 1.2.1.\hfill}
\line{discrete group \ 1.1.17, 1.1.21, 1.5.10, 3.1.26, 3.2.7,\hfill}
\line{\ \ \ \ 4.5, 6.16, 6.17.\hfill}
\line{Dixmier, J. \ \ \cite{DM}, Pf. of 2.1.5, Pf. of 2.1.7, \hfill}
\line{\ \ \ \ Pf. of 6.8, Intro. to App., Pf. of A.8. \hfill}
\line{dominates \ 1.1.8, 1.1.12. \hfill}
\vfill\eject
\line{\hfill}
\line{\hfill}
\line{\hfill}
\line{\hfill}
\line{Du Cloux, F. \ \ \cite{DuC 1-4}, \S 0, 1.3.14, 1.6.1, \hfill}
\line{\ \ \ \ 2.2.4.\hfill}
\line{Dzinotyiweyi, H.A.M. \ \ \cite{Dz 1-2}, 1.1.1, Pf. of  \hfill}
\line{\ \ \ \ 1.1.21, Pf. of 1.2.11, Pf. of 1.5.5. \hfill}
\vskip\baselineskip
\line{$E$ \ (Fr\'echet space) \ 2.1.0, \S 2.1, Appendix.\hfill}
\line{Effros, E. \ \ \cite{EH}, \S 0. \hfill}
\line{Elliott, G. \ \ \cite{ENN}, \S 0, 2.2.4. \hfill}
\line{equivalent scales\ 1.1.8, 1.1.12\hfill}
\line{exponentiated word weight\ 1.1.17, \S 1.6.\hfill}
\vskip\baselineskip
\line{Fr\'echet algebra \ \S 1.3, 1.3.2, 1.3.4, 1.3.13,  \hfill}
\line{\ \ \ \ \S 2.2, 2.2.6, 5.17. \hfill}
\line{Fr\'echet *-algebra \ \S 1.3, 1.3.2, 1.3.5, 1.3.13, \hfill}
\line{\ \ \ \ \S 4, 4.6, 4.9, 5.17. \hfill}
\vskip\baselineskip
\line{$\sona$ (smooth crossed product) \ after 2.2.6. \hfill}
\line{$G \rtimes^{\w} \schwm$ \ 5.15-5.17. \hfill}
\line{$GL(n, \R)$ \ 1.6.10, 5.23, 6.13, 6.15.\hfill}
\line{$G_{0}$ \ (connected comp. of iden. of $G$) \ 1.1.21,\hfill}
\line{\ \ \ \ Intro. to 1.2.1, Pf. of 1.5.10, 1.5.13, 6.12.\hfill}
\line{group Schwartz algebra, $\Cal S(G)$ \ \S 1.3, 1.3.13, \hfill}
\line{\ \ \ \ 1.5.14, \S 3.2, 5.11.\hfill}
\line{guages \  1.1.1-4, 1.1.10, 1.1.14, 1.1.17, 1.5.15.\hfill}
\line{gauge that bounds $Ad$ \ 1.3.9, \S 1.4, 1.4.3,  \hfill}
\line{\ \ \ \ 1.4.7-9, 1.5.12, 1.5.13, 1.5.15.\hfill}
\line{Goodman, R. \ \ \cite{Go}, Intro. to App., Pf. \hfill}
\line{\ \ \ \ of A.2.\hfill}
\line{Greenleaf, F.P. \ \ \cite{Co-Gr}, Pf. of 2.1.5. \hfill}
\vskip\baselineskip
\line{Herb, R.\ \ \cite{HW}, 1.6.4, 6.18.\hfill}
\line{Heisenberg group \ 1.5.15, 3.2.7, 4.5.\hfill}
\line{Hewitt, E. \ \ \cite{HR}, 1.1.20, Pf. of 1.1.21.\hfill}
\line{Horvath. J. \ \ \cite{Hor}, 6.37.\hfill}
\line{Howe, R. E. \ \ \cite{Ho},  \S 0, 1.5.14. \hfill}
\line{ideals \ 2.2.8.\hfill}
\line{inverse (see projective) limit of Fr\'echet algebras \hfill}
\line{\ \ \ \ Intro. to 3.1.7, 3.1.7.\hfill}
\vskip\baselineskip
\line{Jenkins, L.W. \ \ \cite{Je}, 1.4.2, Pf. of 6.13, 6.18. \hfill}
\line{Ji, R. \ \ \cite{Ji}, \S 0, 1.1.1, Intro. to \S 1.5, \S 6 and\hfill}
\line{\ \ \ \  6.13, 6.16.\hfill}
\line{Jolissaint, P. \ \ \cite{Jo}, \S 0, 1.1.1, Intro. to \S 6, \hfill}
\line{\ \ \ \ Pf. of 6.8, 6.16, 6.17, 6.33. \hfill}
\vskip\baselineskip
\line{Kargapolov, M. I. \ \ \cite{KM}, Pf. of 1.4.3. \hfill}
\vskip\baselineskip
\line{$L^{1}(G)$, $L^{r}(G)$\  ($|\varphi(g)|$,
$|\varphi(g)|^{r}$ integrable against \hfill}
\line{\ \ \ \ left Haar measure), but see $L^{k}(G)$ 3.2.6. \hfill}
\line{$L_{1}^{\s}(G)$, $L^{\s}(G)$ \
($\s$-rapidly vanishing $L^{1}$-\hfill}
\line{\ \ \ \ functions) \ 1.2.19. \hfill}
\line{$L_{r}^{\s}(G)$ \ ($\s$-rapidly vanishing $L^{r}$ functions)\
6.3. \hfill}
\line{$L_{1}^{\s}(G, A)$, $L^{\s}(G, A)$\
($\s$-rapidly vanishing $L^{1}$-
\hfill}
\line{\ \ \ \  functions) \ 2.1.3, 2.2.6, 3.1.7. \hfill}
\line{largest gauge, weight \ 1.1.17, 1.1.21.\hfill}
\line{length function \ \S 0, Intro. to \S 1.1, 1.1.1, 1.1.17,\hfill}
\line{\ \ \ \ 6.17. \hfill}
\line{Lie group \  1.2.\hfill}
\line{locally compact group \ 1.1.20, 1.1.21, 1.5.1, 1.5.3, \hfill}
\line{\ \ \ \ 1.5.4, 2.1.0, 2.2.6, 2.2.8, 3.1.7, 3.1.18, 3.1.30,\hfill}
\line{\ \ \ \ 3.2.2, 6.3, 6.13.\hfill}
\line{locally compact space \ \S 5. \hfill}
\line{$L^{r}$-Schwartz functions,
$\Cal S_{r}^{\s}(G)$ \ \S 6, 6.1.\hfill}
\line{Losert, V. \ \ \cite{Lo}, 1.5.10, Intro. and Pf. of 1.5.13. \hfill}
\line{Ludwig, J. \ \ \cite{Lu 1-2}, \S 0, 1.4.4, 1.5.14, 1.5.15.\hfill}
\vskip\baselineskip
\line{$m$-convex Fr\'echet algebra \ \S 3, 3.1.7.\hfill}
\line{$m$-$\s$-tempered\ 3.1.1, 3.1.7, 3.1.18, 3.1.25, 3.1.26,  \hfill}
\line{\ \ \ \ 4.6, 4.9, 5.17. \hfill}
\line{$m$-sub-polynomial scale \ 3.1.6, 3.1.7, \S 3.2, 3.2.8,  \hfill}
\line{\ \ \ \ 3.2.12, 5.17. \hfill}
\line{Mallios, A. \ \ \cite{Ma}, \S 0.\hfill}
\line{Micheal, E. \ \ \cite{Mi}, \S 0, Intro. to 3.1.7, 6.16,\hfill}
\line{\ \ \ \ 6.18. \hfill}
\vskip\baselineskip
\line{$\N$,$\N^{+}$ \ natural numbers with, without zero.\hfill}
\line{Nest, R. \ \ \cite{Ns}, \S 0, \cite{ENN}, \S 0, 2.2.4.\hfill}
\line{nilpotent group \ \S 1.4, \S 1.5, 1.5.10, 1.5.14, \hfill}
\line{\ \ \ \ 6.15.\hfill}
\line{non $m$-convex group Schwartz algebras  \hfill}
\line{\ \ \ \ \S 3.2, 3.2.6, 3.2.8, 3.2.12, 3.2.13. \hfill}
\line{$\pa \quad \pa_{1}$ \ \ ($L^{1}$ norm) \ 1.2.3. \hfill}
\line{$\pa  \quad \pa_{r}$ ($L^{r}$ norm) \ 6.1, 6.3.  \hfill}
\line{$\pa \quad \pa_{m}$ or $\pa \quad \pa_{d}$  (seminorms on
$L_{1}^{\s}(G)$, $E$,  $A$  \hfill}
\line{\ \ \ \ or $C^{\s}(M)$) 1.2.20, 2.1.0, 2.2.4, 3.1.1, 5.2. \hfill}
\line{$\pa \quad \pa_{m, \ga}$\ \  (seminorm on $\Cal S(M)$ or $\Cal S(G)$)
\hfill}
\line{\ \ \ \ 1.2.3, 5.10. \hfill}
\line{$\pa \quad \pa_{l, m}$ (usually
seminorms on $C^{\infty}$-vectors)\hfill}
\line{\ \ \ \ 4.2, A.1, also 3.1.8. \hfill}
\line{$\pa \quad \pa_{m, \ga, d}$ \ (seminorm on $\sonec$) \ 2.1.1,\hfill}
\line{\ \ \ \ also 2.2.7, 3.1.8,  4.10, 5.15. \hfill}
\line{nuclear Fr\'echet space \ \S 6, 6.24, 6.27, 6.33-5.\hfill}
\vskip\baselineskip
\line{Paterson, A. L. T. \ \ \cite{Pa}, 1.4.2, Pf. of 1.4.6, \hfill}
\line{\ \ \ \  3 and 9, Intro. to \S 1.5, 1.5.4, Intro. to \hfill}
\line{\ \ \ \ 1.5.5, 1.5.14, 6.15.\hfill}
\line{Peimbert, H.A. \ \ \cite{Pe}, Intro.  3.1.4, end \S 3.2. \hfill}
\line{Phillips, N.C. \ \ \cite{Ph}, \S 0, Intro. to 3.1.7.\hfill}
\line{Pietsch, A. \ \ \cite{Pi}, Pf. of 6.24, 6.27. \hfill}
\line{polynomial growth groups\ \S 1.5, 3.2.7, 4.5.\hfill}
\line{Poulsen, N. \ \ \cite{Po}, Intro. to App., Pf. of A.2.\hfill}
\line{projective limit of Fr\'echet/Banach algebras \hfill}
\line{\ \ \ \ \S 0, 2.1.3, Intro. to 3.1.7, 3.1.7.\hfill}
\line{projective tensor product (of Fr\'echet spaces)  \hfill}
\line{\ \ \ \ \S 0, Pf. of 2.2.8, 5.20, 6.34, 6.35, Intro. \hfill}
\line{\ \ \ \ to A.8, A.8.\hfill}
\line{Pytlik, T. \ \ \cite{Py}, 1.1.1, 1.1.17. \hfill}
\vskip\baselineskip
\line{$\Q$, \ rational numbers, 3.1.26. \hfill}
\line{questions \ 1.2.4, 1.3.4, 1.4.7, 1.5.3, 3.1.3,  \hfill}
\line{\ \ \ \ 3.1.25, 3.1.30, 3.2.12, 6.33.\hfill}
\line{quotients of algebras \ 2.2.8.\hfill}
\line{quotients of groups \ 1.1.21, 1.4, 1.5.5.\hfill}
\vskip\baselineskip
\line{Rader's algebra \ 6.18.\hfill}
\line{Reed, M. \ \ \cite{RS}, Pf. of A.8. \hfill}
\line{Rolewicz, S. \ \ \cite{RZ}, end \S 3.2, \cite{MRZ}, Intro.\hfill}
\line{\ \ \ \  to 3.1.4.\hfill}
\vskip\baselineskip
\line{$SL(n, \R)$ \ 1.6.4, 6.18.\hfill}
\line{scale, $\s$ \ 1.1.1, \S 5.\hfill}
\line{scaled $(G, \w)$-space \ 5.13, 5.17.\hfill}
\line{Schwartz, L. \ \ \cite{Schw}, 1.2.19, 1.3.3, 2.1.0, \hfill}
\line{\ \ \ \ 2.2.8, 6.3, 6.34. \hfill}
\line{Schwartz functions \ \hfill}
\line{\ \ \ \ Fr\'echet space valued, $\Cal S(G, E)$ \ 2.1.1, 2.1.5.\hfill}
\line{\ \ \ \ $L^{1}$-norm, $\Cal S_{1}^{\s}(G)$ \ \S 1.2, 1.2.2.\hfill}
\line{\ \ \ \ $L^{r}$-norm, $\Cal S_{r}^{\s}(G)$ \  \S 6, 6.1.\hfill}
\line{\ \ \ \ sup norm, $\Cal S_{\infty}^{\s}(G)$ \ \S 5.\hfill}
\line{\ \ \ \ $\schwr$, standard Schwartz functions on $\R$. \hfill}
\line{\ \ \ \ $\Cal S^{k}(G)$ \   3.2.6.\hfill}
\line{\ \ \ \ $G \rtimes^{\w} \schwm$ \ 5.15-5.17. \hfill}
\line{\ \ \ \ $\Cal S^{\w \times \s}(G \times M)$ \ 6.35. \hfill}
\line{Schweitzer, L. B. \ \ \cite{Sc 1-4}, \S 0, Intro. to 1.3.9 \hfill}
\line{\ \ \ \  and \S 1.4, Pf. of 2.2.6, 3.1.26. \hfill}
\line{$\s^{m}(g)$\  ($\s(g)$ to the $m$th power) \ 1.2.3, 1.2.20, \hfill}
\line{\ \ \ \ 2.1.1-2, 5.13, 5.15. \hfill}
\line{$\s_{-}$ (inverse of $\s$) \ 1.3, 1.3.2.\hfill}
\line{$\s$-tempered action \ 2.2.4.\hfill}
\line{$\s_{k}$ \ 3.2.1.\hfill}
\line{smooth action (means differentiable and tempered \hfill}
\line{\ \ \ \ in \cite{ENN}), \S 0, 2.2.4. \hfill}
\line{smooth crossed product, $\sona$ \ \S 2.3, 2.2.6.\hfill}
\line{*-isometric seminorms \ \S 4.\hfill}
\line{strongly dominates \ 1.1.12.\hfill}
\line{strongly equivalent gauges \ 1.1.12, 1.1.14.\hfill}
\line{strongly $m$-tempered action \ 3.1.26-7, 3.1.30.\hfill}
\line{sub-polynomial scale \ 1.3.1.\hfill}
\vskip\baselineskip
\line{$\T$, \ circle group.\hfill}
\line{Taylor, M.E. \ \ \cite{Ta}, Intro. to App., Pf. A.2.\hfill}
\line{tempered action\ 2.2.4.\hfill}
\line{topological group \ 1.1.1.\hfill}
\line{transformation group \ \S 5\hfill}
\line{translationally equivalent scale\ 1.2.7, 1.2.8, \hfill}
\line{\ \ \ \ 1.2.11.\hfill}
\line{Tr\'eves, F. \ \ \cite{Tr}, 2.1.10, 2.2.8, Intro. and Pf. \hfill}
\line{\ \ \ \ of 6.24, Pf. of 6.35, Intro. and Pf. of A.8.\hfill}
\line{triangular matrix group \ 1.4.4, 1.5.15, 3.2.7.\hfill}
\line{Type R Lie group \ \S 1.4, \S 1.5.\hfill}
\vskip\baselineskip
\line{uniformly translationally equivalent scale \hfill}
\line{\ \ \ \ 1.2.10, 1.2.11, 5.7, 5.9, 5.12, \S 5. \hfill}
\vskip\baselineskip
\line{Waelbroeck, L. \ \ \cite{Wa}, \S 1.3.\hfill}
\line{Warner, G. \ \ \cite{War}, \S 0, 4.11, Intro. \hfill}
\line{\ \ \ \ to \S 6, 6.18.\hfill}
\line{weight\  \S 1.1, 1.1.5-7, 1.2.11.\hfill}
\line{word gauge \ 1.1.17, 1.1.21, 1.5.15, 1.6.4.\hfill}
\vskip\baselineskip
\line{$X^{\ga}$ (differential operator) \ 1.2.1. \hfill}
\vskip\baselineskip
\line{$\Z$,  \ integers, 1.1.11, 1.1.22, 5.18, 6.6,7.\hfill}
\line{Zelazko, W. \ \ \cite{Ze}, \S 0, end \S 3.2. \hfill}
\line{zero Schwartz space \ 6.18, 6.23.\hfill}
\vskip\baselineskip
\line{$\scriptstyle \text{Department of Mathematics,
University of Calgary}$\hfill}
\line{$\scriptstyle \text{Calgary, Alberta T2N 1N4 CANADA}$\hfill}
\line{$\scriptstyle \text{\it Current Address:}$ \hfill}
\line{$\scriptstyle \text{Department of Mathematics, University
of Victoria}$ \hfill}
\line{$\scriptstyle \text{Victoria, B.C. V8W 3P4 CANADA}$ \hfill}
\line{$\scriptstyle \text{{\it Email:} lschweit\@ sol.uvic.ca}$\hfill}
\supereject
\if R\lr \null\vfill\eject\fi

\Refs
\widestnumber\key{DuC 1}
%
\ref \key Ar \by R. Arens \paper The space $L^{\omega}$
and convex topological rings \jour Bull. Amer. Math. Soc.
\yr 1946 \pages 931--935 \vol 52
\endref
%
\ref \key BWY \by R.W. Bagley, T. S. Wu, and J.S. Yang \paper
Compactly Generated Subgroups and Open Subgroups of Locally
Compact Groups \jour Proc. Amer. Math. Soc. \vol 103(3) \yr 1988
\pages 969--976 \endref
%
\ref \key Bar \by W.H. Barker \paper $L^{p}$ harmonic
analysis on $SL(2, \R)$ \jour Mem. Amer. Math. Soc., no. 393.\vol 76
\yr 1988 \pages 1--109 \endref
%
\ref \key BC \by P. Baum and A. Connes
\paper Chern Character for Discrete Groups
\inbook Collection: A fete of topology  \publ Academic press
\publaddr Orlando, Florida  \yr 1988 \pages 163--232 \endref
%
\ref \key Bo \by J.B. Bost \paper Principe D'Oka, K-Theorie
et Systems Dynamiques Non-commutative
\jour Invent. Math \vol 101 \pages 261--333
\yr 1990 \endref
%
\ref \key Bou \by N. Bourbaki \book Int\'egration,
no. 1281
\publ Act. Sci. et Ind. \publaddr Paris \yr 1959 \endref
%
\ref \key Br\by R. M. Brooks \paper On Representing F*-algebras
\jour Pacific J. Math \vol 39 \yr 1971 \pages 51--69 \endref
%
\ref \key Ch \by I. M. Chiswell \paper Abstract length
functions in groups \jour Math. Proc. Camb. Phil. Soc. \vol 80
\yr 1976 \pages 451--463 \endref
%
\ref \key Co  \by A. Connes
\paper An Analogue of the Thom Isomorphism for
Crossed Products of a C*-algebra by an Action of $\R$
\jour Adv. in Math.  \vol 39 \pages 31--55
\yr 1981
\endref
%
\ref \key Co-Gr \by L. Corwin, F.P. Greenleaf
\book Representations of Nilpotent Lie Groups and
their Applications, Cambridge studies in advanced
mathematics v. 18 \publ Cambridge Univ Press
\publaddr Cambridge, New York
\yr 1989 \endref
%
\ref \key Da \by H. G. Dales
\paper Automatic Continuity: a survey \jour Bull. London Math. Soc.
\vol 10 \yr 1978 \pages 129--183 \endref
%
\ref  \key  DM \by J. Dixmier and  P. Malliavin
\paper Factorisations de fonctions et de vecteurs id\'efiniment
diff\'erentiables \jour Bull. des sciences Math.
\vol 102 \yr 1978 \pages 305--330 \endref
%
\ref \key DuC 1 \by F. Du Cloux \paper Jets de
Fonctions Diff\'erentiable sur le Dual d'un Groupe
de Lie Nilpotent \jour Invent. Math. \vol 88 \yr 1987
\pages 375--394  \endref
%
\ref \key DuC 2\by F. Du Cloux \paper Repr\'esentations Temp\'er\'ees
des Groupes de Lie Nilpotent \jour J. Funct. Anal.\yr 1989
\vol 85 \pages 420--457 \endref
%
\ref \key DuC 3 \by F. du Cloux \paper Repr\'esentations
de Longueur Finie des Groupes de Lie R\'esoloubles
\jour Mem. Amer. Math. Soc. \vol 80, 407  \yr 1989
\pages 1--78 \endref
%
\ref \key DuC 4 \by F. du Cloux \paper Sur Les Repr\'esentations
Diff\'erentiables des Groupes de Lie Algebriques \jour Preprint
\yr 1989 \endref
%
\ref\key Dz 1\by H.A.M. Dzinotyiweyi \paper Weighted
Function Algebras on Groups and Semigroups \jour Austral. Math. Soc.
\vol 33 \yr 1986 \pages 307--318 \endref
%
\ref \key Dz 2\by H.A.M. Dzinotyiweyi \book The analogue of the
group algebra for topological semigroups \publ Pitman Advanced Publishing
Program  \publaddr Boston, London, Melbourne \yr 1984 \endref
%
\ref \key EH\by E. Effros and F. Hahn
\paper Locally Compact Transformation Groups and
C*-algebras \jour Mem. Amer. Math. Soc. \vol 75
\yr 1967 \endref
%
\ref \key ENN \by G.A. Elliott, T. Natsume,
and R. Nest \paper Cyclic cohomology for one-parameter
smooth crossed products \jour Acta Math.
\vol 160 \yr 1988 \pages 285--305 \endref
%
\ref \key Go\by R. Goodman \paper Analytic and Entire Vectors
for Representations of Lie Groups \jour Trans. Amer. Math. Soc.
\vol143 \yr 1969 \pages 55--76
\endref
%
\ref \key Ha \by N. Harrison \paper Real Length Functions
in Groups \jour Trans. Amer. Math. Soc.
\vol 174 \yr 1972 \pages 77--106 \endref
%
\ref \key HW \by R. Herb and L. Wolf
\paper Rapidly Decreasing Functions on General
Semisimple Lie Groups \jour Comp. Math. \vol 58 \yr 1986
\pages 73--110 \endref
%
\ref \key HR \by E. Hewitt and K.A. Ross \book Abstract Harmonic
Analysis I \publ Springer-Verlag \publaddr Berlin Heidelberg New York
\yr 1979
\endref
%
\ref \key Hor \by J. Horvath
\book Topological Vector Spaces and Distributions
\vol I \publ Addison-Wesley \publaddr Mass
\yr 1966 \endref
%
\ref \key Ho \by R. E. Howe \paper On a Connection Between
Nilpotent Groups and Oscillatory Integrals Associated to
Singularities \jour Pacific J. Math. \yr 1977 \pages 329--363
\vol 73, 2 \endref
%
\ref \key Je \by L. W. Jenkins \paper Growth of Connected Locally
Compact Groups \jour J. Funct. Anal. \vol 12 \pages 113-127
\yr 1973 \endref
%
\ref \key Ji \by R. Ji \paper Smooth Dense Subalgebras
of Reduced Group C*-Algebras, Schwartz Cohomology of
Groups, and Cyclic Cohomology \jour J. Funct. Anal
\vol 107(1)\yr 1992\pages 1--33\endref
%
\ref \key Jo \by P. Jolissaint \paper
Rapidly decreasing functions in reduced
C*-algebras of groups \jour Trans. Amer. Math. Soc.
\vol 317, 1 \yr 1990 \pages 167--196 \endref
%
\ref \key KM \by M. I. Kargapolov and Ju.I. Merzljakov
\book Fundamentals of the Theory of Groups
\publ Springer-Verlag \publaddr New York Heidelberg
Berlin \yr 1979 \endref
%
\ref \key Lo \by V. Losert \paper On the Structure of Groups
with Polynomial Growth \jour Math. Z.\vol 195 \pages 109--117
\yr 1987 \endref
%
\ref \key  Lu 1 \by J. Ludwig \paper Minimal C*-Dense Ideals
and Algebraically Irreducible Representations of the
Schwartz-Algebra of a Nilpotent Lie group \jour Lect.
Notes in Math. \vol 1359 \yr 1987 \pages 209--217 \endref
%
\ref \key  Lu 2 \by J. Ludwig \paper Topologically
irreducible representations of the Schwartz-algebra of
a nilpotent Lie group \jour Arch. Math. \vol 54
\yr 1990 \pages 284--292 \endref
%
\ref \key Ly \by R. C. Lyndon \paper Length Functions in
Groups \jour Math. Scand. \vol 12 \yr 1963 \pages 209--234 \endref
%
\ref \key  Ma \by A. Mallios \book Topological Algebras.
Selected Topics \publ Elsevier Science Publishers B.V.
\publaddr Amsterdam \yr 1986 \endref
%
\ref \key Mi \by E. Micheal \paper Locally multiplicatively
convex topological algebras \jour Mem. Amer. Math. Soc.
\vol 11 \yr 1952 \endref
%
\ref \key MRZ\by B. Mitiagin, S. Rolewicz, and W. Zelazko
\paper Entire functions in $B_{0}$-algebras  \jour Studia Math.
\vol 21 \pages 291--300 \yr 1962 \endref
%
\ref \key Ns \by R. Nest \paper Cyclic cohomology
of crossed products with $\Bbb Z$ \jour J. Funct.
Anal. \vol 80 \yr 1988 \pages 235--283 \endref
%
\ref \key Pa \by A. L. T. Paterson
\book Amenability \publ Math Surveys and Monographs Vol. 29,
 AMS \publaddr Providence, RI
\yr 1988 \endref
%
\ref \key Pe \by H.A. Peimbert \paper Matrix Algebras and
m-Convexity \jour Dem. Math. \vol 17, 3
\pages 711--722 \yr 1984 \endref
%
\ref \key Ph \by N.C. Phillips \paper $K$-thoery
for Fr\'echet Algebras \jour Intl. Jour.  Math. \yr 1991\vol 2(1)
\pages 77--129 \endref
%
\ref \key Pi \by A. Pietsch
\book Nuclear locally convex spaces, Ergebnisse Der Mathematik
und Ihrer Grensgebiete vol. 66
\publ Springer-Verlag
\publaddr New York/Heidleberg/Berlin \yr 1972
\endref
%
\ref \key Po \by N. Poulsen \paper On $C^{\infty}$ vectors
and intertwining bilinear operators for representations
of Lie groups \jour J. Funct. Anal.
\vol 9 \yr 1970 \pages 87--120 \endref
%
\ref \key Pr\by D. Promislow \paper Equivalence Classes
of Length Functions on Groups
\jour Proc. London Math. Soc. \vol 51(3)
\yr 1985 \pages 449-447 \endref
%
\ref \key Py \by T. Pytlik \paper On the spectral radius
of elements in group algebras
\jour Bulletin de l'Acad\'emie Polonaise des Sciences,
astronomiques et physiques \vol 21
\yr 1973 \pages 899-902 \endref
%
\ref \key RS \by M. Reed and B. Simon
\book Methods of Modern Mathematical
Physics \vol 1. Functional Analysis \publ Academic Press, Inc.
 \publaddr Florida \yr 1980 \endref
%
\ref \key RZ \by S. Rolewicz and W. Zelazko \paper
Some Problems Concerning $B_{0}$-algebras
\jour Tensor (New Series) \yr 1963 \vol 13 \pages 269--276 \endref
%
\ref \key Schw \by L. Schwartz
\paper Produits tensoriels topologiques d'espaces
vectoriels topologiques Espaces vectoriels topologiques
nucl\'eaires
\jour Facult\'e des Sciences de Paris
\yr 1953-1954 \endref
%
\ref \key Sc 1 \by L.B. Schweitzer \paper
Representations of Dense Subalgebras of C*-algebras
with Applications to Spectral Invariance \jour U.C. Berkeley thesis
\yr 1991 \endref
%
\ref\key Sc 2 \by L.B.Schweitzer \paper
A short proof that $M_{n}(A)$ is local
if  $A$ is local and Fr\'echet \jour Intl. Jour. Math.
\yr 1992\pages 581--589 \vol 3(4) \endref
%
\ref \key Sc 3 \by L.B. Schweitzer \paper
Spectral Invariance of Dense Subalgebras of Operator Algebras
\jour Intl. Jour. Math.
\yr to appear \endref
%
\ref \key Sc 4 \by L.B. Schweitzer \paper
A Factorization Theorem for Smooth Crossed Products
\jour preprint
\yr 1992 \endref
%
\ref \key Ta \by M.E. Taylor \book  Noncommutative
Harmonic Analysis  \publ  Mathematical Surveys and Monographs
No. 22,  AMS \yr 1986 \endref
%
\ref \key Tr \by F. Tr\'eves \book Topological Vector Spaces,
Distributions, and Kernels
\publ Academic Press \publaddr New York \yr 1967 \endref
%
\ref \key Vi \by Marie-France Vigneras
\paper On Formal dimensions for reductive
$p$-adic groups \jour preprint \yr 1989\endref
%
\ref \key Wa \by L. Waelbroeck \book Topological Vector
Spaces and Algebras \publ Springer-Verlag \publaddr Berlin/
Heidelberg New York \yr 1971 \endref
%
\ref \key War \by G. Warner \book Harmonic Analysis
on Semisimple Lie Groups Vol. I and II \publ Springer-Verlag
\publaddr Berlin and New York \yr 1972 \endref
%
\ref \key Ze \by W. Zelazko \paper On the locally
bounded and $m$-convex topological algebras
\jour Studia Math. \vol T. XIX. \yr 1960 \pages 333--356
\endref
\endRefs
\enddocument